\documentclass{article}[12pt]
\usepackage{amsmath}
\usepackage{amssymb}
\usepackage{extarrows}

\newtheorem{theorem}{Theorem}

\newtheorem{remark}[theorem]{Remark}

\usepackage{amssymb}
\usepackage{graphicx}
\usepackage[dvipsnames,usenames]{color}

\input{epsf.tex}
\usepackage[affil-it]{authblk}

\author[a]{Ioannis Dimitriou \footnote{ idimit@uoi.gr}\footnote{Corresponding author.}}
\author[b]{Dieter Fiems\footnote{Dieter.Fiems@UGent.be}}
\affil[a]{\small Department of Mathematics, 
	University of Ioannina, 
	45110, Ioannina, Greece.}
 \affil[b]{\small Department of Telecommunication and Information Processing, Ghent University,
St-Pietersnieuwstraat 41, 9000 Gent, Belgium.}
\begin{document}
\title{Some reflected autoregressive processes with dependencies}

\maketitle
\begin{abstract}
    Motivated by queueing applications, we study various reflected autoregressive processes with dependencies. Amongst others, we study cases where the interarrival and service times are proportionally dependent with additive and/or subtracting delay, as well as cases cases where interarrival times depends on whether the service duration of the previous arrival exceeds or not a random threshold. Moreover, we study cases where the autoregressive parameter is constant as well as a discrete or a continuous random variable, as well as cases where . More general dependence structures are also discussed. Our primary aim is to investigate a broad class of recursions of autoregressive type for which several independence assumptions are lifted, and for which a detailed exact analysis is provided. We provide expressions for the Laplace transform of the waiting time of a customer in the system in terms of an infinite product of known Laplace transforms. An integer-valued reflected autoregressive process that can be used to model a novel retrial queueing system with orbital searching time to depend on whether the last busy period starts with an empty or a non empty orbit queue, is also discussed. For such a model the probability generating function of the stationary orbit queue length is given as an infinite product of known generating functions. A first attempt towards multidimensional setting is also analyzed. Some additional generalizations with more general dependence structure are also discussed.
    \end{abstract}
    \vspace{2mm}
	
\section{Introduction}
This work focuses on various stochastic recursions of autoregressive type, such as:
\begin{equation}
    W_{n+1}=[V_{n}W_{n}+B_{n}-A_{n}]^{+},\,n=0,1,\ldots,\label{recu1}
\end{equation}
\begin{equation}
    W_{n+1}=\left\{\begin{array}{ll}
          \left[V^{(1)}_{n}W_{n}+B_{n}-A^{(1)}_{n}\right]^{+},&B_{n}\leq T_{n},\\
        \left[V^{(2)}_{n}W_{n}+T_{n}-A^{(2)}_{n}\right]^{+},& B_{n}>T_{n},
    \end{array}\right.\label{rec2}
\end{equation}
\begin{equation}
    W_{n+1}=\left\{\begin{array}{lll}
    \left[W_{n}+B_{n}-A^{(1)}_{n}\right]^{+},&\text{ with probability (w.p.) }p,\\
          \left\{\begin{array}{ll}
              \left[V^{(1)}_{n}W_{n}+\widehat{B}_{n}-A^{(1)}_{n}\right]^{+},&\widehat{B}_{n}\leq T_{n},\\
        \left[V^{(2)}_{n}W_{n}+T_{n}-A^{(2)}_{n}\right]^{+},& \widehat{B}_{n}>T_{n},
          \end{array}\right.,&\text{ with probability (w.p.) }q:=1-p.
    \end{array}\right.\label{rec22}
\end{equation}
\begin{equation}
    W_{n+1}=\left\{\begin{array}{lll}
    \left[V^{(0)}_{n}W_{n}+B_{n}-A_{n}\right]^{+},&\text{ w.p. }p_{1},\,V^{(0)}_{n}\in\{a_{1},\ldots,a_{M}\},a_{k}\in(0,1),k=1,\ldots,M,n\in\mathbb{N}_{0},\\
              \left[V^{(1)}_{n}W_{n}+B_{n}-A_{n}\right]^{+},&\text{ w.p. }p_{2},\,V^{(1)}_{n}\in[0,1),\,n\in\mathbb{N}_{0},\\
        \left[V^{(2)}_{n}W_{n}+B_{n}-A_{n}\right]^{+},&\text{ w.p. }1-p_{1}-p_{2},\,\,V^{(2)}_{n}<0,\,n\in\mathbb{N}_{0},
    \end{array}\right.\label{rec23}
\end{equation}
with $\mathbb{N}_{0}:=\mathbb{N}\cup\{0\}$. We also consider the integer-valued counterpart,
\begin{equation}
    X_{n+1}=\left\{\begin{array}{ll}
         \sum_{k=1}^{X_{n}}U_{k,n}+Z_{n}-Q_{n+1},&X_{n}>0  \\
         Y_{n}-\tilde{Q}_{n+1},&X_{n}=0, 
    \end{array}\right.\label{rec3}
\end{equation}
and a two-dimensional generalization of it, where $x^{+} = max(0, x)$, $x^{-} = min(0, x)$. Moreover, $\{V_{n}\}_{n\in\mathbb{N}_{0}}$, and $\{B_{n}-A_{n}\}_{n\in\mathbb{N}_{0}}$, (similarly $\{\widehat{B}_{n}-A_{n}^{(1)}\}_{n\in\mathbb{N}_{0}}$, $\{\widehat{T}_{n}-A_{n}^{(2)}\}_{n\in\mathbb{N}_{0}}$) are sequences of independent and identically distributed (i.i.d). random variables. For the recursion \eqref{rec2}, the threshold $T_n$ are assumed to be i.i.d. random variables with cumulative distribution function (cdf) $T(\cdot)$ and Laplace–Stieltjes transform (LST) $\tau(\cdot)$. Moreover, $B_{n}$ are i.i.d. random variables with cdf $F_{B} (\cdot)$ and LST $\phi_{B}(\cdot)$.

The ultimate goal of this work is to investigate classes of reflected autoregressive processes described by recursions of the type given above, in which various independence assumptions of $\{B_{n}\}_{n\in\mathbb{N}_{0}}$, $\{A_{n}\}_{n\in\mathbb{N}_{0}}$ are lifted and for which, a detailed exact analysis can be also provided. 

The stochastic recursion \eqref{recu1} where $\{V_{n}\}_{n\in\mathbb{N}_{0}}$ are constant and where $\{B_{n}\}_{n\in\mathbb{N}_{0}}$, $\{A_{n}\}_{n\in\mathbb{N}_{0}}$ are i.i.d. sequences, and also independent on $\{W_{n}\}_{n\in\mathbb{N}_{0}}$ has been treated in \cite{box1}, i.e., the case where $W_{n+1}=[aW_{n}+B_{n}-A_{n}]^{+}$, $n=0,1,\ldots$, with $a\in (0,1)$. The case where $a=1$ corresponds to the classical Lindley recursion describing the waiting time of the classical G/G/1 queue, while the case where $a=-1$ is covered in \cite{vlasiou}. Further progress has been made in \cite{box2}, where additional models described by recursion \eqref{recu1} have been investigated. The work in \cite[Section 3]{box2} is the close to our case, where the authors investigated a recursion where $V$ is either a positive constant with probability $p$, or a random variable taking negative values with probability $1-p$. The fact that $V$ is negative simplified considerably the analysis. 

In \cite{box3}, the authors considered the case where $V_{n}W_{n}$ in \eqref{recu1} was replaced by $S(W_{n})$, where $\{S(t)\}$ is a Levy subordinator (recovering also the case in \cite{box1}, where $S(t)=at$). Note that in \cite{box1,box2,box3} the sequences $\{B_{n}\}_{n\in\mathbb{N}_{0}}$, $\{A_{n}\}_{n\in\mathbb{N}_{0}}$ are assumed to be independent. Recently, in \cite{boxman}, the authors considered Lindley type recursions that arise in queueing and insurance risk models, where the sequences $\{B_{n}\}_{n\in\mathbb{N}_{0}}$, $\{A_{n}\}_{n\in\mathbb{N}_{0}}$ obey a semi-linear dependence. These recursions can also be treated as of autoregressive type. This work is the closer to ours. Moreover, in \cite{adan}, the authors developed a method to study functional equations that arise in a wide range of queueing, autoregressive and branching processes. Finally, the author in \cite{hoo}, considered a generalized version of the model in \cite{box2}, by assuming $V_{n}$ to take values in $(-\infty,1]$. In particular, in \cite{hoo}, the author investigated the recursion \eqref{rec23} for $a=1$. 

The main contribution of this paper is to investigate the transient as well as the stationary behaviour of a wide range of autoregressive processes described in \eqref{recu1}-\eqref{rec3}, by lifting various independence assumptions of the sequences $\{B_{n}\}_{n\in\mathbb{N}_{0}}$, $\{A_{n}\}_{n\in\mathbb{N}_{0}}$. This is accomplished by using Liouville's theorem, and by stating and solving a Wiener-Hopf boundary value problem, or by solving an integral equation, depending on the nature of $\{V_{n}\}_{n\in\mathbb{N}_{0}}$. We have to point out here that to our best knowledge, autoregressive recursions of the form \eqref{rec2}-\eqref{rec3} have not considered in the literature so far.  We also investigate the stationary analysis of $\{X_{n}\}_{n\in\mathbb{N}_{0}}$ in \eqref{rec3}, which represents a novel retrial queueing model. An extension to a two-dimensional case that describes the priority case is also considered. %More precisely, in Section \ref{s2}, we consider versions of recursion \eqref{rec1}, where $V_{n}=a$ with certainty and interarrival times are randomly proportional to service time with additive/subtracting delay, as well a . 
\section{M/G/1-type autoregressive queues with interarrival times randomly proportional to service/system times}\label{s2}
In the following we cope with some autoregressive $M/G/1$-type queueing systems where the interarrival time between the $n$th and the $(n+1)$th job, say $A_{n}$, depends on the service time of the $n$th job, or on the system time after the arrival of the $n$th job.  

\subsection{Interarrival times randomly proportional to service times}
Consider the following variant of the standard autoregressive $M/G/1$ queue: When the service time  equals $x\geq 0$ following the $n$th arrival, then the next interarrival time equals $\beta_{i}x$ (with probability $p_{i}$, $i=1,\ldots,N+M$) increased by an independent  additive delay $J_{n}$. In the following, we consider the recursion \eqref{recu1}, where $P(V_{n}=a)=1$, $a\in(0,1)$.

Let $W_{n}$ the waiting time of the $n$th customer. The interarrival time between the $n$th and the $(n+1)$th customer, say $A_{n}$, satisfies $A_{n}=G_{n}B_{n}+J_{n}$, where $B_{n}$ is the service time of the $n$th customer and $J_{n}$ an additive delay or random jitter. The random variable $G_n$ has finite support. Let $\beta_i$ denote its $i$th value and let $p_i = P(G_{n}=\beta_{i})$ denote the corresponding probability, $i=1,\ldots,N+M$ ($M,N\geq 1$), with $\sum_{i=1}^{N+M}p_{i}=1$. We further assume that the service times and jitter are exponentially distributed: $B_n \sim \exp(\mu)$ and $J_n \sim \exp(\delta)$. Extensions to the case where $J_{n}$ has a rational transform will be also discussed. Thus, the sequence $\{W_{n}\}_{n\in\mathbb{N}_{0}}$ obeys the following recursion:
\begin{equation}
    W_{n+1}=[aW_{n}+(1-G_{n})B_{n}-J_{n}]^{+},
    \label{rec1}
\end{equation}
where $a\in(0,1)$. Without loss of generality and in order to avoid trivial solutions, assume that $1<\beta_{1}<\beta_{2}<\ldots<\beta_{N}$, and $\beta_{N+1}<\beta_{N+1}<\ldots<\beta_{N+M}<1$.
\subsubsection{Transient analysis}
We first focus on the transient distribution, and following the lines in \cite{box1}, let for $|r|<1$,
\begin{displaymath}
Z_{w}(r,s)=\sum_{n=0}^{\infty}r^{n}E(e^{-sW_{n+1}}|W_{0}=w)\, ,\quad U^{-}_{w}(r,s)=\sum_{n=0}^{\infty}r^{n}E(e^{-sU^{-}_{n}}|W_{0}=w),
\end{displaymath}
where $U^{-}_{n}:=[aW_{n}+(1-G_{n})B_{n}-J_{n}]^{-}$. Then, using the property that $1+e^{x}=e^{[x]^{+}}+e^{[x]^{-}}$, \eqref{rec1} leads to
\begin{displaymath}
    \begin{array}{rl}
       E(e^{-sW_{n+1}}|W_{0}=w)= &E(e^{-s(aW_{n}+(1-G_{n})B_{n}+J_{n})}|W_{0}=w)+1-E(e^{-sU_{n}}|W_{0}=w)  \vspace{2mm}\\
        = & E(e^{-saW_{n}}|W_{0}=w)E(e^{sJ_{n}})\sum_{i=1}^{N+M}p_{i}E(e^{-s(1-\beta_{i})B_{n}})+1-E(e^{-sU_{n}}|W_{0}=w)\vspace{2mm}\\
        = & E(e^{-saW_{n}}|W_{0}=w)\frac{\delta}{\delta-s}\sum_{i=1}^{N+M}p_{i}\phi_{B}(\bar{\beta}_{i}s)+1-E(e^{-sU_{n}}|W_{0}=w),
    \end{array}
\end{displaymath}
where $\bar{\beta}_{i}=1-\beta_{i}$, $i=1,\ldots,N+M$, and $\phi_{B}(s)$ being the LST of $B$. Multiplying by $r^{n}$ and taking the sum of $n=0$ to infinity yields
\begin{equation}
    Z_{w}(r,s)-e^{-sw}=r\frac{\delta}{\delta-s}Z_{w}(r,as)\sum_{i=1}^{N+M}p_{i}\phi_{B}(\bar{\beta}_{i}s)+\frac{r}{1-r}-rU^{-}_{w}(r,s).
    \label{f1}
\end{equation}

Assume hereon that $B\sim\exp(\mu)$. Then, $\phi_{B}(\bar{\beta}_{i}s)=\frac{\mu}{\mu+\bar{\beta}_{i}s}=\frac{1}{1-\gamma_{i}s}$, where $\gamma_{i}=\frac{\beta_{i}-1}{\mu}$, $i=1,\ldots,N+M$. Simple calculations imply that
\begin{displaymath}
    \sum_{i=1}^{N+M}\frac{p_{i}}{1-\gamma_{i}s}=\frac{\sum_{i=1}^{N+M}p_{i}\prod_{j\neq i}(1-\gamma_{i}s)}{\prod_{i=1}^{N+M}(1-\gamma_{i}s)}:=\frac{f(s)}{g(s)}.
\end{displaymath}
Note that $g(s)=0$ has $N+M$ distinct and real roots $\gamma_{i}^{-1}$, where $N$ of them are positive and $M$ are negative. In particular, let $s_{i}^{+}=\gamma_{i}^{-1}=\frac{\mu}{\beta_{i}-1}$, $i=1,\ldots,N$ the positive, and $s_{i}^{-}=\gamma_{i}^{-1}=\frac{\mu}{\beta_{i}-1}$, $i=N+1,\ldots,N+M$ the negative roots of $g(s)=0$. Note that
\begin{displaymath}
    g(s)=\prod_{i=1}^{N+M}(1-\gamma_{i}s)=\prod_{i=1}^{N+M}\gamma_{i}(\gamma_{i}^{-1}-s)=\prod_{i=1}^{N+M}(-\gamma_{i})\prod_{i=1}^{N}(s-s_{i}^{+})\prod_{i=N+1}^{N+M}(s-s_{i}^{-}):=g^{+}(s)g^{-}(s),
\end{displaymath}
where $g^{+}(s):=\prod_{i=1}^{N}(s-s_{i}^{+})$, $g^{-}(s):=\prod_{i=1}^{N+M}(-\gamma_{i})\prod_{i=N+1}^{N+M}(s-s_{i}^{-})$. Now \eqref{f1} becomes for $Re(s)=0$:
\begin{equation}
   (\delta-s)g^{+}(s)[ Z_{w}(r,s)-e^{-sw}]-r\delta \frac{f(s)}{g^{-}(s)}Z_{w}(r,as)=(\delta-s)g^{+}(s)[\frac{r}{1-r}-rU_{w}(r,s)].\label{f2}
\end{equation}
Now we make the following observations:
\begin{itemize}
    \item The left-hand side is analytic in $Re(s)>0$, and continuous in $Re(s)\geq 0$,
    \item The right-hand side is analytic in $Re(s)<0$, and continuous in $Re(s)\leq 0$,
    \item $Z_{w}(r,s)$ (resp. $U_{w}(r,s)$) is for $Re(s)\geq 0$ (resp. $Re(s)\leq 0$) bounded by $(1-r)^{-1}$.    
\end{itemize}
Thus, \eqref{f2} represents an entire function. Generalized Liouville's theorem \cite[Theorem 10.52]{tit} states that in their respective half planes, both the left and the right-hand side can be written as the same $(N+1)-$th order polynomial in $s$, for large $s$, i.e.,
\begin{equation}
    (\delta-s)g^{+}(s)[ Z_{w}(r,s)-e^{-sw}]-r\delta\frac{f(s)}{g^{-}(s)}K_{w}(r,s)Z_{w}(r,as)=\sum_{i=0}^{N+1}s^{i}C_{i,w}(r),\,Re(s)\geq 0.\label{f3}
\end{equation}
Note that for $s=0$ \eqref{f3} yields
\begin{displaymath}
    \delta\prod_{i=1}^{N}(-s_{i}^{+})(\frac{1}{1-r}-1)-r\delta
\frac{f(0)}{g^{-}(0)}\frac{1}{1-r}=C_{0,w}(r).
\end{displaymath}
Having in mind that $\frac{f(0)}{g(0)}=1$, so that $\frac{f(0)}{g^{-}(0)}=\prod_{i=1}^{N}(-s_{i}^{+})$, we easily realize that $C_{0.w}(r)=0$. Moreover, setting $s=\delta$, and $s=s^{+}_{i}$, $i=1,\ldots,N$ we obtain the following system of equations for the rest of unknowns $C_{i,w}(r)$, $i=1,\ldots,N$:
\begin{equation}
    \begin{array}{rl}
    -r\delta \frac{f(s^{+}_{i})}{g^{-}(s^{+}_{i})}Z_{w}(r,as^{+}_{i})=&\sum_{i=1}^{N+1}(s^{+}_{i})^{i}C_{i,w}(r),\,i=1,\ldots,N,\vspace{2mm}\\
    -r \delta\frac{f(\delta)}{g^{-}(\delta)}Z_{w}(r,a\delta)=&\sum_{i=1}^{N+1}\delta^{i}C_{i,w}(r).
    \end{array}\label{sys1}
\end{equation}
It remains to obtain $Z_{w}(r,as^{+}_{i})$, $i=1,\ldots,N$, $Z_{w}(r,a\delta)$. These terms are evaluated as follows: Expression \eqref{f3} is now written as
\begin{equation}
     Z_{w}(r,s)=K_{w}(r,s)Z_{w}(r,as)+L_{w}(r,s),\label{iter1}
\end{equation}
where
\begin{displaymath}\begin{array}{lr}
      K_{w}(r,s):= r\frac{\delta}{\delta-s}\frac{f(s)}{g(s)},& L_{w}(r,s):=\frac{\sum_{i=1}^{N+1}s^{i}C_{i,w}(r)}{(\delta-s)g^{+}(s)}+e^{-sw}. 
      \end{array}
\end{displaymath}
Iterating \eqref{iter1} yields
\begin{equation}
    Z_{w}(r,s)=\sum_{n=0}^{\infty}L_{w}(r,a^{n}s)\prod_{j=0}^{n-1}K_{w}(r,a^{j}s),
    \label{solt1}
\end{equation}
with the convention of empty product to be equal to one. Setting $s=\alpha\delta$, and $s=as^{+}_{i}$ in \label{sol1}, we obtain expressions for the $Z_{w}(r,a\delta)$, $Z_{w}(r,as^{+}_{i})$, $i=1,\ldots,N$, respectively. Substituting back in \eqref{sys1}, we obtain an $(N+1)\times (N+1)$ system of equations with $C_{i,w}(r)$, $i=1,\ldots,N+1$ as unknowns. 
\begin{remark}
    It is easily realised in \eqref{solt1} that $Z_{w}(r,s)$ appears to have singularities in $s=\delta/a^{j}$, and $s=s_{i}^{+}/a^{j}$, $i=1,\ldots,N$, $j=0,1,\ldots$. We can show that these are removable singularities. Let us show this for $s=\delta$, and $s=s_{i}^{+}$. We write \eqref{solt1} as follows to isolate the singularities for $s=\delta$, and $s=s_{i}^{+}$:
    \begin{displaymath}
        \begin{array}{rl}
           Z_{w}(r,s)=  &\frac{\sum_{i=1}^{N+1}s^{i}C_{i,w}(r)}{(\delta-s)g^{+}(s)}+e^{-sw}  +\sum_{n=1}^{\infty}(\frac{\sum_{i=1}^{N+1}(a^{n}s)^{i}C_{i,w}(r)}{(\delta-a^{n}s)g^{+}(a^{n}s)}+e^{-sa^{n}w})r^{n}\frac{\delta}{\delta-s}\frac{f(s)}{g(s)}\prod_{j=1}^{n-1}\frac{\delta}{\delta-a^{j}s}\frac{f(a^{j}s)}{g(a^{j}s)}\vspace{2mm}\\
            = & e^{-sw}+\frac{1}{(\delta-s)g^{+}(s)}\left[\sum_{i=1}^{N+1}s^{i}C_{i,w}(r)+r\delta\frac{f(s)}{g^{-}(s)}\sum_{n=1}^{\infty}(\frac{\sum_{i=1}^{N+1}(a^{n}s)^{i}C_{i,w}(r)}{(\delta-a^{n}s)g^{+}(a^{n}s)}+e^{-sa^{n}w})r^{n-1}\prod_{j=1}^{n-1}\frac{\delta}{\delta-a^{j}s}\frac{f(a^{j}s)}{g(a^{j}s)}\right]\vspace{2mm}\\
            =&e^{-sw}+\frac{1}{(\delta-s)g^{+}(s)}\left[\sum_{i=1}^{N+1}s^{i}C_{i,w}(r)+r\delta\frac{f(s)}{g^{-}(s)}Z_{w}(r,as)\right].
        \end{array}
    \end{displaymath}
    \end{remark}
   It is easily realised by using \eqref{sys1}, that the term inside the brackets in the last line vanishes for $s=\delta$, and $s=s_{i}^{+}$, confirming that $s=\delta$, and $s=s_{i}^{+}$, $i=1,\ldots,N$ are not poles of $Z_{w}(r,s)$. Similarly, we can show using \eqref{f3} that $Z_{w}(r,s)$ has no singularity at $s=\delta/a$, $s=s_{i}^{+}/a$, and so on.
\subsubsection{Stationary analysis I}\label{sf1a}
One can follow an analogous analysis as the one presented in the previous sub-subsection or the lines in \cite{box1} and obtain a functional equation similar to \eqref{iter1}. In particular, by applying Abel's theorem, or by considering the relation $ W_{\infty}=[aW_{\infty}+(1-G)B-J]^{+}$, leads to
\begin{equation}
    Z(s)=Z(as)\frac{\delta}{\delta-s}\frac{f(s)}{g(s)}+1-U^{-}(s).\label{nj1}
\end{equation}

 Since $a\in(0,1)$, the stability condition can be ensured as long as $E(log(1+(1-G)B))<\infty$ (see also \cite{box2}). % (due to the standard Loynes argument \cite{loynes}) when $E((1-G)B-J)<0$, i.e.,
%\begin{displaymath}
 %   \begin{array}{c}
  %     E((1-G)B-J)<0 \Leftrightarrow\sum_{i=1}^{N+M}p_{i}\frac{(1-\beta_{i})}{\mu}-\frac{1}{\delta}<0\Leftrightarrow \sum_{i=1}^{N+M}p_{i}\beta_{i}+\frac{\mu}{\delta}>1.
  %  \end{array}
%\end{displaymath}

    Note also that 
    \begin{displaymath}
        [aW_{\infty}+(1-G)B-J]^{-}=\left\{\begin{array}{ll}
             aW_{\infty}+(1-G)B-J,&aW_{\infty}+(1-G)B-J<0,  \\
             0,& aW_{\infty}+(1-G)B-J\geq 0,
        \end{array}\right.
    \end{displaymath}
    thus,
    \begin{displaymath}
        \begin{array}{rl}
            U^{-}(s)= &E(e^{-s(aW_{\infty}+(1-G)B-J)}|aW_{\infty}+(1-G)B-J<0)P(aW_{\infty}+(1-G)B-J<0)\\&+P(aW_{\infty}+(1-G)B-J\geq 0)\vspace{2mm}  \\
             =&\frac{\delta}{\delta-s}P(aW_{\infty}+(1-G)B-J<0)+P(aW_{\infty}+(1-G)B-J\geq 0)\vspace{2mm}\\
             =&1+\frac{s}{\delta-s}P(aW_{\infty}+(1-G)B-J<0),
        \end{array}
    \end{displaymath}
where we used the fact that $E(e^{-s(aW_{\infty}+(1-G)B-J)}|aW_{\infty}+(1-G)B-J<0)$ is the LST of the pdf characterized by:
\begin{displaymath}\begin{array}{rl}
    P(aW_{\infty}+(1-G)B-J\leq x|aW_{\infty}+(1-G)B-J<0)= &P(J\geq aW_{\infty}+(1-G)B-x|J>aW_{\infty}+(1-G)B)  \\
     =&P(J\geq -x)=P(-J\leq x), 
\end{array}
\end{displaymath}
and thus, $E(e^{-s(aW_{\infty}+(1-G)B-J)}|aW_{\infty}+(1-G)B-J<0)=\frac{\delta}{\delta-s}$. Let $P:=P(aW_{\infty}+(1-G)B-J<0)$. Then, \eqref{nj1} is now written as
\begin{equation}
\begin{array}{rl}
    Z(s)=&Z(as)\frac{\delta}{\delta-s}\frac{f(s)}{g(s)}-\frac{s}{\delta-s}P\vspace{2mm}\\
    =&-\frac{Ps}{\delta-s}+\frac{\delta}{\delta-s}\frac{f(s)}{g(s)}[-\frac{Pas}{\delta-as}+\frac{\delta}{\delta-as}\frac{f(as)}{g(as)}Z(a^{2}s)]\\
    =&\ldots\\
    =&-\sum_{n=0}^{\infty}\frac{Pa^{n}s}{\delta-a^{n}s}\prod_{j=0}^{n-1}\frac{f(a^{j}s)\delta}{g(a^{j}s)(\delta-a^{j}s)}+\lim_{n\to\infty}Z(a^{n}s)\prod_{j=0}^{n-1}\frac{f(a^{j}s)\delta}{g(a^{j}s)(\delta-a^{j}s)}\vspace{2mm}\\
    =&-\sum_{n=0}^{\infty}\frac{Pa^{n}s}{\delta-a^{n}s}\prod_{j=0}^{n-1}\frac{f(a^{j}s)\delta}{g(a^{j}s)(\delta-a^{j}s)}+\prod_{j=0}^{\infty}\frac{f(a^{j}s)\delta}{g(a^{j}s)(\delta-a^{j}s)},
\end{array}\label{nj2}
\end{equation}
since $\lim_{n\to\infty}Z(a^{n}s)=Z(0)=1$. Note that $P=P(W_{\infty}=0)$, and is derived by multiplying \eqref{nj2} with $\delta-s$ (i.e., the functional equation before the iterations), and setting $s=\delta$, so that
\begin{displaymath}
    P=Z(a\delta)\frac{f(\delta)}{g(\delta)}.
\end{displaymath}
Setting $s=a\delta$ in \eqref{nj2} (so that to obtain $Z(a\delta)$), and substituting back, yields,
\begin{displaymath}
    P=\frac{\frac{f(\delta)}{g(\delta)}\prod_{j=0}^{\infty}\frac{f(a^{j+1}\delta)}{g(a^{j}s)(1-a^{j+1})}}{1+\frac{f(\delta)}{g(\delta)}\sum_{n=0}^{\infty}\frac{a^{n+1}}{1-a^{n+1}}\prod_{j=0}^{n-1}\frac{f(a^{j+1}\delta)}{g(a^{j+1}\delta)(1-a^{j+1})}}.
\end{displaymath}
%\subsubsection{Asymptotic expansions}
%We now focus on the performance measures $P(W=0)$, $E(W^{k})$, $k=1,2,\ldots$ (in the following to improve readability we drop the the index $\infty$ when considering the stationary version of $W_{n}$) in the regime that $a\to 1$, i.e., a perturbation of the model in \cite[Section 5]{cidon}. For such values of $a$, an arrival does have a huge impact on the work already existed in the queue. In view of the complicated exact expression of $P:=P(W=0)$ as given above, it would be rather interesting to obtain asymptotic expansions of them in the parameter $a$.

\subsubsection{Stationary analysis II}
We now focus on the stationary distribution of $W_{n}$ as $n\to \infty$: $ W_{\infty}=[aW_{\infty}+(1-G)B-J]^{+}$. In the following, we follow a different approach compared to the previous sub-subsection: Denote by $Z(s):=\lim_{n\to\infty}Z_{n}(s)$, where $Z_{n}(s):=E(e^{-sW_{n}})$ and let $f_{W_{n}}(w)$ the probability density function (pdf) of $W_{n}$ (for convenience, assume also $W_{0}=0$). Assume now that $\beta_{i}\in(0,1)$, $i=1,\ldots,K$, where $K=N+M$ (so that $Re(s\bar{\beta}_{i})\geq 0$ for $Re(s)\geq 0$), and service time is arbitrarily distributed with pdf $f_{B}(.)$, and LST $\phi_{B}(.)$. Then, 
\begin{equation}
    \begin{array}{rl}
         Z_{n+1}(s)=&E(e^{-s[aW_{n}+(1-G_{n})B_{n}-J_{n}]^{+}})  \\
         =&\sum_{i=1}^{K}p_{i}\int_{w=0}^{\infty}f_{W_{n}}(w)dw\int_{x=0}^{\infty}f_{B}(x)dx \left\{\int_{y=0}^{aw+\bar{\beta}_{i}x}\delta e^{-\delta y}e^{-s(aw+\bar{\beta}_{i}x-y)}dy+\int_{y=aw+\bar{\beta}_{i}x}^{\infty}\delta e^{-\delta y}dy\right\}\vspace{2mm}\\
    =&\sum_{i=1}^{K}p_{i}\int_{w=0}^{\infty}f_{W_{n}}(w)\int_{x=0}^{\infty}f_{B}(x)[\frac{\delta e^{-s(aw+\bar{\beta}_{i}x)}-se^{-\delta(aw+\bar{\beta}_{i}x)}}{\delta-s}]dwdx\vspace{2mm}\\
    =&\frac{\delta}{\delta-s}Z_{n}(as)\sum_{i=1}^{K}p_{i}\phi_{B}(s\bar{\beta}_{i})-\frac{s}{\delta-s}Z_{n}(a\delta)\sum_{i=1}^{K}p_{i}\phi_{B}(\delta\bar{\beta}_{i}).
    \end{array}
    \label{stat1}
\end{equation}
Letting $n\to\infty$ we get:
\begin{equation}\begin{array}{rl}
    Z(s)= &\frac{\delta}{\delta-s}Z(as)\sum_{i=1}^{K}p_{i}\phi_{B}(s\bar{\beta}_{i})-\frac{s}{\delta-s}Z(a\delta)\sum_{i=1}^{K}p_{i}\phi_{B}(\delta\bar{\beta}_{i}).
\end{array}
    \label{stat2}
\end{equation}
Iterating \eqref{stat2} we get:
\begin{displaymath}
    \begin{array}{rl}
       Z(s)= &-Z(a\delta)\sum_{i=1}^{K}p_{i}\phi_{B}(\bar{\beta}_{i}\delta)\sum_{n=0}^{\infty}\frac{a^{n}s}{\delta-a^{n}s}\prod_{j=0}^{n-1}\frac{\delta}{\delta-a^{j}s}\sum_{i=1}^{K}p_{i}\phi_{B}(a^{j}\bar{\beta}_{i}s) \vspace{2mm}\\
         &+\lim_{n\to\infty}Z(a^{n}s) \prod_{j=0}^{n-1}\frac{\delta}{\delta-a^{j}s}\sum_{i=1}^{K}p_{i}\phi_{B}(a^{j}\bar{\beta}_{i}s)\vspace{2mm}\\
         =&-Z(a\delta)\sum_{i=1}^{K}p_{i}\phi_{B}(\bar{\beta_{i}}\delta)\sum_{n=0}^{\infty}\frac{a^{n}s}{\delta-a^{n}s}\prod_{j=0}^{n-1}\frac{\delta}{\delta-a^{j}s}\sum_{i=1}^{K}p_{i}\phi_{B}(a^{j}\bar{\beta}_{i}s) \vspace{2mm}\\
         &+\prod_{j=0}^{\infty}\frac{\delta}{\delta-a^{j}s}\sum_{i=1}^{K}p_{i}\phi_{B}(a^{j}\bar{\beta}_{i}s),
    \end{array}
\end{displaymath}
since $\lim_{n\to\infty}Z(a^{n}s)=Z(0)=1$. Setting $s=a\delta$ we obtain
\begin{equation}
    Z(a\delta)=\frac{\prod_{j=0}^{\infty}\frac{1}{1-a^{j+1}}\sum_{i=1}^{K}p_{i}\phi_{B}(a^{j+1}\delta\bar{\beta}_{i})}{1+\sum_{i=1}^{K}p_{i}\phi_{B}(\bar{\beta_{i}}\delta)\sum_{n=0}^{\infty}\frac{a^{n+1}}{1-a^{n+1}}\prod_{j=0}^{n-1}\frac{1}{1-a^{j+1}}\sum_{i=1}^{K}p_{i}\phi_{B}(a^{j+1}\delta\bar{\beta_{i}})}.\label{qe}
    \end{equation}

Differentiating \eqref{stat2} with respect to $s$, and setting $s=0$ yields after some algebra
\begin{displaymath}
    E(W_{\infty}):=-\frac{d}{ds}Z(s)|_{s=0}=\frac{E(B)\sum_{i=1}^{K}p_{i}\bar{\beta}_{i}-\frac{1}{\delta}(1-Z(a\delta)\sum_{i=1}^{K}p_{i}\phi_{B}(\delta\bar{\beta}_{i}))}{1-a},
\end{displaymath}
where $Z(a\delta)$ as given in \eqref{qe}, and $E(B)$ the mean service duration.
\begin{remark}
    Note that the analysis can be considerably adapted to consider the case where $J_{n}$ follows a hyperexponential distribution with $L$ phases, i.e., with density function $f_{J}(x):=\sum_{i=1}^{L}q_{i}\delta_{i}e^{-\delta_{i}x}$, $x\geq 0$, $\sum_{i=1}^{L}q_{i}=1$. In such a case
    \begin{displaymath}
        \begin{array}{rl}
             Z(s)=&Z(as)\sum_{i=1}^{L}q_{i}\frac{\delta_{i}}{\delta_{i}-s}\sum_{i=1}^{K}p_{i}\phi_{B}(s\bar{\beta}_{i})-P(aW+(1-\Omega)B-J<0)(1-\sum_{i=1}^{L}q_{i}\frac{\delta_{i}}{\delta_{i}-s}),
        \end{array}
    \end{displaymath}
    where 
    \begin{displaymath}
        \begin{array}{rl}
          P(aW+(1-\Omega)B-J<0)=   &  \sum_{i=1}^{K}p_{i}\int_{0}^{\infty}f_{W_{n}}(w)dw\int_{0}^{\infty}f_{B}(x)dx\int_{aw+\bar{\beta}_{i}x}^{\infty}\sum_{j=1}^{L}q_{j}\delta_{i}e^{-\delta_{i}y}dy\\
             =& \sum_{i=1}^{K}p_{i}\sum_{j=1}^{L}q_{j}\phi_{B}(\delta_{j}\bar{\beta}_{i})Z(\alpha\delta_{j}),
        \end{array}
    \end{displaymath}
    so that
     \begin{equation}
        \begin{array}{rl}
             Z(s)=&Z(as)V(s)\sum_{i=1}^{K}p_{i}\phi_{B}(s\bar{\beta}_{i})-\sum_{i=1}^{K}p_{i}\sum_{j=1}^{L}q_{j}\phi_{B}(\delta_{j}\bar{\beta}_{i})Z(\alpha\delta_{j})(1-V(s)),\label{cvbb}
        \end{array}
     \end{equation}
or equivalently,
\begin{equation}
        \begin{array}{l}
             \prod_{j=1}^{L}(\delta_{j}-s)Z(s)=Z(as)\sum_{j=1}^{L}q_{j}\delta_{j}\prod_{m\neq j}(\delta_{m}-s)\sum_{i=1}^{K}p_{i}\phi_{B}(s\bar{\beta}_{i})\vspace{2mm}\\-\sum_{i=1}^{K}p_{i}\sum_{j=1}^{L}q_{j}\phi_{B}(\delta_{j}\bar{\beta}_{i})Z(\alpha\delta_{j})[\prod_{j=1}^{L}(\delta_{j}-s)-\sum_{j=1}^{L}q_{j}\delta_{j}\prod_{m\neq j}(\delta_{m}-s)],\label{cvbbn}
        \end{array}
     \end{equation}
where 
\begin{displaymath}
    V(s):=\frac{\sum_{j=1}^{L}q_{j}\delta_{j}\prod_{m\neq j}(\delta_{m}-s)}{\prod_{j=1}^{L}(\delta_{j}-s)}.
\end{displaymath}
Note that we first have to derive expressions for the $Z(\alpha\delta_{j})$, $j=1,\ldots,L$. Iterating \eqref{cvbb} yields
\begin{equation}
    \begin{array}{rl}
         Z(s)=&\sum_{i=1}^{K}p_{i}\sum_{j=1}^{L}q_{j}\phi_{B}(\delta_{j}\bar{\beta}_{i})Z(a\delta_{j})\sum_{n=0}^{\infty}\Psi(a^{n}s)\prod_{l=0}^{n-1}\Phi(a^{l}s)+\prod_{l=0}^{\infty}\Phi(a^{l}s),
    \end{array}\label{cvbba}
\end{equation}
where 
\begin{displaymath}
    \begin{array}{lr}
         \Phi(s):= \sum_{i=1}^{K}p_{i}\phi_{B}(s\bar{\beta}_{i})V(s), &\Psi(s):=V(s)-1. 
    \end{array}
\end{displaymath}
Setting in \eqref{cvbba}, $s=a\delta_{p}$, $p=1,\ldots,L$ we obtain an $L\times L$ system of equations with $L$ unknowns, $Z(a\delta_{p})$, $p=1,\ldots,L$:
\begin{displaymath}
    \begin{array}{l}
       Z(a\delta_{p})(1-\sum_{i=1}^{K}p_{i}q_{p}\phi_{B}(\delta_{p}\bar{\beta}_{i})\sum_{n=0}^{\infty}\Psi(a^{n}s)\prod_{l=0}^{n-1}\Phi(a^{l}s))\vspace{2mm}\\-\sum_{i=1}^{K}p_{i}\sum_{j\neq p}q_{j}\phi_{B}(\delta_{j}\bar{\beta}_{i})Z(a\delta_{j})\sum_{n=0}^{\infty}\Psi(a^{n+1}\delta_{p})\prod_{l=0}^{n-1}\Phi(a^{l+1}\delta_{p})=\prod_{l=0}^{\infty}\Phi(a^{l+1}\delta_{p}).
    \end{array}
\end{displaymath}
\end{remark}

\begin{remark}
   Consider the case of a reflected  autoregressive M/G/1 type queue where interarrival times are deterministic proportional dependent to service times with additive delay. We consider the case where $I_{n+1}=bB_{n}+J_{n}$, where $b\in(0,1)$ and $J_{n}\sim \exp(\delta)$. The sequence $(W_{n})_{n}$ obeys the following recursion:
\begin{equation}
    W_{n+1}=[aW_{n}+(1-b)B_{n}-J_{n}]^{+},
    \label{rec1xz}
\end{equation}
where $a\in(0,1)$. Note for $a=1-b$, the sequence \eqref{rec1xz} was investigated in \cite[Section 2]{boxman}. Here we cope with the general case ($a\neq 1-b$), although the analysis follows the lines in \cite{box1}.
\end{remark}
\subsection{Proportional dependency with additive and subtracting delay}
We now focus on the case where the interarrival times are such that $A_{n}=cB_{n}+J_{n}$, with
\begin{equation}
    J_{n}:=\left\{\begin{array}{ll}
         \tilde{J}_{n}&,\text{ with probability }p,  \\
         -\widehat{J}_{n}&,  \text{ with probability }q:=1-p,
    \end{array}\right.\label{form}
\end{equation}
where $\tilde{J}_{n}\sim exp(\delta)$, $\widehat{J}_{n}\sim exp(\nu)$. Now the sequence $(W_{n})_{n}$ obeys $W_{n+1}=[aW_{n}+B_{n}-(cB_{n}+J_{n})^{+}]^{+}$. With probability $p$, $J_{n}=\tilde{J}_{n}$, and thus $(cB_{n}+J_{n})^{+}=cB_{n}+\tilde{J}_{n}$, while with probability $q$, $J_{n}=\widehat{J}_{n}$, and thus $(cB_{n}+J_{n})^{+}=(cB_{n}-\widehat{J}_{n})^{+}$. Therefore,
\begin{equation}
         E(e^{-sW_{n+1}})=pE(e^{-s[aW_{n}+\bar{c}B_{n}-\tilde{J}_{n}]^{+}})+q E(e^{-s[aW_{n}+B_{n}-(cB_{n}-\widehat{J}_{n})^{+}]}),\label{sa}
         \end{equation}
where $\bar{c}:=1-c$. Assume that the steady-state waiting time distribution exists and denote by $Z(s)$ its LST. Following similar steps as in the derivation of \eqref{stat1}, we can obtain:
\begin{displaymath}
    E(e^{-s[aW_{n}+\bar{c}B_{n}-\tilde{J}_{n}]^{+}})=\frac{\delta}{\delta-s}Z(as)\phi_{B}(s\bar{c})-\frac{s}{\delta-s}Z(a\delta)\phi_{B}(\delta\bar{c}).
\end{displaymath}
Now
\begin{displaymath}
    \begin{array}{rl}
        E(e^{-s[aW_{n}+B_{n}-(cB_{n}-\widehat{J}_{n})^{+}]}= &\int_{w=0}^{\infty}\int_{x=0}^{\infty}f_{B}(x)\left[\int_{y=0}^{cx}e^{-s(aw+\bar{c}x+y)}\nu e^{-\nu y}dy+\int_{y=cx}^{\infty}e^{-s(aw+x)}\nu e^{-\nu y}dy\right]dxdP(W<w)\vspace{2mm}\\  \\
         =&\frac{\nu}{\nu+s}Z(as)\phi_{B}(s\bar{c})+\frac{s}{\nu+s}Z(as)\phi_{B}(s+\nu c).
    \end{array}
\end{displaymath}
Thus, \eqref{sa} reads 
\begin{displaymath}
    Z(s)=H(s)Z(as)+L(s),
\end{displaymath}
where
\begin{displaymath}
    \begin{array}{rl}
        H(s)=&\phi_{B}(s\bar{c})(\frac{p\delta}{\delta-s}+\frac{\nu q}{\nu+s})+ \frac{qs}{\nu+s}Z(as)\phi_{B}(s+\nu c),\vspace{2mm}  \\
        L(s)= & -\frac{s}{\delta-s}pZ(a\delta)\phi_{B}(\delta\bar{c}):=-\frac{s}{\delta-s}P.
    \end{array}
\end{displaymath}
Iterating as in subsection \ref{sf1a}, and having in mind that $\lim_{n\to\infty}Z(a^{n}s)=1$, we arrive at
\begin{displaymath}
    Z(s)=-P\sum_{n=0}^{\infty}\frac{a^{n}s}{\delta-a^{n}s}\prod_{j=0}^{n-1}H(a^{j}s)+\prod_{j=0}^{\infty}H(a^{j}s).
\end{displaymath}
Setting $s=a\delta$, and substitute back, we obtain
\begin{displaymath}
    P=\frac{p\phi_{B}(\delta\bar{c})\prod_{j=0}^{\infty}H(a^{j+1}\delta)}{1+p\phi_{B}(\delta\bar{c})\sum_{n=0}^{\infty}\frac{a^{n+1}}{1-a^{n+1}}\prod_{j=0}^{n-1}H(a^{j+1}\delta)}.
\end{displaymath}
\begin{remark}
    One may also consider the case where the interarrival times are related to the previous service time as follows: $A_{n}=G_{n}B_{n}+J_{n}$, where $J_{n}$, as given in \eqref{form}, and $G_{n}$ are iid random variables with probability mass function given by $P(G_{n}=c_{k})=p_{k}$, $k=1,\ldots,N$, $\sum_{k=1}^{N}p_{k}=1$, $c_{k}\in(0,1)$. In particular, \eqref{sa} now becomes
    \begin{equation}
         E(e^{-sW_{n+1}})=pE(e^{-s[aW_{n}+(1-G_{n})B_{n}-\tilde{J}_{n}]^{+}})+q E(e^{-s[aW_{n}+B_{n}-(G_{n}B_{n}-\widehat{J}_{n})^{+}]}),\label{sa1}
         \end{equation}
    and following the same arguments as above, we again have
    \begin{displaymath}
    Z(s)=H(s)Z(as)+L(s),
\end{displaymath}
where now
\begin{displaymath}
    \begin{array}{rl}
        H(s):=&\sum_{k=1}^{N}p_{k}[\phi_{B}(s\bar{c}_{k})(\frac{p\delta}{\delta-s}+\frac{\nu q}{\nu+s})+ \frac{qs}{\nu+s}\phi_{B}(s+\nu c_{k})],\vspace{2mm}  \\
        L(s):= & -\frac{s}{\delta-s}pZ(a\delta)\sum_{k=1}^{N}p_{k}\phi_{B}(\delta\bar{c}_{k}):=-\frac{s}{\delta-s}P.
    \end{array}
\end{displaymath}
Following the lines in subsection \ref{sf1a}, and having in mind that $\lim_{n\to\infty}Z(a^{n}s)=1$, we obtain the desired expression for $Z(s)$.
\end{remark}
\begin{remark}
    The case where $\tilde{J}_{n}$, $\widehat{J}_{n}$ are i.i.d. random variables following a distribution with rational LST can also be treated similarly. In particular, assume that $\tilde{J}_{n}$, $\widehat{J}_{n}$ follow hyperexponential distributions, i.e., their pdfs are $f_{\tilde{J}}(x):=\sum_{i=1}^{L}q_{i}\delta_{i}e^{-\delta_{i}x}$, and $f_{\widehat{J}}(x):=\sum_{i=1}^{M}h_{i}\nu_{i}e^{-\nu_{i}x}$, respectively. Then, following similar arguments as above, and assuming $A_{n}=G_{n}B_{n}+J_{n}$, where $J_{n}$, as given in \eqref{form}, we obtain after lengthy computations:
    \begin{equation}
          Z(s)=H(s)Z(as)+L(s),\label{xde}
    \end{equation}
where now
\begin{displaymath}
    \begin{array}{rl}
        H(s):=&\sum_{k=1}^{N}p_{k}[\phi_{B}(s\bar{c}_{k})(p\sum_{j=1}^{L}\frac{\delta_{j}q_{j}}{\delta_{j}-s}+q\sum_{j=1}^{M}\frac{\nu_{j}h_{j}}{\nu_{j}+s})+ qs\sum_{j=1}^{M}\frac{h_{j}}{\nu_{j}+s}\phi_{B}(s+\nu_{j} c_{k})],\vspace{2mm}  \\
        L(s):= &-sp \sum_{k=1}^{N}p_{k}\sum_{j=1}^{L}\frac{q_{j}}{\delta_{j}-s}Z(a\delta_{j})\phi_{B}(\delta_{j}\bar{c}_{k}).
    \end{array}
\end{displaymath}
Iterating \eqref{xde} as in subsection \ref{sf1a}, and having in mind that $\lim_{n\to\infty}Z(a^{n}s)=1$, we obtain the desired expression for $Z(s)$.
\end{remark}
\subsection{Interarrival times randomly proportional to system time}\label{cc1}
Consider the following variant of the standard $M/G/1$ queue: When the workload just after the $n$th arrival equals $x\geq 0$, then the next interarrival time equals $\beta_{i}x$ (with probability $p_{i}$) increased by a random jitter $J_{n}\sim \exp(\delta)$. Thus, $A_{n}=G_{n}(W_{n}+B_{n})+J_{n}$, where $P(G_{n}=\beta_{i}):=p_{i}$, $i=1,\ldots,K$, $\beta_{i}\in(0,1)$. Note that our model generalizes the one in \cite[Section 2]{boxman}, in which $P(G_{n}=c):=1$, i.e., $\beta_{1}=c\in(0,1)$, $\beta_{i}=0$, $i\neq 1$.
\begin{equation}
    W_{n+1}=[(1-G_{n})W_{n}+(1-G_{n})B_{n}-J_{n}]^{+},\label{hu}
\end{equation}
Note that the recursion $\eqref{hu}$ is a special case of the recursion \eqref{recu1} with $V_{n}:=1-G_{n}$.

In the following, we assume that the steady-state waiting time distribution exists. Letting $n\to\infty$ in \eqref{hu}, and study the limiting random variable $W$ we have
\begin{displaymath}
    \begin{array}{rl}
        Z(s):=E(e^{-sW_{n+1}})= &\sum_{i=1}^{K}p_{i}\int_{0}^{\infty}\int_{0}^{\infty}f_{B_{n}}(x)dx[\int_{0}^{\bar{\beta}_{i}(w+x)}e^{-s(\bar{\beta}_{i}(w+x)-y)}\delta e^{-\delta y}dy+\int_{\bar{\beta}_{i}(w+x)}^{\infty}\delta e^{-\delta y}dy]dP(W_{n}<w)\vspace{2mm} \\
         =&\frac{\delta}{\delta-s}\sum_{i=1}^{K}p_{i}\phi_{B}(s\bar{\beta}_{i})Z(s\bar{\beta}_{i})-\frac{s}{\delta-s}\sum_{i=1}^{K}p_{i}\phi_{B}(\delta\bar{\beta}_{i})Z(\delta\bar{\beta}_{i}). 
    \end{array}
\end{displaymath}
It is easy to show that $P(J>\bar{\beta}_{i}(W+B))=\phi_{B}(\delta\bar{\beta}_{i})Z(\delta\bar{\beta}_{i})$. Thus, $P(W=0)=\sum_{i=1}^{K}p_{i}P(J>\bar{\beta}_{i}(W+B))$, and therefore, 
\begin{equation}
    Z(s)=\frac{\delta}{\delta-s}\sum_{i=1}^{K}p_{i}\phi_{B}(s\bar{\beta}_{i})Z(s\bar{\beta}_{i})-\frac{s}{\delta-s}P(W=0).\label{klo}
\end{equation}
Following \cite{adan}, we can obtain
\begin{displaymath}
    \begin{array}{c}
        Z(s)=\sum_{k=0}^{\infty}\sum_{i_{1}+\ldots+i_{K}=k}p_{1}^{i_{1}}\ldots p_{K}^{i_{K}}L_{i_{1},\ldots,i_{K}}(s)K(\bar{\beta}_{1}^{i_{1}}\ldots\bar{\beta}_{K}^{i_{K}}s)\vspace{2mm}\\
        +\lim_{k\to\infty}\sum_{i_{1}+\ldots+i_{K}=k}p_{1}^{i_{1}}\ldots p_{K}^{i_{K}}L_{i_{1},\ldots,i_{K}}(s),
    \end{array}
\end{displaymath}
where $K(s):=\frac{s}{\delta-s}P(W=0)$, and $L_{0,0,\ldots,0,1,0,\ldots,0}(s):=\phi_{B}(\bar{\beta}_{k}s)$, with 1 in position $k$, and $k=1,\ldots, K$,
\begin{displaymath}
    L_{i_{1},\ldots,i_{K}}(s)=\phi_{B}(\bar{\beta}_{1}^{i_{1}}\ldots\bar{\beta}_{K}^{i_{K}}s)\sum_{j=1}^{K}L_{i_{1},\ldots,i_{j}-1,\ldots,i_{K}}(s).
\end{displaymath}
\begin{remark}
    A similar analysis can be applied in order to investigate recursions of the form $W_{n+1}=[V_{n}W_{n}+(1-G_{n})B_{n}-J_{n}]^{+}$, where $V_{n}$ i.i.d. random variables with $P(V_{n}=\gamma_{i})=q_{i}$, $i=1,\ldots,K$, $\gamma_{i}\in (0,1)$.
\end{remark}
\subsubsection{Asymptotic expansions}
In the following we focus on deriving asymptotic expansions of the basic performance metrics $P(W=0)$, $E(W^{l})$, $l=1,2,\ldots$, by perturbating $\beta_{i}$s, i.e., by letting in \eqref{klo} $\beta_{i}$ to be equal to $\beta_{i}\epsilon$ with $\epsilon$ very small. The, \eqref{klo} is written:
\begin{displaymath}
    (\delta-s)Z(s)=\delta\sum_{i=1}^{K}p_{i}\phi_{B}(s(1-\beta_{i}\epsilon))Z(s(1-\beta_{i}\epsilon))-sP(W=0).
\end{displaymath}
Note that for $\epsilon=0$, the above equation provides the LST of the waiting time (say $\tilde{W}$) of the classical $M/G/1$ queue where arrivals occur according to a Poisson process with rate $\delta$. So when $\epsilon\to 0$, there is a weak dependence between soji=ourn time and the subsequent interarrival time. Following \cite[subsection 2.3]{boxman}, consider the Taylor series development of $P(W=0)$, $E(W^{l})$, $l=1,\ldots,L$ up to $\epsilon^{m}$ terms for $m\in \mathbb{N}$. So for $\epsilon\to 0$:
\begin{equation}
    \begin{array}{rl}
         P(W=0)=&P(\tilde{W}=0)+\sum_{h=1}^{m}R_{0,h}\epsilon^{m}+O(e^{m}),\\
         E(W^{l})  = &E(\tilde{W}^{l})+\sum_{h=1}^{m}R_{l,h}\epsilon^{m}+O(e^{m}).
    \end{array}\label{vol1}
\end{equation}
Differentiating the functional equation with respect to $s$, setting $s=0$ yields for $\rho=\delta E(B)$
\begin{displaymath}
    E(W)=\frac{P(W=0)-(1-\rho)-\rho\epsilon\sum_{i=1}^{K}p_{i}\beta_{i}}{\delta\epsilon\sum_{i=1}^{K}p_{i}\beta_{i}}.
\end{displaymath}
Simple calculations imply that 
\begin{displaymath}
    \begin{array}{rl}
        R_{0,1}= &(\delta E(\tilde{W}) +\rho)\sum_{i=1}^{K}p_{i}\beta_{i}, \\
         \delta\sum_{i=1}^{K}p_{i}\beta_{i}R_{1,h-1}=&R_{0,h},\,h=2,3,\ldots. 
    \end{array}
\end{displaymath}
Assuming that the first $L$ moments of $W$ are well defined, we subsequently differentiate the functional equation above $l=2,\ldots,L$ times with respect to $s$, and set $s=0$. Then, for $l=2,3,\ldots,L$ we have:
\begin{equation}
    \begin{array}{c}
   \delta(1-\sum_{i=1}^{K}p_{i}(1-\beta_{i}\epsilon)^{l})E(W^{l})=-lE(W^{l-1})+\delta\sum_{i=1}^{K}p_{i}(1-\beta_{i}\epsilon)^{l}\sum_{j=0}^{l-1}\binom{l}{j}E(W^{j})E(B^{l-j}).
\end{array}\label{vol}
\end{equation}
Setting $\epsilon=0$, and having in mind that $\sum_{i=1}^{K}p_{i}=1$ we recover the recursive formula to obtain the moments of the standard $M/G/1$ queue:
\begin{displaymath}
    0=-lE(\tilde{W}^{l-1})+\delta\sum_{j=0}^{l-1}\binom{l}{j}E(\tilde{W}^{j})E(B^{l-j}),\,l=2,3,\ldots,L.
\end{displaymath}
Then, substituting \eqref{vol1} in \eqref{vol} we have 
\begin{equation}
    \begin{array}{rl}
        \delta(1-\sum_{i=1}^{K}p_{i}(1-\beta_{i}\epsilon)^{l})\sum_{h=1}^{m}R_{l,h}\epsilon^{h}= &-l\sum_{h=1}^{m}R_{l-1,h}\epsilon^{h}+\delta\sum_{i=1}^{K}p_{i}(1-\beta_{i}\epsilon)^{l})\sum_{j=0}^{l-1}\binom{l}{j}E(B^{l-j})\sum_{h=1}^{m}R_{j,h}\epsilon^{h}  \vspace{2mm}\\
         &+\delta(\sum_{i=1}^{K}p_{i}(1-\beta_{i}\epsilon)^{l}-1)\sum_{j=0}^{l}\binom{l}{j} E(\tilde{W}^{j})E(B^{l-j}).
    \end{array}\label{nhb}
\end{equation}
Equating $\epsilon$ factors on both sides we obtain $R_{l-1,1}$ in terms of $R_{l-2,1},\ldots,R_{0,1}$, as well as in terms of $E(\widehat{W}^{n})$ obtained above. Thus, since $R_{0,1}$ is known, all $R_{l,1}$ can be derived by:
\begin{displaymath}
    R_{l-1,1}=\frac{1}{1-\delta E(B)}[\frac{\delta}{l}\sum_{n=0}^{l-2}\binom{l}{n}E(B^{l-n})R_{n,1}-\delta\sum_{i=1}^{K}p_{i}\beta_{i}\sum_{n=0}^{l}\binom{l}{n} E(\tilde{W}^{n})E(B^{l-n})].
\end{displaymath}
Similarly, for $h=2$,
\begin{displaymath}
    \begin{array}{rl}
         R_{l-1,2}=&\frac{1}{1-\delta E(B)}[ \frac{\delta}{l}\sum_{n=0}^{l-2}\binom{l}{n}E(B^{l-n})R_{n,2} -\frac{\delta}{l}\sum_{i=1}^{K}p_{i}\beta_{i}\sum_{n=0}^{l}\binom{l}{n}E(B^{l-n})R_{n,1}\vspace{2mm}\\
         &+\delta\frac{l-1}{2}\sum_{i=1}^{K}p_{i}\beta_{i}^{2} \sum_{n=0}^{l}\binom{l}{n} E(\tilde{W}^{n})E(B^{l-n})].
    \end{array}
\end{displaymath}
Similarly, we can obtain $R_{k-1,h}$ in terms of $R_{k,h-1}$ and $R_{n,l}$, $n+l\leq l-2+h$. The procedure we follow to recursively obtain $R_{l,h}$ is the same as the one given in \cite[subsection 2.3]{boxman}, so further details are omitted. 

\section{The single-server queue with service time randomly dependent on waiting time}
Customers arrive according to a Poisson process with rate $\lambda$. If the waiting time $W_n$ of
the $n$th arriving customer, then her service time equals $[B_n - \Omega_{n} W_{n}]^{+}$, with $P(\Omega_{n}=a_{l})=g_{l}$, $l=1,\ldots,K$, with $a_{l}\in(0,1)$, where $\{B_n\}_{n\in\mathbb{N}_{0}}$ is a sequence of independent, exponentially distributed random variables with rate $\mu$, independent of anything else. Note that when the waiting time is very large the service requirement tends to zero, which can explained as an abandonment.

Let $Z(s):=E(e^{-sW_{n}})$, and $A_{n}$ i.i.d. random variables from $exp(\lambda)$. Then,
\begin{equation}
    \begin{array}{rl}
       Z(s):= &E(e^{-sW_{n+1}})=E(e^{-s[W_{n}+[B_n - \Omega_{n} W_{n}]^{+}-A_{n}]^{+}})  \vspace{2mm}\\
       = & E(e^{-s[W_{n}+[B_n - \Omega_{n} W_{n}]^{+}-A_{n}]}) +1-E(e^{-s[W_{n}+[B_n - \Omega_{n} W_{n}]^{+}-A_{n}]^{-}})\vspace{2mm}\\
       =&\sum_{l=1}^{K}g_{l}E(e^{sA_{n}})E(e^{-s[W_{n}+[B_n - a_{l} W_{n}]^{+}]})+1-E(e^{-sU_{n}}),
    \end{array}\label{r1}
\end{equation}
where $U_{n}:=[W_{n}+[B_n - \Omega_{n} W_{n}]^{+}-A_{n}]^{-}$. Assuming $n\to\infty$ and studying the limiting variable $W$ we have
\begin{displaymath}
    \begin{array}{rl}
         E(e^{-s[W_{n}+[B_n - a_{l} W_{n}]^{+}]})=&\int_{w=0}^{\infty}\left[\int_{x=0}^{a_{l}w}\mu e^{-\mu x}e^{-sw}dx+\int_{x=a_{l}w}^{\infty}e^{-s(x+(1-a_{l})w)}\mu e^{-\mu x}dx\right]dP(W<w) \vspace{2mm} \\
        = & \int_{w=0}^{\infty}(e^{-sw}-e^{-(a_{l}\mu+s)w})dP(W<w)+\frac{\mu}{\mu+s}\int_{w=0}^{\infty}e^{-w(s+a_{l}\mu)}dP(W<w) \vspace{2mm} \\
        =&Z(s)-\frac{s}{\mu+s}Z(s+a_{l}\mu).
    \end{array}
\end{displaymath}

Moreover, since
\begin{displaymath}
    [W+(B-a_{l}W)-A]^{-}=\left\{\begin{array}{ll}
         W+(B-a_{l}W)-A,&W+(B-a_{l}W)-A<0,  \\
         0,&W+(B-a_{l}W)-A\geq 0, 
    \end{array}\right.
\end{displaymath}
\begin{displaymath}
    \begin{array}{rl}
         E(e^{-sU_{n}})=& E(e^{-s[W+[B - a_{l} W]^{+}-A]}|A>W+[B - a_{l} W]^{+})P(A>W+[B - a_{l} W]^{+})\vspace{2mm}\\
         & +P(A\leq W+[B- a_{l} W]^{+})\vspace{2mm} \\
         =&\frac{\lambda}{\lambda-s}P(A>W+[B - a_{l} W]^{+})
         +P(A\leq W+[B - a_{l} W]^{+}) \vspace{2mm}\\
         =&1+\frac{s}{\lambda-s}P(A>   W+[B - a_{l} W]^{+}).
    \end{array}
\end{displaymath}
Note that 
\begin{displaymath}
   \begin{array}{rl}
        P(A> W+[B - a_{l} W]^{+})=&\int_{w=0}^{\infty}\left(\int_{x=0}^{a_{l}w}\mu e^{-\mu x}dx\int_{y=w}^{\infty}\lambda e^{-\lambda y}dydx+\right.\vspace{2mm}\\
        &\left.\int_{x=a_{l}w}^{\infty}\mu e^{-\mu x}dx\int_{y=x+(1-a_{l})w}^{\infty}\lambda e^{-\lambda y}dydx\right)dP(W<w)  \vspace{2mm}\\
        =& \int_{w=0}^{\infty}\left(e^{-\lambda w}(1-e^{-\mu a_{l}w})+\frac{\mu}{\mu+\lambda}e^{-(\lambda+\mu a_{l}w)}\right)dP(W<w)\vspace{2mm}\\
        =&Z(\lambda)-\frac{\lambda}{\mu+\lambda}Z(\lambda+\mu a_{l}).
   \end{array}
\end{displaymath}
\begin{remark}
    Note that $P(A> W+[B - \Omega W]^{+})=P(W=0):=\pi_{0}$.
\end{remark}
Thus, substituting back in \eqref{r1} we arrive after simple calculations in:
\begin{equation}
    Z(s)=\frac{\lambda}{\mu+s}\sum_{l=1}^{K}g_{l}Z(s+a_{l}\mu)+C,
    \label{r2}
\end{equation}
where $C:=Z(\lambda)-\frac{\lambda}{\mu+\lambda}\sum_{l=1}^{K}g_{l}Z(\lambda+\mu a_{l})=\pi_{0}$. For $s=0$, \eqref{r2} yields $\sum_{l=1}^{K}g_{l}Z(\mu a_{l})=\frac{\mu}{\lambda}(1-\pi_{0})$. Note also that $Z(\mu a_{l})=P(B>a_{l}W)$, and $\sum_{l=1}^{K}g_{l}Z(\mu a_{l})=P(B>\Omega W)$. 

We need to iterate \eqref{r2} and having in mind that as $s\to\infty$, $Z(s)\to 0$ (needs some work). Note that such kind of recursions were treated in \cite{adan}, since the cummutativity of $\zeta_{l}(s):=s+a_{l}\mu$ and $\zeta_{m}(s):=s+a_{m}\mu$, i.e., $\zeta_{l}(\zeta_{m}(s))=\zeta_{m}(\zeta_{l}(s))$ makes the recursion \eqref{r2} relatively easy to handle, although in each iteration, any term gives rise to new $K$ terms; see also \cite[Remark 5.3]{boxman}. Extensions to the case where service time distributions have rational LST are doable; e.g. hyperexponential.

\section{Threshold-type dependence among interarrival and service times}\label{thre}
\subsection{The simple case}
Customers arrive with a service request at a single server.
Service requests of successive customers are independent, identically distributed (i.i.d.) random variables $B_i$,
$i=1,2,\ldots$ with distribution $B(.)$, mean $b$ and Laplace-Stieltjes transform $(LST)$ $\beta(.)$. Upon arrival, the
service request is registered. If the service request $B_i$ is less than a threshold $T_i$, then the next interarrival
interval, say $J_{i}^{(0)}$, is exponentially distributed with rate $\lambda_{0}$; otherwise, the service time becomes exactly equal to $T_i$ (is
cut off at $T_i$), and the next interarrival interval, say $J_{i}^{(1)}$, is exponentially distributed with rate $\lambda_{1}$. We assume that an arrival makes obsolete a fixed fraction $1-a_{0}$ (resp. $1-a_{1}$) of the work that is already present, with $a_{k}\in(0,1)$, $k=0,1$. We assume the threshold $T_i$ to be i.i.d. random variables with general distribution $T(\cdot)$, that has mean $\tau$ and LST $\tau(\cdot)$. Let also for $Re(s)\geq 0$
\begin{displaymath}
    \begin{array}{rl}
        \chi(s):= &E(e^{-sB}1(B<T))=\int_{0}^{\infty}e^{-sx}(1-T(x))dB(x),  \vspace{2mm}\\
        \psi(s):= &E(e^{-sT}1(B\geq T))=\int_{0}^{\infty}e^{-sx}(1-B(x))dT(x),
    \end{array}
\end{displaymath}
with 
\begin{displaymath}
    \chi(s)+\psi(s)=E(e^{-smin(B,T)}).
\end{displaymath}
Let $W_{n}$ the waiting time of the nth arriving customer, $n=1,2,\ldots$. Then,
\begin{equation}
    W_{n+1}=\left\{\begin{array}{ll}
         \left[a_{0}W_{n}+B_{n}-J^{(0)}_{n}\right]^{+},&\,B_{n}<T_{n},\\
         \left[a_{1}W_{n}+T_{n}-J^{(1)}_{n}\right]^{+},&\,B_{n}\geq T_{n},
    \end{array}\right.
\end{equation}
with $J^{(k)}_{n}\sim exp(\lambda_{k})$, $k=0,1$. Assume that $W_{1}=w$, and let $E_{w}(e^{-sW_{n}}):=E(e^{-sW_{n}}|W_{1}=w)$. Then,
\begin{equation}
    \begin{array}{rl}
        E_{w}(e^{-sW_{n+1}}) = & E_{w}(e^{-s[a_{0}W_{n}+B_{n}-J^{(0)}_{n}]^{+}}1(B_{n}<T_{n}))+ E_{w}(e^{-s[a_{1}W_{n}+T_{n}-J^{(1)}_{n}]^{+}}1(B_{n}\geq T_{n}))\vspace{2mm}\\
         =&  E_{w}(e^{-s[a_{0}W_{n}+B_{n}-J^{(0)}_{n}]}1(B_{n}<T_{n})+ E_{w}(e^{-s[a_{1}W_{n}+T_{n}-J^{(1)}_{n}]}1(B_{n}\geq T_{n}))\vspace{2mm}\\
         &+1- E_{w}(e^{-s[a_{0}W_{n}+B_{n}-J^{(0)}_{n}]^{-}}1(B_{n}<T_{n}))-E_{w}(e^{-s[a_{1}W_{n}+T_{n}-J^{(1)}_{n}]^{-}}1(B_{n}\geq T_{n}))\vspace{2mm}\\
         =&E_{w}(e^{-sa_{0} W_{n}})E(e^{sJ^{(0)}_{n}})E(e^{-sB_{n}}1(B_{n}<T_{n}))\vspace{2mm}\\&+E_{w}(e^{-sa_{1} W_{n}})E(e^{sJ^{(1)}_{n}})E(e^{-sT_{n}}1(B_{n}\geq T_{n}))+1-U^{-}_{w,n}(s),
    \end{array}\label{o1}
\end{equation}
where $U^{-}_{w,n}(s):=E_{w}(e^{-s[a_{0}W_{n}+B_{n}-J^{(0)}_{n}]^{-}}1(B_{n}<T_{n}))+E_{w}(e^{-s[a_{1}W_{n}+T_{n}-J^{(1)}_{n}]^{-}}1(B_{n}\geq T_{n}))$. Note that $U^{-}_{w,n}(s)$ is analytic in $Re(s)\leq 0$. Let
\begin{displaymath}
\begin{array}{rl}
     Z_{w}(r,s):=&\sum_{n=1}^{\infty}r^{n}E_{w}(e^{-sW_{n}}),\,Re(s)\geq 0, \vspace{2mm} \\
     M_{w}(r,s):=&\sum_{n=1}^{\infty}r^{n}U^{-}_{w,n}(s),\,Re(s)\leq 0. 
\end{array}
\end{displaymath}
Then, \eqref{o1} leads for $Re(s)=0$ to:
\begin{equation}
    Z_{w}(r,s)-re^{-sw}=r\frac{\lambda_{0}}{\lambda_{0}-s}\chi(s)Z_{w}(r,a_{0}s)+r\frac{\lambda_{1}}{\lambda_{1}-s}\psi(s)Z_{w}(r,a_{1}s)+\frac{r^{2}}{1-r}-rM_{w}(r,s).\label{o2}
\end{equation}
Multiplying \eqref{o2} by $\prod_{k=0}^{1}(\lambda_{k}-s)$, we obtain
\begin{equation}
\begin{array}{r}
    \prod_{k=0}^{1}(\lambda_{k}-s)(Z_{w}(r,s)-re^{-sw})-r(\lambda_{0}(\lambda_{1}-s)\chi(s)Z_{w}(r,a_{0}s)+\lambda_{1}(\lambda_{0}-s)\psi(s)Z_{w}(r,a_{1}s))\vspace{2mm}\\=\prod_{k=0}^{1}(\lambda_{k}-s)\left(\frac{r^{2}}{1-r}-rM_{w}(r,s)\right).
\end{array}\label{o3}
\end{equation}
Our objective is to obtain $Z_{w}(r,s)$, $M_{w}(r,s)$ by formulating and solving a Wiener-Hopf boundary value problem. A few observations:
\begin{itemize}
    \item The LHS in \eqref{o3} is analytic in $Re(s) > 0$, and continuous in $Re(s)\geq  0$.
    \item The RHS in \eqref{o3} is analytic in $Re(s)<0$, and continuous in $Re(s)\leq  0$.
\item $Z_{w}(r,s)$ is for $Re(s)\geq 0$ bounded by $|\frac{r}{1-r}|$, so by the generalized Liouville's theorem \cite[Theorem 10.52]{tit}, the LHS is at most a quadratic polynomial in $s$ (dependent on $r$) for large $s$, $Re(s)>0$.
\item  $M_{w}(r,s)$ is for $Re(s)\leq 0$ bounded by $|\frac{r}{1-r}|$, so by the generalized Liouville's theorem \cite[Theorem 10.52]{tit}, the RHS is at most a quadratic polynomial in $s$ (dependent on $r$) for large $s$, $Re(s) < 0$.
\end{itemize}

Thus, 
\begin{equation}
    \begin{array}{l}
        \prod_{k=0}^{1}(\lambda_{k}-s)(Z_{w}(r,s)-re^{-sw})-r(\lambda_{0}(\lambda_{1}-s)\chi(s)Z_{w}(r,a_{0}s)\vspace{2mm}\\+\lambda_{1}(\lambda_{0}-s)\psi(s)Z_{w}(r,a_{1}s))= C_{0,w}(r)+sC_{1,w}(r)+s^{2}C_{2,w}(r),\,Re(s)\geq 0,
\end{array}\label{cv1}
\end{equation}
\begin{equation}
     \prod_{k=0}^{1}(\lambda_{k}-s)\left(\frac{r^{2}}{1-r}-rM_{w}(r,s)\right)= C_{0,w}(r)+sC_{1,w}(r)+s^{2}C_{2,w}(r),\,Re(s)\leq 0,
\end{equation}
with $C_{i,w}(r)$, $i=0,1,2,$ functions of $r$ to be determined.

Taking $s=0$ in \eqref{cv1} yields
\begin{displaymath}
    \lambda_{0}\lambda_{1}(\frac{r}{1-r}-r)-r(\chi(0)+\psi(0))\lambda_{0}\lambda_{1}\frac{r}{1-r}=C_{0,w}(r),
\end{displaymath}
and having in mind that $\chi(0)+\psi(0)=1$, $C_{0,w}(r)=0$. Substituting $s=\lambda_{0}$ in \eqref{cv1} leads to
\begin{equation}
    -r(\lambda_{1}-\lambda_{0})\chi(\lambda_{0})Z_{w}(r,\alpha_{0}\lambda_{0})=C_{1,w}(r)+\lambda_{0}C_{2,w}(r).\label{fr1}
\end{equation}
Similarly, for $s=\lambda_{1}$,
\begin{equation}
    -r(\lambda_{0}-\lambda_{1})\psi(\lambda_{1})Z_{w}(r,\alpha_{1}\lambda_{1})=C_{1,w}(r)+\lambda_{1}C_{2,w}(r).\label{fr2}
\end{equation}
Thus, $C_{1,w}(r)$, $C_{2,w}(r)$, needs to be determined, since $Z_{w}(r,\alpha_{i}\lambda_{i})$, $i=0,1$ are still unknown. This difficulty can be overcome by deriving an expression for $Z_{w}(r,\alpha_{i}\lambda_{i})$, $i=0,1$, by deriving  $Z_{w}(r,s)$ after successive iterations of \eqref{cv1}, which now becomes
\begin{equation}
    Z_{w}(r,s)=r\sum_{k=0}^{1}h_{k}(s)Z_{w}(r,a_{k}s)+L_{w}(r,s),\label{cv2}
\end{equation}
where,
\begin{equation}
    \begin{array}{rl}
        L_{w}(r,s)=&\frac{sC_{1,w}(r)+s^{2}C_{2,w}(r)}{(\lambda_{0}-s)(\lambda_{1}-s)}+re^{-sw},\vspace{2mm}\\
        h_{k}(s)=&\frac{\lambda_{k}}{\lambda_{k}-s}(\chi(s)1_{\{k=0\}}+\psi(s)1_{\{k=1\}}),\,k=0,1.
\end{array}
\end{equation}
After $n-1$ iterations, we obtain
\begin{equation}
    \begin{array}{rl}
       Z_{w}(r,s)= &r^{n}\sum_{k=0}^{n}K_{k,n-k}(s)Z_{w}(r,a_{0}^{k}a_{1}^{n-k}s)+\sum_{i=0}^{n-1}r^{i} \sum_{k=0}^{i}K_{k,i-k}(s)L_{w}(r,a_{0}^{k}a_{1}^{i-k}s),
    \end{array}
\end{equation}
where $K_{k,n-k}(s)$ are recursively defined as follows: $K_{0,0}(s)=1$, $K_{.,-1}(s)=0=K_{-1,.}(s)$, $K_{1,0}(s)=h_{0}(s)$, $K_{0,1}(s)=h_{1}(s)$ and
\begin{displaymath}
    \begin{array}{rl}
        K_{k+1,n-k}(s)=&K_{k,n-k}(s)h_{0}(a_{0}^{k}a_{1}^{n-k}s)+K_{k+1,n-k-1}(s)h_{1}(a_{0}^{k+1}a_{1}^{n-k-1}s),\,n-k\geq k+1,\vspace{2mm}  \\
         K_{k,n-k+1}(s)=&K_{k,n-k}(s)h_{1}(a_{0}^{k}a_{1}^{n-k}s)+K_{k-1,n-k+1}(s)h_{0}(a_{0}^{k-1}a_{1}^{n-k+1}s),\,n-k\leq k-1.
    \end{array}
\end{displaymath}
Therefore,
\begin{equation}
    Z_{w}(r,s)=\sum_{i=0}^{\infty}r^{i} \sum_{k=0}^{i}K_{k,i-k}(s)L_{w}(r,a_{0}^{k}a_{1}^{i-k}s)+\lim_{n\to\infty}r^{n}\sum_{k=0}^{n}K_{k,n-k}(s)Z_{w}(r,a_{0}^{k}a_{1}^{n-k}s).\label{it0}
\end{equation}
The second term in the RHS of \eqref{it0} converges to zero due to the fact that $|r|<1$, thus,
\begin{equation}
    Z_{w}(r,s)=\sum_{i=0}^{\infty}r^{i} \sum_{k=0}^{i}K_{k,i-k}(s)L_{w}(r,a_{0}^{k}a_{1}^{i-k}s).\label{it01}
\end{equation}
Setting in \eqref{it01} $s=a_{i}\lambda_{i}$ we obtain expressions for the $Z_{w}(r,a_{i}\lambda_{i})$, $i=0,1$. Note that these expressions are given in terms of the unknowns $C_{i,w}(r)$, $i=0,1$. Substituting back in \eqref{fr1}, \eqref{fr2}, we obtain a linear system of two equations with two unknowns $C_{i,w}(r)$, $i=0,1$.
\subsubsection{Stationary analysis}
Using Abel's theorem, or the relation $W_{\infty}=\left\{\begin{array}{ll}
         \left[a_{0}W_{\infty}+B-J^{(0)}\right]^{+},&\,B<T,\\
         \left[a_{1}W_{\infty}+T-J^{(1)}\right]^{+},&\,B\geq T,\end{array}\right.$ leads for $Re(s)=0$ to
         \begin{equation}
             \begin{array}{rl}
                  E(e^{-sW_{\infty}})=& \frac{\lambda_{0}}{\lambda_{0}-s}\chi(s)E(e^{-sa_{0}W_{\infty}})+\frac{\lambda_{1}}{\lambda_{1}-s}\psi(s)E(e^{-sa_{1}W_{\infty}})+1-M(s),
             \end{array}\label{lt}
         \end{equation}
         where
         \begin{displaymath}
             M(s):=E(e^{-s[a_{0}W_{\infty}+B-J^{(0)}]^{-}}1(B<T))+E(e^{-s[a_{1}W_{\infty}+T-J^{(1)}]^{-}}1(B\geq T)).
         \end{displaymath}
Setting $Z(s)=E(e^{-sW_{\infty}})$ Liouville's theorem implies that
\begin{equation}
    \begin{array}{l}
        \prod_{k=0}^{1}(\lambda_{k}-s)Z(s)\vspace{2mm}\\-(\lambda_{0}(\lambda_{1}-s)\chi(s)Z(a_{0}s)+\lambda_{1}(\lambda_{0}-s)\psi(s)Z(a_{1}s))= C_{0}+sC_{1}+s^{2}C_{2},\,Re(s)\geq 0,
\end{array}\label{ctv1}
\end{equation}
Setting $s=0$ leads to $C_{0}=0$. Moreover, for $s=\lambda_{0}$ 
\begin{equation}
    -(\lambda_{1}-\lambda_{0})\chi(\lambda_{0})Z(\alpha_{0}\lambda_{0})=C_{1}+\lambda_{0}C_{2}.\label{ftr1}
\end{equation}
Similarly, for $s=\lambda_{1}$,
\begin{equation}
    -(\lambda_{0}-\lambda_{1})\psi(\lambda_{1})Z(\alpha_{1}\lambda_{1})=C_{1}+\lambda_{1}C_{2}.\label{ftr2}
\end{equation}
Now,
\begin{equation}
    Z(s)=\sum_{k=0}^{1}h_{k}(s)Z(a_{k}s)+L(s),\label{ctv2}
\end{equation}
where, $h_{k}(s)$, $k=0,1$ as above and
\begin{equation*}
    \begin{array}{rl}
        L(s)=&\frac{sC_{1}+s^{2}C_{2}}{(\lambda_{0}-s)(\lambda_{1}-s)}.
\end{array}
\end{equation*}
After $n-1$ iterations, we obtain
\begin{equation}
    \begin{array}{rl}
       Z(s)= &\sum_{k=0}^{n}K_{k,n-k}(s)Z(a_{0}^{k}a_{1}^{n-k}s)+\sum_{i=0}^{n-1} \sum_{k=0}^{i}K_{k,i-k}(s)L(a_{0}^{k}a_{1}^{i-k}s),
    \end{array}
\end{equation}
where $K_{k,n-k}(s)$ as above. Note that $L(0)=0$ and \cite[Theorem 2]{adan} applies. Thus,
\begin{equation}
    Z(s)= \lim_{n\to\infty}\sum_{k=0}^{n}K_{k,n-k}(s)+\sum_{i=0}^{\infty} \sum_{k=0}^{i}K_{k,i-k}(s)L(a_{0}^{k}a_{1}^{i-k}s).\label{fin}
\end{equation}
Setting in \eqref{fin}, $s=a_{k}\lambda_{k}$, $k=0,1$, we obtain, respectively:
\begin{equation}
    \begin{array}{rl}
         Z(a_{0}\lambda_{0})= \lim_{n\to\infty}\sum_{k=0}^{n}K_{k,n-k}(a_{0}\lambda_{0})+\sum_{i=0}^{\infty} \sum_{k=0}^{i}K_{k,i-k}(a_{0}\lambda_{0})L(a_{0}^{k+1}a_{1}^{i-k}\lambda_{0}),\vspace{2mm}  \\
        Z(a_{1}\lambda_{1})= \lim_{n\to\infty}\sum_{k=0}^{n}K_{k,n-k}(a_{1}\lambda_{1})+\sum_{i=0}^{\infty}  \sum_{k=0}^{i}K_{k,i-k}(a_{1}\lambda_{1})L(a_{0}^{k}a_{1}^{i-k+1}\lambda_{1}). 
    \end{array}\label{ftr3}
\end{equation}
Using \eqref{ftr3} in \eqref{ftr1}, \eqref{ftr2}, we obtain a linear system of equations for $C_{1}$, $C_{2}$.

\begin{remark}
             It would be interesting to consider  the performance measures $P(W = 0)$ and $E(W^{k})$, $k = 1, 2,\ldots$ in the regime that $a_{k}\to 1$ (see also \cite[Section 2.3]{boxman}), i.e., , a perturbation of the model in \cite{boxperry}.
         \end{remark}

         Differentiating \eqref{ctv2} with respect to $s$ and letting $s=0$ yields after some algebra that,
         \begin{displaymath}
             E(W_{\infty}):=-\frac{d}{ds}Z(s)|_{s=0}=\frac{\frac{\chi(0)}{\lambda_{0}}+\frac{\psi(0)}{\lambda_{0}}+\chi^{\prime}(0)+\psi^{\prime}(0)-\frac{C_{1}+2C_{2}}{\lambda_{0}\lambda_{1}}}{1-a_{0}\chi(0)-a_{1}\psi(0)},
         \end{displaymath}
         where $f^{\prime}(.)$ denotes the derivative of a function $f(.)$ and $C_{1}$, $C_{2}$ are derived as shown above.
\subsubsection{The case $a_{0}\in(0,1)$, $a_{1}=1$}
We now consider the stationary version of the special case where $a_{1}=1$, i.e., we assume that when $B_{n}\geq T_{n}$, the next arrival does not make obsolete a fixed fraction of the already present work. This can be considered natural if we think that in such a case the service time is cut-off since exceeds the threshold $T_{n}$. Following similar arguments as above, we obtain
\begin{equation}
    \begin{array}{rl}
        Z(s)= & \frac{\lambda_{0}}{\lambda_{0}-s}\chi(s)Z(a_{0}s)+\frac{\lambda_{1}}{\lambda_{1}-s}\psi(s)Z(s)+M^{-}(s)\Leftrightarrow\vspace{2mm}\\
         (\lambda_{0}-s)Z(s)-\lambda_{0}\beta(s)Z(a_{0}s)=&(\lambda_{0}-s)\frac{\beta(s)}{\chi(s)}M^{-}(s), 
    \end{array}\label{bu}
\end{equation}
where $\beta(s):=\frac{\chi(s)}{1-\frac{\lambda_{1}\psi(s)}{\lambda_{1}-s}}$, $M^{-}(s):=1-E(e^{-s[a_{0}W_{\infty}+B-J^{(0)}]^{-}}1(B<T))-E(e^{-s[W_{\infty}+T-J^{(1)}]^{-}}1(B\geq T))$. Note that $\beta(s)$ is the LT of the distribution of the random variable $\tilde{B}$, which is the time elapsed from the epoch a service request arrives until the epoch the registered service is of threshold type:
\begin{displaymath}
    \beta(s)=E(e^{-s B}1(B\leq T))+E(e^{-s(T-J^{(1)})}1(B\geq T))\beta(s)\Leftrightarrow \beta(s)=\frac{E(e^{-s B}1(B\leq T))}{1-E(e^{-s(T-J^{(1)})}1(B\geq T))}.
\end{displaymath}
Thus, following the lines in \cite{box1}, Liouville's theorem \cite{tit} states that
\begin{equation}
    (\lambda_{0}-s)Z(s)-\lambda_{0}\beta(s)Z(a_{0}s)=C_{0}+sC_{1}.\label{bu1}
\end{equation}
For $s=0$, \eqref{bu1} implies that $C_{0}=0$. Thus,
\begin{equation}
    Z(s)=\frac{\lambda_{0}}{\lambda_{0}-s}\beta(s)Z(a_{0}s)+\frac{sC_{1}}{\lambda_{0}-s},\label{bu2}
\end{equation}
which has a solution similar to the one in \cite[Theorem 2.2]{box1}, so further details are omitted.
\subsubsection{The case $a_{0}=a_{1}:=a\in (0,1)$}
Now consider the case where the fraction of work that becomes obsolete because of an arrival, is independent on whether $B<T$ or $B\geq T$. In such a scenario,
\begin{equation}
    \begin{array}{l}
        \prod_{k=0}^{1}(\lambda_{k}-s)Z(s)\vspace{2mm}\\-[\lambda_{0}(\lambda_{1}-s)\chi(s)+\lambda_{1}(\lambda_{0}-s)\psi(s)]Z(as)= \prod_{k=0}^{1}(\lambda_{k}-s)(1-M(s)),\,Re(s)=0,
\end{array}\label{ctv11}
\end{equation}
Now,
\begin{itemize}
    \item The LHS of \eqref{ctv11} is analytic in $Re(s)>0$, and continuous in $Re(s)\geq 0$.
    \item The RHS of \eqref{ctv11} is analytic in $Re(s)<0$, and continuous in $Re(s)\leq 0$.
    \item $Z(s)$ is for $Re(s)\geq 0$ bounded by 1, and hence the LHS of \eqref{ctv11} behaves at most as a quadratic polynomial in $s$ fro large $s$, with $Re(s)> 0$.
    \item $M(s)$ is for $Re(s)\leq 0$ bounded by 1, and hence the RHS of \eqref{ctv11} behaves at most as a quadratic polynomial in $s$ fro large $s$, with $Re(s)< 0$.
\end{itemize}
Liouville's theorem \cite{tit} implies that both sides in \eqref{ctv11} are equal to the same quadratic polynomial in $s$, in their respective half-planes. Therefore,
\begin{equation}
    \prod_{k=0}^{1}(\lambda_{k}-s)Z(s)\\-[\lambda_{0}(\lambda_{1}-s)\chi(s)+\lambda_{1}(\lambda_{0}-s)\psi(s)]Z(as)=C_{0}+sC_{1}+s^{2}C_{2}.\label{sf1}
\end{equation}
Setting $s=0$ in \eqref{sf1}, and having in mind that $\chi(0)+\psi(0)=1$, we obtain $C_{0}=0$. Setting $s=\lambda_{i}$, $i=0,1$, we obtain
\begin{equation}
    \begin{array}{rl}
        -(\lambda_{1}-\lambda_{0})\chi(\lambda_{0})Z(a\lambda_{0})= &C_{1}+\lambda_{0}C_{2},\vspace{2mm}  \\
         -(\lambda_{0}-\lambda_{1})\psi(\lambda_{1})Z(a\lambda_{1})= &C_{1}+\lambda_{1}C_{2}. 
    \end{array}\label{con0}
\end{equation}
We further need to obtain $Z(a\lambda_{i})$, $i=0,1$. Note that $Z(a\lambda_{i})=P(A^{(i)}>a W)$, $i=0,1$. Now \eqref{sf1} is rewritten as
\begin{equation}
    Z(s)=H(s)Z(as)+L(s),\label{bh}
\end{equation}
where $H(s):=\frac{\lambda_{0}}{\lambda_{0}-s}\chi(s)+\frac{\lambda_{1}}{\lambda_{1}-s}\psi(s)$. Iterating \eqref{bh} and having in mind that $Z(a^{n}s)\to 1$, as $n\to\infty$, we obtain,
\begin{equation}
    Z(s)=\prod_{n=0}^{\infty}H(a^{n}s)+\sum_{n=0}^{\infty}L(a^{n}s)\prod_{j=0}^{n-1}H(a^{j}s).\label{bhh}
\end{equation}
Note that in \eqref{bhh}, $Z(s)$ appears to have singularities in $s=\lambda_{k}/a^{j}$, $j=0,1,\ldots$, $k=0,1$, but following \cite[see Remark 2.5]{box1}, it can be seen that these are removable singularities.

Setting $s=a\lambda_{0}$,
\begin{equation}\begin{array}{rl}
         Z(a\lambda_{0})=&\prod_{n=0}^{\infty}\frac{(\lambda_{1}-a^{n+1}\lambda_{0})\chi(a^{n+1}\lambda_{0})+\lambda_{1}(1-a^{n+1})\psi(a^{n+1}\lambda_{0})}{(\lambda_{1}-a^{n+1}\lambda_{0})(1-a^{n+1})}  \vspace{2mm}\\
         &+\sum_{n=0}^{\infty}\frac{a^{n+1}(C_{1}+C_{2}\lambda_{0}a^{n+1})}{(\lambda_{1}-a^{n+1}\lambda_{0})(1-a^{n+1})} \prod_{j=0}^{n-1}\frac{(\lambda_{1}-a^{j+1}\lambda_{0})\chi(a^{j+1}\lambda_{0})+\lambda_{1}(1-a^{j+1})\psi(a^{j+1}\lambda_{0})}{(\lambda_{1}-a^{j+1}\lambda_{0})(1-a^{j+1})}.
    \end{array}\label{con1}
\end{equation}
    Similarly, for $s=a\lambda_{1}$,
\begin{equation}\begin{array}{rl}
         Z(a\lambda_{1})=&\prod_{n=0}^{\infty}\frac{(\lambda_{0}-a^{n+1}\lambda_{1})\psi(a^{n+1}\lambda_{1})+\lambda_{0}(1-a^{n+1})\chi(a^{n+1}\lambda_{1})}{(\lambda_{0}-a^{n+1}\lambda_{1})(1-a^{n+1})}  \vspace{2mm}\\
         &+\sum_{n=0}^{\infty}\frac{a^{n+1}(C_{1}+C_{2}\lambda_{1}a^{n+1})}{(\lambda_{0}-a^{n+1}\lambda_{1})(1-a^{n+1})} \prod_{j=0}^{n-1}\frac{(\lambda_{0}-a^{j+1}\lambda_{1})\psi(a^{j+1}\lambda_{1})+\lambda_{0}(1-a^{j+1})\chi(a^{j+1}\lambda_{1})}{(\lambda_{0}-a^{j+1}\lambda_{1})(1-a^{j+1})}.
    \end{array}\label{con2}
\end{equation}
    Substituting \eqref{con1}, \eqref{con2} in \eqref{con0}, we obtain a linear system of equations for $C_{1}$, $C_{2}$.
\begin{remark}
Assume now that the interarrival times deterministic proportionally dependent to service times. More precisely, let $J^{(k)}_{n}=c_{k}U^{(k)}_{n}+X^{(k)}_{n}$, $k=0,1$, $c_{k}\in (0,1)$ where $U^{(0)}_{n}:=B_{n}$, $U^{(1)}_{n}:=T_{n}$, and $X^{(k)}_{n}\sim exp(\delta_{k})$. Thus,
\begin{equation}
    W_{n+1}=\left\{\begin{array}{ll}
         \left[a_{0}W_{n}+(1-c_{0})B_{n}-X^{(0)}_{n}\right]^{+},&\,B_{n}<T_{n},\\
         \left[a_{1}W_{n}+(1-c_{1})T_{n}-X^{(1)}_{n}\right]^{+},&\,B_{n}\geq T_{n},
    \end{array}\right.
\end{equation}
Following similar arguments as in the previous section, we arrive for $Re(s)=0$, at
\begin{equation}
    Z_{w}(r,s)-re^{-sw}=r\frac{\delta_{0}}{\delta_{0}-s}\chi(s(1-c_{0}))Z_{w}(r,a_{0}s)+r\frac{\delta_{1}}{\delta_{1}-s}\psi(s(1-c_{1}))Z_{w}(r,a_{1}s)+\frac{r^{2}}{1-r}-rM_{w}(r,s),\label{y2}
\end{equation}
where now $M_{w}(r,s)=\sum_{n=1}^{\infty}r^{n}U^{-}_{w,n}(s)$ with
\begin{displaymath}
    U^{-}_{w,n}(s):=E_{w}(e^{-s[a_{0}W_{n}+(1-c_{0})B_{n}-X^{(0)}_{n}]^{-}}1(B_{n}<T_{n}))+E_{w}(e^{-s[a_{1}W_{n}+(1-c_{1})T_{n}-X^{(1)}_{n}]^{-}}1(B_{n}\geq T_{n})).
\end{displaymath}
Using similar arguments as above, Liouville's theorem \cite{tit} implies that
\begin{equation}
    \begin{array}{r}
        \prod_{k=0}^{1}(\delta_{k}-s)(Z_{w}(r,s)-re^{-sw})-r(\delta_{0}(\delta_{1}-s)\chi(s(1-c_{0}))Z_{w}(r,a_{0}s)+\delta_{1}(\delta_{0}-s)\psi(s(1-c_{1}))Z_{w}(r,a_{1}s))\vspace{2mm}\\ =C_{0,w}(r)+sC_{1,w}(r)+s^{2}C_{2,w}(r),\,Re(s)\geq 0.
\end{array}\label{cv1b}
\end{equation}

The rest of the analysis follows as the one in the previous section. Similar steps as those in the previous section can be followed to cope with the stationary analysis, so further details are omitted.
\end{remark}
\subsection{Interarrival times random proportionally dependent to service times}
Assume that $J^{(k)}_{n}=G^{(k)}_{n}U^{(k)}_{n}+X^{(k)}_{n}$, $k=0,1$, where $U^{(0)}_{n}:=B_{n}$, $U^{(1)}_{n}:=T_{n}$, and $X_{n}^{(k)}$ i.i.d. random variables with distribution that have rational LST:
\begin{displaymath}
    \phi_{X_{k}}(s)=\frac{N_{k}(s)}{D_{k}(s)},
\end{displaymath}
where $D_{k}(s):=\prod_{i=1}^{L_{k}}(s+t_{i}^{(k)})$ with $N_{k}(s)$ is a polynomial at most $L_{k}-1$, not sharing zeros with $D_{k}(s),$ $k=0,1$. Moreover, assume that $Re(t_{i}^{(k)})>0$, $i=1,\ldots,L_{k}$. Thus,
\begin{equation}
    W_{n+1}=\left\{\begin{array}{ll}
         \left[a_{0}W_{n}+(1-G^{(0)}_{n})B_{n}-X^{(0)}_{n}\right]^{+},&\,B_{n}<T_{n},\\
         \left[a_{1}W_{n}+(1-G^{(1)}_{n})T_{n}-X^{(1)}_{n}\right]^{+},&\,B_{n}\geq T_{n},
    \end{array}\right.
\end{equation}
where $P(G_{n}^{(0)}=\beta_{i})=p_{i}$, $i=1,\ldots,K$, $P(G_{n}^{(1)}=\gamma_{i})=q_{i}$, $i=1,\ldots,M$. Assume that $\beta_{i}\in(0,1)$, $i=1,\ldots,K$, $\gamma_{i}\in(0,1)$, $i=1,\ldots,M$. Following similar arguments as in the previous section, we arrive for $Re(s)=0$, at
\begin{equation}
\begin{array}{rl}
     Z_{w}(r,s)-re^{-sw}=& r\frac{N_{0}(-s)}{D_{0}(-s)}\sum_{i=1}^{K}p_{i}\chi(s(1-\beta_{i}))Z_{w}(r,a_{0}s)+r\frac{N_{1}(-s)}{D_{1}(-s)}\sum_{i=1}^{M}q_{i}\psi(s(1-\gamma_{i}))Z_{w}(r,a_{1}s)\vspace{2mm}\\
     &+\frac{r^{2}}{1-r}-rM_{w}(r,s), 
\end{array}
    \label{y3}
\end{equation}
where now $M_{w}(r,s)=\sum_{n=1}^{\infty}r^{n}U^{-}_{w,n}(s)$ with
\begin{displaymath}
    U^{-}_{w,n}(s):=E_{w}(e^{-s[a_{0}W_{n}+(1-G_{n}^{(0)})B_{n}-X^{(0)}_{n}]^{-}}1(B_{n}<T_{n}))+E_{w}(e^{-s[a_{1}W_{n}+(1-G_{n}^{(1)})T_{n}-X^{(1)}_{n}]^{-}}1(B_{n}\geq T_{n})).
\end{displaymath}
Then, for $Re(s)=0$,
\begin{equation}
    \begin{array}{l}
        D_{0}(-s)D_{1}(-s)[ Z_{w}(r,s)-re^{-sw}]-rN_{0}(-s)D_{1}(-s)\sum_{i=1}^{K}p_{i}\chi(s(1-\beta_{i}))Z_{w}(r,a_{0}s)\vspace{2mm}\\
        -rN_{1}(-s)D_{0}(-s)\sum_{i=1}^{M}q_{i}\psi(s(1-\gamma_{i}))Z_{w}(r,a_{1}s)=D_{0}(-s)D_{1}(-s)[\frac{r^{2}}{1-r}-rM_{w}(r,s)].
    \end{array}\label{eqw}
\end{equation}
Now,
\begin{itemize}
    \item The LHS of \eqref{eqw} is analytic in $Re(s)>0$, and continuous in $Re(s)\geq 0$.
    \item The RHS of \eqref{eqw} is analytic in $Re(s)<0$, and continuous in $Re(s)\leq 0$.
    \item For large $s$, both sides in \eqref{eqw} are $O(s^{L_{0}+L_{1}})$ in their respective half-planes.
\end{itemize}
Thus, Liouville's theorem implies that for $Re(s)\geq 0$,
\begin{equation}
\begin{array}{l}
    D_{0}(-s)D_{1}(-s)[ Z_{w}(r,s)-re^{-sw}]-rN_{0}(-s)D_{1}(-s)\sum_{i=1}^{K}p_{i}\chi(s(1-\beta_{i}))Z_{w}(r,a_{0}s)\vspace{2mm}\\
    -rN_{1}(-s)D_{0}(-s)\sum_{i=1}^{M}q_{i}\psi(s(1-\gamma_{i}))Z_{w}(r,a_{1}s)=\sum_{l=0}^{L_{1}+L_{2}}c_{l}(r)s^{l},\end{array}\label{qq1}
\end{equation}
and for $Re(s)\leq 0$,
\begin{displaymath}
    D_{0}(-s)D_{1}(-s)[\frac{r^{2}}{1-r}-rM_{w}(r,s)]=\sum_{l=0}^{L_{1}+L_{2}}c_{l}(r)s^{l}.
\end{displaymath}
For $s=0$, \eqref{qq1} implies after simple computations that $c_{0}(r)=0$. For $s=t_{j}^{(0)}$, $j=1,\ldots,L_{0}$, \eqref{qq1} implies that
\begin{equation}
    -rN_{0}(-t_{j}^{(0)})D_{1}(-t_{j}^{(0)})\sum_{i=1}^{K}p_{i}\chi(t_{j}^{(0)}(1-\beta_{i}))Z_{w}(r,a_{0}t_{j}^{(0)})=\sum_{l=1}^{L_{1}+L_{2}}c_{l}(r)(t_{j}^{(0)})^{l}.\label{cva1}
\end{equation}
Similarly, for $s=t_{j}^{(1)}$, $j=1,\ldots,L_{1}$, we have,
\begin{equation}
    -rN_{1}(-t_{j}^{(1)})D_{0}(-t_{j}^{(1)})\sum_{i=1}^{M}q_{i}\psi(t_{j}^{(1)}(1-\gamma_{i}))Z_{w}(r,a_{1}t_{j}^{(1)})=\sum_{l=0}^{L_{1}+L_{2}}c_{l}(r)(t_{j}^{(1)})^{l},\label{cva2}
\end{equation}
Note that $N_{k}(-t_{j}^{(k)})\neq 0$, $k=0,1$. Then, \eqref{cva1}, \eqref{cva2} constitutes a system of equations to obtain the remaining of the coefficients $c_{l}(r)$, $l=1,\ldots,L_{0}+L_{1}$. However, we still need to obtain $Z_{w}(r,a_{1}t_{j}^{(k)})$, $k=0,1$, $j=1,\ldots,L_{k}$. Note that \eqref{qq1} has the same form as in \eqref{cv2} but now,
\begin{equation}
    \begin{array}{rl}
        L_{w}(r,s)=&\frac{\sum_{l=1}^{L_{1}+L_{2}}s^{l}c_{i}(r)}{D_{0}(-s)D_{1}(-s)}+re^{-sw},\vspace{2mm}\\
        h_{k}(s)=&\frac{N_{k}(-s)}{D_{k}(-s)}(\sum_{i=1}^{K}p_{i}\chi(s(1-\beta_{i}))1_{\{k=0\}}+\sum_{i=1}^{M}q_{i}\psi(s(1-\gamma_{i}))1_{\{k=1\}}),\,k=0,1.
\end{array}\label{bnhy}
\end{equation}
Thus, the expression for $Z_{w}(r,s)$ is the same as in \eqref{it01}, where the expressions $K_{k,i-k}(s)$, $L_{w}(r,a_{0}^{k}a_{1}^{i-k}s)$ are obtained analogously using \eqref{bnhy}. Having this expression, we can obtain  $Z_{w}(r,a_{1}t_{j}^{(k)})$, $k=0,1$, $j=1,\ldots,L_{k}$, and substituting back in \eqref{cva1}, \eqref{cva2}, we can derive the remaining coefficients $c_{l}(r)$, $l=1,\ldots,L_{0}+L_{1}$.

\subsubsection{A more general case}
We now consider the case where the interarrival times are also depended on the system time. More precisely, we assume that $J^{(k)}_{n}=G_{n}^{(k)}(U^{(k)}_{n}+W_{n})+X^{(k)}_{n}$. 
\begin{equation}
    W_{n+1}=\left\{\begin{array}{ll}
         \left[(1-G_{n}^{(0)})W_{n}+(1-G_{n}^{(0)})B_{n}-X^{(0)}_{n}\right]^{+},&\,B_{n}<T_{n},\\
         \left[(1-G_{n}^{(1)})W_{n}+(1-G_{n}^{(1)})T_{n}-X^{(1)}_{n}\right]^{+},&\,B_{n}\geq T_{n},
    \end{array}\right.
\end{equation}
Thus,
\begin{equation}
    \begin{array}{rl}
        E_{w}(e^{-sW_{n+1}}) = & E_{w}(e^{-s[(1-G_{n}^{(0)})W_{n}+(1-G_{n}^{(0)})B_{n}-X^{(0)}_{n}]^{+}}1(B_{n}<T_{n}))\\&+ E_{w}(e^{-s[(1-G_{n}^{(1)})W_{n}+(1-G_{n}^{(1)})T_{n}-X^{(1)}_{n}]^{+}}1(B_{n}\geq T_{n}))\vspace{2mm}\\
         =&E(e^{sJ^{(0)}_{n}})\sum_{i=1}^{K}p_{i}E_{w}(e^{-s\bar{\beta}_{i}W_{n}})E(e^{-s\bar{\beta}_{i}B_{n}}1(B_{n}<T_{n}))\vspace{2mm}\\&+E(e^{sJ^{(1)}_{n}})\sum_{i=1}^{M}q_{i}E_{w}(e^{-s\bar{\gamma}_{i}W_{n}})E(e^{-s\bar{\gamma}_{i}T_{n}}1(B_{n}\geq T_{n}))+1-U^{-}_{w,n}(s),
    \end{array}\label{op1}
\end{equation}
Then, by using \eqref{op1} and using similar arguments as above we obtain for $Re(s)=0$,
\begin{equation*}
    \begin{array}{r}
        Z_{w}(r,s)-re^{-sw}=r\phi_{X_{0}}(-s)\sum_{i=1}^{K}p_{i}\chi(s\bar{\beta}_{i})Z_{w}(r,s\bar{\beta}_{i})+r\phi_{X_{1}}(-s)\sum_{i=1}^{M}q_{i}\psi(s\bar{\gamma}_{i})Z_{w}(r,s\bar{\gamma}_{i})\vspace{2mm}\\
        +\frac{r^{2}}{1-r}-rM_{w}(r,s),
    \end{array}\label{eee}
\end{equation*}
or equivalently,
\begin{equation}
    \begin{array}{l}
        D_{0}(-s)D_{1}(-s)[ Z_{w}(r,s)-re^{-sw}]-rN_{0}(-s)D_{1}(-s)\sum_{i=1}^{K}p_{i}\chi(s\bar{\beta}_{i})Z_{w}(r,s\bar{\beta}_{i})\vspace{2mm}\\
        -rN_{1}(-s)D_{0}(-s)\sum_{i=1}^{M}q_{i}\psi(s\bar{\gamma}_{i})Z_{w}(r,s\bar{\gamma}_{i})=D_{0}(-s)D_{1}(-s)[\frac{r^{2}}{1-r}-rM_{w}(r,s)].
    \end{array}\label{eqw1}
\end{equation}
Now,
\begin{itemize}
    \item The LHS of \eqref{eqw1} is analytic in $Re(s)>0$, and continuous in $Re(s)\geq 0$.
    \item The RHS of \eqref{eqw1} is analytic in $Re(s)<0$, and continuous in $Re(s)\leq 0$.
    \item For large $s$, both sides in \eqref{eqw} are $O(s^{L_{0}+L_{1}})$ in their respective half-planes.
\end{itemize}
Thus, Liouville's theorem implies that for $Re(s)\geq 0$,
\begin{equation}
\begin{array}{l}
    D_{0}(-s)D_{1}(-s)[ Z_{w}(r,s)-re^{-sw}]-rN_{0}(-s)D_{1}(-s)\sum_{i=1}^{K}p_{i}\chi(s\bar{\beta}_{i})Z_{w}(r,s\bar{\beta}_{i})\vspace{2mm}\\
        -rN_{1}(-s)D_{0}(-s)\sum_{i=1}^{M}q_{i}\psi(s\bar{\gamma}_{i})Z_{w}(r,s\bar{\gamma}_{i})=\sum_{l=0}^{L_{1}+L_{2}}c_{l}(r)s^{l},\end{array}\label{qq11}
\end{equation}
For $s=0$, $c_{0}(r)=0$. For convenience, set $\bar{\beta}_{i}=a_{i}$, $i=1,\ldots,K$, $\bar{\gamma}_{i}=a_{K+i}$, $q_{i}=p_{K+i}$, $i=1,\ldots,M$, and
\begin{displaymath}
    f(a_{i}s):=\left\{\begin{array}{ll}
         \phi_{X_{0}}(-s)\chi(sa_{i}),&i=1,\ldots,K,  \\
         \phi_{X_{1}}(-s)\psi(sa_{i}),&i=K+1,\ldots,K+M.  
    \end{array}\right.
\end{displaymath}
Then, \eqref{qq11} can be written as
\begin{equation}
    Z_{w}(r,s)=r\sum_{i=1}^{K+M}p_{i}f(a_{i}s)Z_{w}(r,a_{i}s)+L_{w}(r,s),
\end{equation}
where $L_{w}(r,s):=\frac{\sum_{l=1}^{L_{1}+L_{2}}s^{l}c_{i}(r)}{D_{0}(-s)D_{1}(-s)}+re^{-sw}$. Therefore,
\begin{equation}
\begin{array}{l}
    Z_{w}(r,s)=\sum_{i=0}^{\infty}r^{i} \sum_{i_{1}+\ldots+i_{K+M}=i}p_{1}^{i_{1}}\ldots p_{K+M}^{i_{K+M}}L_{i_{1},\ldots,i_{K+M}}(s)\vspace{2mm}\\+\lim_{n\to\infty}r^{n}\sum_{i_{1}+\ldots+i_{K+M}=n}p_{1}^{i_{1}}\ldots p_{K+M}^{i_{K+M}}L_{i_{1},\ldots,i_{K+M}}(s)Z_{w}(r,a_{1}^{i_{1}}\ldots a_{K+M}^{i_{K+M}}s),\end{array}\label{ait0}
\end{equation}
where $L_{0,0,\ldots,0,1,0,\ldots,0}(s):=f(a_{k}s)$, with 1 in position $k$, and $k=1,\ldots, K+M$,
\begin{displaymath}
    L_{i_{1},\ldots,i_{K+M}}(s)=f(a_{1}^{i_{1}}\ldots a_{K+M}^{i_{K+M}}s)\sum_{j=1}^{K+M}L_{i_{1},\ldots,i_{j}-1,\ldots,i_{K+M}}(s).
\end{displaymath}
The second term in the RHS of \eqref{ait0} converges to zero due to the fact that $|r|<1$, thus,
\begin{equation}
    Z_{w}(r,s)=\sum_{i=0}^{\infty}r^{i} \sum_{i_{1}+\ldots+i_{K+M}=i}p_{1}^{i_{1}}\ldots p_{K+M}^{i_{K+M}}L_{i_{1},\ldots,i_{K+M}}(s).\label{ait01}
\end{equation}
Setting $s=t_{j}^{(k)}$, $j=1,\ldots,L_{k}$, $k=0,1$ in \eqref{qq11} we obtain a system of equations for the remaining coefficients $c_{l}(r)$, $l=1,\ldots,L_{0}+L_{1}$. Specifically for $s=t_{j}^{(0)}$, $j=1,\ldots,L_{0}$
\begin{equation}
    -rN_{0}(-t_{j}^{(0)})D_{1}(-t_{j}^{(0)})\sum_{i=1}^{K}p_{i}\chi(t_{j}^{(0)}\bar{\beta}_{i})Z_{w}(r,\bar{\beta}_{i}t_{j}^{(0)})=\sum_{l=1}^{L_{1}+L_{2}}c_{l}(r)(t_{j}^{(0)})^{l},\label{cva11}
\end{equation}
and for $s=t_{j}^{(1)}$, $j=1,\ldots,L_{1}$, we have,
\begin{equation}
    -rN_{1}(-t_{j}^{(1)})D_{0}(-t_{j}^{(1)})\sum_{i=1}^{M}q_{i}\psi(t_{j}^{(1)}\bar{\gamma}_{i})Z_{w}(r,\bar{\gamma}_{i}t_{j}^{(1)})=\sum_{l=0}^{L_{1}+L_{2}}c_{l}(r)(t_{j}^{(1)})^{l},\label{cva22}
\end{equation}
where we have further use the expression in \eqref{ait01}.

\subsection{A mixed case}
Consider the following recursion:
\begin{equation}
    W_{n+1}=\left\{\begin{array}{ll}
         \left[aW_{n}+(1-G^{(0)}_{n})B_{n}-X^{(0)}_{n}\right]^{+},&\,B_{n}<T_{n},\\
         \left[V_{n}W_{n}+(1-G^{(1)}_{n})T_{n}-X^{(1)}_{n}\right]^{+},&\,B_{n}\geq T_{n},
    \end{array}\right.
\end{equation}
where $V_{n}<0$, and $a\in(0,1)$. Then, following a similar procedure we obtain
\begin{displaymath}
    \begin{array}{l}
         Z_{w}(r,s)-re^{-sw}=Z_{w}(r,as)r\frac{\delta_{0}}{\delta_{0}-s}\sum_{i=1}^{K}p_{i}\chi(s\bar{\beta}_{i})+\int_{-\infty}^{0}Z_{w}(r,sy)P(V\in dy)r\frac{\delta_{1}}{\delta_{1}-s}\sum_{i=1}^{M}q_{i}\psi(s\bar{\gamma}_{i})+\frac{r^{2}}{1-r}-rM_{w}(r,s),
    \end{array}
\end{displaymath}
$M_{w}(r,s)=\sum_{n=1}^{\infty}r^{n}U^{-}_{w,n}(s)$ with
\begin{displaymath}
    U^{-}_{w,n}(s):=E_{w}(e^{-s[aW_{n}+(1-G_{n}^{(0)})B_{n}-X^{(0)}_{n}]^{-}}1(B_{n}<T_{n}))+E_{w}(e^{-s[V_{n}W_{n}+(1-G_{n}^{(1)})T_{n}-X^{(1)}_{n}]^{-}}1(B_{n}\geq T_{n})).
\end{displaymath}
Equivalently, we have
\begin{equation}
    \begin{array}{l}
       \prod_{j=0}^{1}(\delta_{j}-s)(Z_{w}(r,s)-re^{-sw})-Z_{w}(r,as)r\delta_{0}(\delta_{1}-s)\sum_{i=1}^{K}p_{i}\chi(s\bar{\beta}_{i})   \vspace{2mm}\\
         =\int_{-\infty}^{0}Z_{w}(r,sy)P(V\in dy)r\delta_{1}(\delta_{0}-s)\sum_{i=1}^{M}q_{i}\psi(s\bar{\gamma}_{i})+\prod_{j=0}^{1}(\delta_{j}-s)(\frac{r^{2}}{1-r}-rM_{w}(r,s)).
    \end{array}
    \label{vby}
\end{equation}
Clearly,
\begin{itemize}
     \item the LHS of \eqref{vby} is analytic in $Re(s)>0$ and continuous in $Re(s)\geq 0$,
    \item the RHS of \eqref{vby} is analytic in $Re(s)<0$ and continuous in $Re(s)\leq 0$,
    \item for large $s$, both sides are $O(s^{2})$ in their respective half planes. 
\end{itemize}
Thus, Liouville's theorem now states that 
\begin{displaymath}
    \prod_{j=0}^{1}(\delta_{j}-s)(Z_{w}(r,s)-re^{-sw})-Z_{w}(r,as)r\delta_{0}(\delta_{1}-s)\sum_{i=1}^{K}p_{i}\chi(s\bar{\beta}_{i}) =c_{0}+c_{1}s+c_{2}s^{2},\,Re(s)\geq 0.
\end{displaymath}
For $s=0$, we have $c_{0}=\frac{r\delta_{0}\delta_{1}}{1-r}(1-\chi(0))$. Setting $s=\delta_{1}$, and $s=\delta_{0}$, we respectively have the following linear system:
\begin{displaymath}
    \begin{array}{rl}
        c_{2}\delta_{1}^{2}+c_{1}\delta_{1}= &-c_{0},\vspace{2mm}  \\
         c_{2}\delta_{0}^{2}+c_{1}\delta_{0}= &-c_{0}-r\delta_{0}(\delta_{1}-\delta_{0})\chi(0)Z_{w}(r,a\delta_{0}),
    \end{array}
\end{displaymath}
from which,
\begin{equation}
    \begin{array}{rl}
        c_{1}=&-\frac{r}{1-r}((\delta_{0}+\delta_{1})(1-\chi(0))+\delta_{1}\chi(\delta_{0})(1-r)Z_{w}(r,a\delta_{0})), \vspace{2mm} \\
         c_{2}=&\frac{r}{1-r}(1-\chi(0)+\chi(\delta_{0})(1-r)Z_{w}(r,a\delta_{0})).
    \end{array}\label{cxa}
\end{equation}

It remains to find $Z_{w}(r,a\delta_{0})$. This can be done by iteratively solving 
\begin{displaymath}
    Z_{w}(r,s)=r\frac{\delta_{0}}{\delta_{0}-s}\chi(s)\sum_{i=1}^{K}p_{i}\chi(s\bar{\beta}_{i})Z_{w}(r,as)+\frac{c_{0}+sc_{1}+s^{2}c_{2}}{ \prod_{j=0}^{1}(\delta_{j}-s)}+re^{-sw}.
\end{displaymath}
In particular, 
\begin{equation}
    Z_{w}(r,s)=\sum_{n=0}^{\infty}L_{w}(r,a^{n}s)\prod_{j=0}^{n-1}K_{w}(r,a^{j}s),
    \label{solt1s}
\end{equation}
where $L_{w}(r,s):=\frac{c_{0}+sc_{1}+s^{2}c_{2}}{ \prod_{j=0}^{1}(\delta_{j}-s)}+re^{-sw}$, $K_{w}(r,s):=r\frac{\delta_{0}}{\delta_{0}-s}\chi(s)\sum_{i=1}^{K}p_{i}\chi(s\bar{\beta}_{i})$. Thus,
\begin{equation}
    Z_{w}(r,a\delta_{0})=\sum_{n=0}^{\infty}L_{w}(r,a^{n}\delta_{0})\prod_{j=0}^{n-1}K_{w}(r,a^{j}\delta_{0}).
    \label{solt1cs}
\end{equation}
Substituting \eqref{solt1cs} in \eqref{cxa} we obtain a linear system for $c_{1}$, $c_{2}$.
\subsubsection{A more general case}
Assume the case where the Laplace-Stieltjes transform of the distribution of $X_{n}^{(k)}$ are rational and such that:
\begin{displaymath}
    A_{0}(s)=\frac{\widehat{A}_{0}(s)}{\prod_{i=1}^{L_{0}}(s+\delta_{i})},\,\,A_{1}(s)=\frac{\widehat{A}_{1}(s)}{\prod_{i=1}^{L_{1}}(s+\zeta_{i})},
\end{displaymath}
with $\widehat{A}_{k}(s)$ is a polynomial at most $L_{k}-1$, not sharing zeros with the corresponding denominators of $A_{k}(s),$ $k=0,1$. Moreover, assume that $Re(\delta_{i})>0$, $i=1,\ldots,L_{0}$, and $Re(\zeta_{i})<0$, $i=1,\ldots,L_{1}$. Moreover, assume that $\gamma_{i}>1$, $i=1,\ldots,M$. Then, \eqref{vby} becomes now for $Re(s)=0$:
\begin{equation}
    \begin{array}{l}
       \prod_{j=1}^{L_{0}}(\delta_{j}-s)\prod_{j=1}^{L_{1}}(\zeta_{j}-s)(Z_{w}(r,s)-re^{-sw})-Z_{w}(r,as)r\widehat{A}_{0}(-s)\prod_{j=1}^{L_{1}}(\zeta_{j}-s)\sum_{i=1}^{K}p_{i}\chi(s\bar{\beta}_{i})   \vspace{2mm}\\
         =\int_{-\infty}^{0}Z_{w}(r,sy)P(V\in dy)r\widehat{A}_{1}(-s)\prod_{j=1}^{L_{0}}(\delta_{j}-s)\sum_{i=1}^{M}q_{i}\psi(s\bar{\gamma}_{i})+\prod_{j=1}^{L_{0}}(\delta_{j}-s)\prod_{j=1}^{L_{1}}(\zeta_{j}-s)(\frac{r^{2}}{1-r}-rM_{w}(r,s)).
    \end{array}
    \label{vby1}
\end{equation}
Again,
\begin{itemize}
     \item the LHS of \eqref{vby1} is analytic in $Re(s)>0$ and continuous in $Re(s)\geq 0$,
    \item the RHS of \eqref{vby1} is analytic in $Re(s)<0$ and continuous in $Re(s)\leq 0$,
    \item for large $s$, both sides are $O(s^{L_{0}+L_{1}})$ in their respective half planes. 
\end{itemize}
Thus, Liouville's theorem now states that for $Re(s)\geq 0$,
\begin{equation}
  \prod_{j=1}^{L_{0}}(\delta_{j}-s)\prod_{j=1}^{L_{1}}(\zeta_{j}-s)(Z_{w}(r,s)-re^{-sw})-Z_{w}(r,as)r\widehat{A}_{0}(-s)\prod_{j=1}^{L_{1}}(\zeta_{j}-s)\sum_{i=1}^{K}p_{i}\chi(s\bar{\beta}_{i})=\sum_{j=0}^{L_{0}+L_{1}}c_{j}(r)s^{j},
  \label{pop1}
\end{equation}
and for $Re(s)\leq 0$,
\begin{equation}
\int_{-\infty}^{0}Z_{w}(r,sy)P(V\in dy)r\widehat{A}_{1}(-s)\prod_{j=1}^{L_{0}}(\delta_{j}-s)\sum_{i=1}^{M}q_{i}\psi(s\bar{\gamma}_{i})+\prod_{j=1}^{L_{0}}(\delta_{j}-s)\prod_{j=1}^{L_{1}}(\zeta_{j}-s)(\frac{r^{2}}{1-r}-rM_{w}(r,s))=\sum_{j=0}^{L_{0}+L_{1}}c_{j}(r)s^{j}.  \label{pop2}  
\end{equation}
Setting $s=0$, and using either \eqref{pop1}, or \eqref{pop2} we get after straightforward computations that
\begin{displaymath}
    c_{0}(r)=\frac{r^{2}(1-\chi(0))}{1-r}\prod_{j=1}^{L_{0}}\delta_{j}\prod_{j=1}^{L_{1}}\zeta_{j}.
\end{displaymath}

For $s=\delta_{j}$, $j=1,\ldots,L_{0}$, \eqref{pop1} gives,
\begin{displaymath}
    \sum_{k=1}^{L_{0}+L_{1}}c_{k}(r)\delta_{j}^{k}=-r\widehat{A}_{0}(-\delta_{j})\prod_{m=1}^{L_{1}}(\zeta_{m}-\delta_{j})\sum_{i=1}^{K}p_{i}\chi(\delta_{j}\bar{\beta}_{i})Z_{w}(r,a\delta_{j}).
\end{displaymath}
Note that we further need other $L_{1}$ equations to obtain all $c_{m}(r)$. Note that for $s=\zeta_{j}$, $j=1,\ldots,L_{1}$, \eqref{pop2} gives:
\begin{displaymath}
    \sum_{k=1}^{L_{0}+L_{1}}c_{k}(r)\zeta_{j}^{k}=-r\widehat{A}_{1}(-\zeta_{j})\prod_{m=1}^{L_{0}}(\delta_{m}-\zeta_{j})\sum_{i=1}^{M}q_{i}\psi(\zeta_{j}\bar{\gamma}_{i})\int_{-\infty}^{0}Z_{w}(r,sy)P(V\in dy).
\end{displaymath}
Note that \eqref{pop1} can be rewritten as
\begin{displaymath}
    Z_{w}(r,s)=K(r,s)Z_{w}(r,as)+L_{w}(r,s),
\end{displaymath}
with
\begin{displaymath}
    K(r,s)=r\frac{A_{0}(s)}{\prod_{m=1}^{L_{1}}(\zeta_{m}-s)}\sum_{i=1}^{K}p_{i}\chi(s\bar{\beta}_{i}),\,\,\,L_{w}(r,s)=\frac{\sum_{k=0}^{L_{0}+L_{1}}c_{k}(r)s^{k}}{\prod_{m=1}^{L_{1}}(\zeta_{m}-s)\prod_{m=1}^{L_{0}}(\delta_{m}-s)}+re^{-sw}.
\end{displaymath}
Iteration implies that
\begin{equation}
    Z_{w}(r,s)=\sum_{n=0}^{\infty}L_{w}(r,a^{n}s)\prod_{m=0}^{n-1}K(r,a^{m}s),\label{ale}
\end{equation}
where the convergence of the infinite sum can be proven with the aid of D'Alembert's test, since $a\in(0,1)$, and
\begin{displaymath}
    \lim_{k\to\infty}|\frac{L_{w}(r,a^{k}s)}{L_{w}(r,a^{k+1}s)K(r,a^{k}s)}|=|\frac{\prod_{j=1}^{L_{1}}\zeta_{j}}{r\chi(0)}|.
\end{displaymath}
Setting $s=a\delta_{j}$, $j=1,\ldots,L_{0}$ in \eqref{ale} we obtain $Z_{w}(r,a\delta_{j})$.

\section{The uniform proportional case with dependence}
In the following, we consider recursions of the form
\begin{equation}
     W_{n+1}=[V_{n}W_{n}+B_{n}-A_{n}]^{+},\label{ddd}
\end{equation}
with $V_{n}\sim U(0,1)$, and dependence among $B_{n}$, $A_{n}$.
\subsection{Deterministic proportional dependency with additive and subtracting delay}
    We consider the case where
\begin{displaymath}
    W_{n+1}=[V_{n}W_{n}+B_{n}-A_{n}]^{+},
\end{displaymath}
with $V_{n}\sim U(0,1)$ and for $c_{0},c_{1}\in (0,1)$, $\tilde{J}_{n}\sim exp(\delta)$, $\widehat{J}_{n}\sim exp(\nu)$:
\begin{displaymath}
    A_{n}=\left\{\begin{array}{ll}
         A_{n}^{(0)}:=c_{0}B_{n}+\tilde{J}_{n},& \text{ w.p. }p,  \\
        A_{n}^{(1)}:=[c_{1}B_{n}-\widehat{J}_{n}]^{+}, & \text{ w.p. }q:=1-p,
    \end{array}\right.
\end{displaymath}
The case of \textit{independent} $\{A_{i}\}$, $\{B_{i}\}$ was treated in \cite{box2}. Stability is ensured when $E(log|V|)<0$.
Note that
\begin{displaymath}
    \begin{array}{rl}
      E(e^{-sA_{n}^{(0)}}|B_{n}=t)=   &\frac{\delta}{\delta-s}e^{-sc_{0}t},\vspace{2mm}  \\
         E(e^{-sA_{n}^{(1)}}|B_{n}=t)=   &\frac{\nu e^{-sc_{1}t}-se^{-\nu c_{1}t}}{\nu+s},
    \end{array}
\end{displaymath}
so that
\begin{displaymath}
    \begin{array}{rl}
        E(e^{-sA_{n}^{(0)}-z B_{n}})=   &\frac{\delta}{\delta+s}\phi_{B}(z+sc_{0}),\,Re(z+sc_{0})>0,  \vspace{2mm}\\
         E(e^{-sA_{n}^{(1)}-zB_{n}})=   &\frac{\nu \phi_{B}(z+sc_{1})-s\phi_{B}(z+\nu c_{1})}{\nu-s},\,Re(z+sc_{1})>0.
    \end{array}
\end{displaymath}
Then,
\begin{displaymath}\begin{array}{rl}
     Z_{n+1}(s)=&E(e^{-sW_{n+1}})=E(e^{-s[V_{n}W_{n}+B_{n}-A_{n}]^{+}})  \vspace{2mm}\\
     =&pE(e^{-s[V_{n}W_{n}+B_{n}-A_{n}^{(0)}]^{+}})+qE(e^{-s[V_{n}W_{n}+B_{n}-A_{n}^{(1)}]^{+}}) 
\end{array}
\end{displaymath}
Since for $c_{1}\in(0,1)$, 
\begin{displaymath}\begin{array}{rl}
    [V_{n}W_{n}+B_{n}-A_{n}^{(1)}]^{+}&=[V_{n}W_{n}+B_{n}-[c_{1}B_{n}-\widehat{J}_{n}]^{+}]^{+}\vspace{2mm}\\&=V_{n}W_{n}+B_{n}-[cB_{n}-\widehat{J}_{n}]^{+}.
    \end{array}
\end{displaymath}
    Therefore, for $n\in \mathbb{N}$:
 {\small{\begin{displaymath}
        \begin{array}{rl}
            Z_{n+1}(s)=&p\left(E(e^{-sV_{n}W_{n}})E(e^{-sB_{n}+sA_{n}^{(0)}})+1-E(e^{-s[V_{n}W_{n}+B_{n}-A_{n}^{(0)}]^{-}})\right) \vspace{2mm} \\
             &+q E(e^{-sV_{n}W_{n}}) E(e^{-sB_{n}+sA_{n}^{(1)}})\vspace{2mm}\\
             =&E(e^{-sV_{n}W_{n}})\left(p\frac{\delta}{\delta-s}\phi_{B}(s\bar{c}_{0})+q\frac{\nu \phi_{B}(s\bar{c}_{1})+s\phi_{B}(s+\nu c_{1})}{\nu+s}\right)\vspace{2mm}\\
             &+p\left(1-\left[P(V_{n}W_{n}+B_{n}-A_{n}^{(0)}\geq 0)\right.\right.\vspace{2mm}\\
             &\left.\left.+P(V_{n}W_{n}+B_{n}-A_{n}^{(0)}< 0)\frac{\delta}{\delta-s}\right]\right)\vspace{2mm}\\
             =&\frac{1}{s}\int_{0}^{s}Z_{n}(y)dy\left(p\frac{\delta}{\delta-s}\phi_{B}(s\bar{c}_{0})+q\frac{\nu \phi_{B}(s\bar{c}_{1})+s\phi_{B}(s+\nu c_{1})}{\nu+s}\right)-\frac{spd_{n+1}}{\delta-s},
        \end{array}
    \end{displaymath}}}
    where $d_{n}:=P(W_{n}=0)$ and used the fact: 
    \begin{displaymath}
        \begin{array}{rl}
       E(e^{-sV_{n}W_{n}})=&\int_{0}^{1}E(e^{-svW_{n}})dv=\frac{1}{s}\int_{0}^{s}E(e^{-yW_{n}})dy\vspace{2mm}\\
       =&   \frac{1}{s}\int_{0}^{s}Z_{n}(y)dy.  
        \end{array}
    \end{displaymath}
 If $W_{0}=w$, then $E(e^{-sW_{0}})=e^{-sw_{0}}$, the last expression allows to recursively determine all the transforms $Z_{n}(s)$, $n\in \mathbb{N}$.
    
 Multiplying with $\delta-s$, and setting $s=\delta$:
\begin{displaymath}
    p_{n+1}=\frac{\phi_{B}(\delta\bar{c}_{0})}{\delta}\int_{0}^{\delta}Z_{n}(y)dy.
\end{displaymath}

Let $U_{W}(r,s):=\sum_{n=0}^{\infty}r^{n}Z_{n}(s)$, $|r|<1$, then:
\begin{equation}
    U_{W}(r,s)=r\frac{\Psi(s)}{s(\delta-s)}\int_{0}^{s}U_{W}(r,y)dy+K(s),\label{hgmf}
\end{equation}
where
\begin{displaymath}
\begin{array}{rl}
  \Psi(s)=   & \frac{p\delta\phi_{B}(s\bar{c}_{0})+q(\delta-s)\frac{\nu \phi_{B}(s\bar{c}_{1})+s\phi_{B}(s+\nu c_{1})}{\nu+s}}{\delta-s}\vspace{2mm} \\
    K(s)=&e^{-sw_{0}}-\frac{sp}{\delta-s}(U_{W}(r,\infty)-p_{0}).\end{array}
\end{displaymath}
Letting $I(s)=\int_{0}^{s}U_{W}(r,y)dy$, \eqref{hgmf} becomes:
{\small{\begin{displaymath}
    I^{\prime}(s)=r\frac{\Psi(s)}{s(\delta-s)}I(s)+K(s).
\end{displaymath}}}
The solution of such kind of first-order differential equation is obtained by following the lines in \cite[Section 5]{box2}, so further details are omitted.
\subsection{Randomly proportional dependency with additive delay}
In the following we consider the case where $A_{n}=G_{n}B_{n}+J_{n}$, with $P(G_{n}=\beta_{i})=p_{i}$, $i=1,\ldots,K$, and $J_{n}$ to follow a hyperexponential distribution with pdf $f(x)=\sum_{j=1}^{L}q_{j}\delta_{j}e^{-\delta_{j}x}$ (the analysis can be further generalized in the case of a distribution with a rational Laplace transform). Then,
\begin{displaymath}\begin{array}{l}
     Z_{n+1}(s)=E(e^{-sW_{n+1}})=E(e^{-s[V_{n}W_{n}+(1-G_{n})B_{n}-J_{n}]^{+}})  \vspace{2mm}\\
     =\sum_{j=1}^{L}q_{j}\sum_{l=1}^{K}p_{l}\int_{v=0}^{1}\int_{w=0}^{\infty}\int_{x=0}^{\infty}f_{B}(x)\left[\int_{y=0}^{vw+\bar{\beta}_{l}x}e^{-s(vw+\bar{\beta}_{l}x-y)}\delta_{j}e^{-\delta_{j}y}+\int_{y=vw+\bar{\beta}_{l}x}^{\infty}\delta_{j}e^{-\delta_{j}y}dy\right]dxdP(W<w)dv\vspace{2mm}\\
     =\sum_{j=1}^{L}q_{j}\sum_{l=1}^{K}p_{l}\int_{v=0}^{1}\int_{w=0}^{\infty}\int_{x=0}^{\infty}f_{B}(x)\left[\frac{\delta_{j}e^{-s(vw+\bar{\beta}_{l}x)}-se^{-\delta_{j}(vw+\bar{\beta}_{l}x)}}{\delta_{j}-s}\right]dxdP(W<w)dv\vspace{2mm}\\
     =\sum_{j=1}^{L}q_{j}(\frac{\delta_{j}}{\delta_{j}-s})\sum_{l=1}^{K}p_{l}\phi_{B}(s\bar{\beta}_{l})\frac{1}{s}\int_{0}^{s}Z_{n}(y)dy-s\sum_{j=1}^{L}(\frac{q_{j}}{\delta_{j}-s})\sum_{l=1}^{K}p_{l}\phi_{B}(\delta_{j}\bar{\beta}_{l})\frac{1}{\delta_{j}}\int_{0}^{\delta_{j}}Z_{n}(y)dy\vspace{2mm}\\
     =\frac{\sum_{j=1}^{L}q_{j}\delta_{j}\prod_{m\neq j}(\delta_{m}-s)\sum_{l=1}^{K}p_{l}\phi_{B}(s\bar{\beta}_{l})}{s\prod_{m=1}^{L}(\delta_{m}-s)}\int_{0}^{s}Z_{n}(y)dy-s\sum_{j=1}^{L}\frac{q_{j}}{\delta_{j}-s}c_{j,n+1},
\end{array}
\end{displaymath}
where 
\begin{displaymath}
    c_{j,n+1}:=\frac{\sum_{l=1}^{K}p_{l}\phi_{B}(\delta_{j}\bar{\beta}_{l})}{\delta_{j}}\int_{0}^{\delta_{j}}Z_{n}(y)dy=P(W_{n+1}=0|Q=j),\,j=1,\ldots,L,
\end{displaymath}
where $Q$ the type of the arrival process. Then, multiplying with $r^{n+1}$ and sum over $n$ (with $W_{0}=w$) results in
\begin{displaymath}
    U_{W}(r,s)=r\frac{N(s)}{sD(s)}\sum_{l=1}^{K}p_{l}\phi_{B}(s\bar{\beta}_{l})\int_{0}^{s}U_{W}(r,y)dy+e^{-sw}-s\sum_{j=1}^{M}\frac{q_{j}}{\delta_{j}-s}[U^{(j)}_{W}(r,\infty)-c_{j,0}],
    \end{displaymath}
    where
    $N(s):=\sum_{j=1}^{L}q_{j}\delta_{j}\prod_{m\neq j}(\delta_{m}-s)$, $D(s):=\prod_{j=1}^{M}(\delta_{j}-s)$, and $U^{(j)}_{W}(r,s)=E(e^{-sW}|Q=j)$. Note that $U^{(j)}_{W}(r,s)$, $j=1,\ldots,L$. Then,
\begin{equation}
I^{\prime}(s)=r\frac{N(s)}{sD(s)}\sum_{l=1}^{K}p_{l}\phi_{B}(s\bar{\beta}_{l})I(s)+K(r,s),\label{bbb}
\end{equation}
where
\begin{displaymath}
    K(r,s):=e^{-sw}-s\sum_{j=1}^{M}\frac{q_{j}}{\delta_{j}-s}[U^{(j)}_{W}(r,\infty)-c_{j,0}].
\end{displaymath}
The form of \eqref{bbb} is that same as in \cite[Section 3]{box3}, and the analysis can be performed similarly.
\subsection{Proportionally dependent on system time}
We now consider the case where $A_{n}=c(W_{n}+B_{n})+J_{n}$, $c\in(0,1)$. We assume that $(B_{n},J_{n})$ are i.i.d. sequences of random vectors. Thus, the quantities $(\bar{c}B_{n}-J_{n})$ are i.i.d. random variables, however, within a pair $B_{n}$, $J_{n}$ are dependent. Here, we assume that a non-negative random vector $(B,J)$ has a bivariate matrix-exponential distribution with LST $E(e^{-sB-zJ}):=\frac{G(s,z)}{D(s,z)}$, where $G(s,z)$, $D(s,z)$ are polynomial functions in $s,z$. A consequence of this definition is that the transform of $Y:=\bar{c}B-J$ is also a rational function; the distribution of $Y$ is called a bilateral matrix-exponential \cite[Theorem 3.1]{bla}. This class of distributions, under which we model the dependence structure, belongs to the class of multivariate matrix-exponential distributions, which was introduced in \cite{blaa}. For ease of notation, let $E(e^{-sY}):=h(s)=\frac{f(s)}{g(s)}$, and assume that $g(s)$ has $L$ zeros, say $t_{j}$ such that $Re(t_{j})>0$, $j=1,\ldots,L$, and $M$ zeros, say $h_{j}$, such that $Re(h_{j})<0$, $j=1,\ldots,M$, whereas $f(s)$ is a polynomial of degree at most $M+L-1$, not sharing the same zeros with $g(s)$.

Note that now recursion \eqref{ddd} becomes
\begin{displaymath}
    W_{n+1}=[(V_{n}-c)W_{n}+\bar{c}B_{n}-J_{n}]^{+},
\end{displaymath}
so that $V_{n}-c\sim U(-c,\bar{c})$. Then, for $H_{n}=[(V_{n}-c)W_{n}+\bar{c}B_{n}-J_{n}]^{-}$, $Re(s)=0$
\begin{displaymath}
    \begin{array}{rl}
        E(e^{-sW_{n+1}})= &\frac{f(s)}{g(s)}[\int_{-c}^{0}E(e^{-svW_{n}})dv+\int_{0}^{\bar{c}}E(e^{-svW_{n}})dv]+1-E(e^{-sH_{n}}).
    \end{array}
\end{displaymath}
Multiplying with $r^{n+1}$ ($0<r<1$) and summing for $n=0$ to infinity, and having in mind that $W_{0}=w$, we obtain
\begin{equation}
    g(s)(Z_{w}(r,s)-e^{-sw})-rf(s)\int_{0}^{\bar{c}} Z_{w}(r,sy_{1})dy_{1}=      rf(s)\int_{-c}^{0} Z(r,sy_{1})dy_{1}+rg(s)(\frac{1}{1-r}-H(r,s)),\label{bnl}
\end{equation}
where $H(r,s)=\sum_{n=0}^{\infty}r^{n}E(e^{-sH_{n}})$. Now we have:
\begin{enumerate}
    \item the LHS of \eqref{bnl} is analytic in $Re(s)>0$ and continuous in $Re(s)\geq 0$,
    \item the RHS of \eqref{bnl} is analytic in $Re(s)<0$ and continuous in $Re(s)\leq 0$,
    \item for large $s$, both sides are $O(s^{M+L})$ in their respective half planes.
\end{enumerate}
Thus, for $Re(s)\geq 0$,
\begin{equation}
    g(s)(Z_{w}(r,s)-e^{-sw})-rf(s)\int_{0}^{\bar{c}} Z_{w}(r,sy_{1})dy_{1}=\sum_{l=0}^{M+L}C_{l}(r)s^{l},\label{hu1}
\end{equation}
and for $Re(s)\leq 0$,
\begin{equation}
    rf(s)\int_{-c}^{0} Z(r,sy_{1})dy_{1}+rg(s)(\frac{1}{1-r}-H(r,s))=\sum_{l=0}^{M+L}C_{l}(r)s^{l}.\label{hu2}
\end{equation}
For $s=0$, \eqref{hu1} yields
\begin{displaymath}
    C_{0}(r)=g(0)(\frac{1}{1-r}-1)-rf(0)\int_{0}^{\bar{c}}\frac{dy_{1}}{1-r}=\frac{rc}{1-r}g(0),
\end{displaymath}
having in mind that $f(0)=g(0)$. The same value for $C_{0}(r)$ can be derived from \eqref{hu2} by setting $s=0$. We can also obtain other $L$ equations for the remaining coefficients. For $s=t_{j}$, $j=1,\ldots,L$ in \eqref{hu1} we obtain:
\begin{equation}
    -rf(t_{j})\int_{0}^{\bar{c}} Z_{w}(r,t_{j}y_{1})dy_{1}=\sum_{l=0}^{M+L}C_{l}(r)(t_{j})^{l}.\label{hu3}
\end{equation}

% For $s=\delta$,
%\begin{displaymath}
 %   c_{1}=-r\frac{\phi_{B}(\delta\bar{c})}{\delta}\int_{0}^{\bar{c}\delta}Z(r,y)dy.
%\end{displaymath}
Proceeding in a similarly as in \cite{box2,hoo},
\begin{equation}
   Z_{w}(r,s)=r\frac{f(s)}{g(s)}\int_{0}^{\bar{c}}Z_{w}(r,sy_{1})dy_{1}+L(r,s),\, Re(s)\geq 0,\label{poi}
\end{equation}
%\begin{equation}
 %   Z(r,s)=r\frac{\delta\phi_{B}(s\bar{c})}{s(\delta-s)}\int_{0}^{\bar{c}s}Z(r,y)dy+e^{-sw}+\frac{s}{\delta-s}c_{1}.\label{poi}
%\end{equation}
where $L(r,s):=e^{-sw}+\frac{\sum_{l=0}^{M+L}C_{l}(r)s^{l}}{g(s)}$.
Next, we follow the lines in \cite{hoo}. Note that for $r\in[0,1)$, $|K(r,s)|:=|r\frac{f(s)}{g(s)}|\leq r<1$ as $s\to 0$. Iterating \eqref{poi} $n$ times we obtain
\begin{displaymath}
    \begin{array}{rl}
Z_{w}(r,s)=&\int\ldots\int_{[0,\bar{c}]^{n+1}}K(r,s)\prod_{h=1}^{n}K(r,sy_{1}\ldots y_{h})Z_{w}(r,sy_{1}\ldots y_{n+1})dy_{1}\ldots dy_{n+1}\vspace{2mm}\\
&+L(r,s)+\sum_{j=1}^{n}\int\ldots\int_{[0,\bar{c}]^{j}}K(r,s)\prod_{h=1}^{j-1}K(r,sy_{1}\ldots y_{h})L(r,sy_{1}\ldots y_{j})dy_{1}\ldots dy_{j}.
\end{array}
\end{displaymath}
Since we will let $n$ to tend to $\infty$ we are interested in investigating the convergence of the summation in the previous expression as well as in obtaining the limit of the first term in the right hand-side of the previous expression. Since the expressions of $K(r,s)$, $L(r,s)$ share the same properties as those in \cite{hoo}, we can show that
\begin{equation}
    Z_{w}(r,s)=L(r,s)+\sum_{n=1}^{\infty}\int\ldots\int_{[0,\bar{c}]^{n}}K(r,s)\prod_{j=1}^{n-1}K(r,sy_{1}\ldots y_{j})L(r,sy_{1}\ldots y_{n})dy_{1}\ldots dy_{n}.\label{njic}
\end{equation}
We still need $M$ more equations to obtain a system for the coefficients $C_{l}(r)$. Substituting $s=h_{j}$, $j=1,\ldots,L$ in \eqref{hu2} and using \eqref{njic} we obtain
\begin{equation}
    rf(h_{j})\int_{-c}^{0} Z(r,h_{j}y_{1})dy_{1}=\sum_{l=0}^{M+L}C_{l}(r)(h_{j}))^{l}.\label{hu4}
\end{equation}
Finally, by using \eqref{hu3}, \eqref{hu4} we can derive the remaining coefficients $C_{l}(r)$, $l=1.\ldots,L+M$.
\begin{remark}
An alternative way to solve \eqref{poi} is by performing the transformation $v_{1}=sy_{1}$, so that \eqref{poi} becomes:
\begin{equation} 
   Z_{w}(r,s)=r\int_{0}^{\bar{c}s}h(s)Z_{w}(r,v_{1})dv_{1}+L(r,s),\, Re(s)\geq 0.\label{poi1}
\end{equation}
Note that \eqref{poi1} is a Fredholm equation \cite{masu}; therefore, a natural way to proceed is
by successive substitutions. Define now iteratively the function
\begin{displaymath}
    L^{i^*}(r,s):=r\int_{0}^{\bar{c}s}h(s)L^{(i-1)^*}(r,v)dv,\,i\geq 1,
\end{displaymath}
with $L^{0^*}(r,s):=L(r,s)$. Then, after $n$ iterations we
have that
\begin{displaymath}
    Z_{w}(r,s)=\sum_{i=0}^{n+1}L^{i^*}(r,s)+r^{n+1}\int_{v_{1}=0}^{\bar{c}s}\int_{v_{2}=0}^{\bar{c}v_{1}}\ldots\int_{v_{n+1}=0}^{\bar{c}v_{n}}h(s)\prod_{j=1}^{n}h(v_{j})Z_{w}(r,v_{n+1})dv_{n+1}\ldots dv_{2}dv_{1}.
\end{displaymath}
Note that
\begin{displaymath}
    \lim_{n\to\infty}r^{n+1}\int_{v_{1}=0}^{\bar{c}s}\int_{v_{2}=0}^{\bar{c}v_{1}}\ldots\int_{v_{n+1}=0}^{\bar{c}v_{n}}h(s)\prod_{j=1}^{n}h(v_{j})Z_{w}(r,v_{n+1})dv_{n+1}\ldots dv_{2}dv_{1}=0.
\end{displaymath}
To see this, observe that
\begin{displaymath}
    |h(v_{n})\int_{v_{n+1}=0}^{\bar{c}v_{n}}Z_{w}(r,v_{n+1})dv_{n+1}|<|\int_{v_{n+1}=0}^{1}Z_{w}(r,v_{n+1})dv_{n+1}|\leq \frac{1}{1-r}<1.
\end{displaymath}
Thus, the above limit is less than or equal to

\begin{displaymath}
    \lim_{n\to\infty}r^{n+1}\frac{1}{1-r}=0.
\end{displaymath}
Therefore
    \begin{equation}
        Z_{w}(r,s)=\sum_{i=0}^{\infty}L^{i^*}(r,s).\label{plo}
    \end{equation}
Now for $Re(s)\geq 0$ $M_{2}(r,s)=max_{v\in[0,\bar{c}s]}|L(r,s)|<\infty$. Then, 
\begin{displaymath}
    |L^{i^*}(r,s)|<|\int_{0}^{\bar{c}s}L^{(i-1)^*}(r,s)|\leq \bar{c}s max_{v\in[0,\bar{c}s]}|L(r,s)|=\bar{c}sM_{2}(r,s)<\infty,
\end{displaymath}
which ensures the convergence of the infinite sum in \eqref{plo}.
\end{remark}

%Equation \eqref{poi} is similar to a linear Volterra integral equation of the second kind:
%\begin{equation}
 %   I(s)=\delta F(s)\int_{0}^{\psi(s)}I(y)dy+K(s),\label{vca}
%\end{equation}
%where $\psi(s):=\bar{c}s$, $I(s):=Z(r,s)$, $F(s):=r\frac{\phi_{B}(s\bar{c})}{s(\delta-s)}$, $K(s):=e^{-sw}+\frac{s}{\delta-s}c_{1}$. Successive iterations imply that
%\begin{displaymath}
 %   I(s)=\sum_{n=0}^{\infty}\delta^{n}u_{n}(s),
%\end{displaymath}
%where $u_{0}(s):=K(s)$, and for $n=1,2,\ldots,$
%\begin{displaymath}
 %   u_{n}(s):=F(s)\int_{y=0}^{\bar{c}s}u_{n-1}(y)dy
%\end{displaymath}
%The constant $c_{1}$ can be obtained by solving the equation
%\begin{displaymath}
 %   c_{1}=-r\frac{\phi_{B}(\delta\bar{c})}{\delta}\int_{0}^{\bar{c}\delta}\sum_{n=0}^{\infty}u_{n}(y)dy.
%\end{displaymath}
\subsection{A Bernoulli dependent structure}
Consider the following (simpler) case of recursion in \eqref{rec2} where $V_{n}^{(1)}<0$ a.s., and $V_{n}^{(2)}=U_{n}^{1/a}$, with $U_{n}\sim U(0,1)$, $a\geq 2$:
\begin{equation}
    W_{n+1}=\left\{\begin{array}{ll}
          \left[V^{(1)}_{n}W_{n}+B_{n}-A^{(1)}_{n}\right]^{+},&\text{w.p. } p,\\
        \left[U_{n}^{1/a}W_{n}+T_{n}-A^{(2)}_{n}\right]^{+},& \text{w.p. } q:=1-p,
    \end{array}\right.\label{recz2}
\end{equation}
where the LST of $B_{n}$, say $\phi_{B}(s):=\frac{N_{B}(s)}{D_{B}(s)}$ is rational with poles at $s_{1},\ldots,s_{l}$, with $Re(s_{j})<0$, $j=1,\ldots,l$. Then, for $Re(s)=0$,
\begin{displaymath}
    \begin{array}{rl}
       E(e^{-sW_{n+1}})=  &p E(e^{V_{n}W_{n}})\phi_{B}(s)\phi_{A_{1}}(-s)+qE(e^{U^{1/a}_{n}W_{n}})\phi_{T}(s)\phi_{A_{2}}(-s) +1-J_{n}(s),
    \end{array}
\end{displaymath}
where for $n=0,1,\ldots,$
\begin{displaymath}
    J_{n}(s):=pE(e^{\left[V^{(1)}_{n}W_{n}+B_{n}-A^{(1)}_{n}\right]^{-}})+qE(e^{\left[U_{n}^{1/a}W_{n}+T_{n}-A^{(2)}_{n}\right]^{-}}).
\end{displaymath}
Note that for $u=sv^{1/a}$, we have,
\begin{displaymath}
    E(e^{U^{1/a}_{n}W_{n}})=\int_{0}^{1}E(e^{v^{1/a}_{n}W_{n}})dv=\frac{a}{s^{a}}\int_{0}^{s}u^{a-1}E(e^{uW_{n}})du.
\end{displaymath}
Setting $Z_{w}^{(a)}(r,s):=s^{a-1}Z_{w}(r,s)$ and proceed as in \cite{box2} we obtain,
\begin{equation}
\begin{array}{l}
   D_{B}(s)[Z_{w}^{(a)}(r,s)-s^{a-1}e^{-sw}-rq\frac{a}{s}\phi_{T}(s)\phi_{A_{2}}(-s)\int_{0}^{s}Z^{(a)}_{w}(r,u)du]\vspace{2mm}\\
   =rs^{a-1}[pN_{B}(s)\phi_{A_{1}}(-s)\int_{-\infty}^{0}Z_{w}(r,sv)P(V^{(1)}\in dy)+D_{B}(s)(\frac{1}{1-r}-M_{w}(r,s))],
   \end{array}\label{mop}
    \end{equation}
    where $M_{w}(r,s):=\sum_{n=0}^{\infty}r^{n}J_{n}(s)$.
Note that
\begin{itemize}
    \item The LHS in \eqref{mop} is analytic for $Re(s)>0$, and continuous for $Re(s)\geq 0$,
    \item The RHS in \eqref{mop} is analytic for $Re(s)<0$, and continuous for $Re(s)\leq 0$,
    \item For large $s$, both sides are $O(s^{l})$ in their respective half-planes.
\end{itemize}
It follows by Liouville’s theorem \cite[p. 85]{tit} that
\begin{equation}
     D_{B}(s)[Z_{w}^{(a)}(r,s)-s^{a-1}e^{-sw}-rq\frac{a}{s}\phi_{T}(s)\phi_{A_{2}}(-s)\int_{0}^{s}Z^{(a)}_{w}(r,u)du]=    \sum_{k=0}^{l}c_{k}(r)s^{k},\,\,Re(s)\geq 0,\label{pq1}
\end{equation}
  \begin{equation}
      rs^{a-1}[pN_{B}(s)\phi_{A_{1}}(-s)\int_{-\infty}^{0}Z_{w}(r,sv)P(V^{(1)}\in dy)+D_{B}(s)(\frac{1}{1-r}-M_{w}(r,s))]=      \sum_{k=0}^{l}c_{k}(r)s^{k},\,\,Re(s)\leq 0.\label{pq2}
  \end{equation}
Setting $s=0$ in either \eqref{pq1} or \eqref{pq2} yields $c_{0}(r)=0$. For $s=s_{j}$, $D_{B}(s_{j})=0$, $j=1,\ldots, l$. Substituting in \eqref{pq2} yields
\begin{equation}
    rs_{j}^{a-1}[pN_{B}(s_{j})\phi_{A_{1}}(-s_{j})\int_{-\infty}^{0}Z_{w}(r,s_{j}v)P(V^{(1)}\in dy)=      \sum_{k=1}^{l}c_{k}(r)s_{j}^{k}.
\end{equation}
Note that from \eqref{pq1}
\begin{displaymath}
    Z_{w}^{(a)}(r,s)=rq\frac{a}{s}\phi_{T}(s)\phi_{A_{2}}(-s)\int_{0}^{s}Z^{(a)}_{w}(r,u)du+s^{a-1}e^{-sw}+\frac{\sum_{k=1}^{l}c_{k}(r)s^{k}}{D_{B}(s)},
\end{displaymath}
or equivalently, if $I^{(a)}(s):=\int_{0}^{s}Z_{w}^{(a)}(r,u)du$,
\begin{equation}
    I^{(a)\prime}(s)=rq\frac{a}{s}\phi_{T}(s)\phi_{A_{2}}(-s)I^{(a)}(s)+\frac{\sum_{k=1}^{l}c_{k}(r)s^{k}}{D_{B}(s)}+s^{a-1}e^{-sw}.\label{qerw}
\end{equation}

Assume now that $A$ Thus, following standard techniques from the theory of ordinary differential equations, we have for a positive number $c$, such that $c<s$,
\begin{displaymath}
    I^{(a)}(s)=e^{\int_{c}^{s}rq\frac{a}{u}\phi_{T}(u)\phi_{A_{2}}(-u)du}\left(I^{(a)}(c)+\int_{c}^{s}e^{-\int_{c}^{t}rq\frac{a}{u}\phi_{T}(u)\phi_{A_{2}}(-u)du}\left(\frac{\sum_{k=1}^{l}c_{k}(r)t^{k}}{D_{B}(t)}+t^{a-1}e^{-tw}\right)dt\right).
\end{displaymath}
Note that
\begin{displaymath}
  \int_{c}^{s}rq\frac{a}{u}\phi_{T}(u)\phi_{A_{2}}(-u)du=-(1+o(1))rqa\ln(c),\text{ as }c\to 0.  
\end{displaymath}
Since $I^{(a)\prime}(s)=s^{a-1}Z{w}(r,s)$, $I^{(a)\prime}(0)=0$, and thus,
\begin{displaymath}
    I^{(a)}(s)=\int_{0}^{s}e^{\int_{t}^{s}rq\frac{a}{u}\phi_{T}(u)\phi_{A_{2}}(-u)du}\left(\frac{\sum_{k=1}^{l}c_{k}(r)t^{k}}{D_{B}(t)}+t^{a-1}e^{-tw}\right)dt.
\end{displaymath}
Combining the above with \eqref{qerw}, and having in mind that $I^{(a)\prime}(s)=Z_{w}^{(a)}(r,s)=s^{a-1}Z_{w}(r,s)$ we have that
\begin{displaymath}
    Z_{w}(r,s)=\frac{\sum_{k=1}^{l}c_{k}(r)s^{k}}{s^{a-1}D_{B}(s)}+e^{-sw}+rq\frac{a}{s^{a}}\phi_{T}(s)\phi_{A_{2}}(-s)\int_{0}^{s}e^{\int_{t}^{s}rq\frac{a}{u}\phi_{T}(u)\phi_{A_{2}}(-u)du}\left(\frac{\sum_{k=1}^{l}c_{k}(r)t^{k}}{D_{B}(t)}+t^{a-1}e^{-tw}\right)dt.
\end{displaymath}

Now substituting the derived expression for $Z_{w}(r,s)$ in \eqref{pq2} we can obtain a system of equations to obtain the remaining unknowns $c_{k}(r)$, $k=1,\ldots,l$.

\subsection{Other generalizations}
We consider the case where
\begin{displaymath}
    W_{n+1}=[V_{n}W_{n}+B_{n}-A_{n}]^{+},
\end{displaymath}
with $V_{n}\sim U(0,1)$, and $E(e^{-s A_{n}}|B_{n}=t)=\chi(s)e^{-\psi(s)t}$, and $B_{n}\sim exp(\mu)$. Thus, the interarrival times distribution depends on the service time of the previous customer, so that
\begin{displaymath}
    E(e^{-sA_{n}-zB_{n}})=\int_{0}^{\infty}\mu e^{-\mu t}e^{-zt}\chi(s)e^{-\psi(s)t}dt=\frac{\mu\chi(s)}{\mu+\psi(s)+z},
\end{displaymath}
when $Re(\mu+\psi(s)+z)>0$. Since for $s=0$ the $E(e^{-s A_{n}}|B_{n}=t)$ should be equal
to one, we have to implicitly assume that $\psi(0)=0$ and $\chi(0)=1$. Therefore, by denoting $Z_{n}(s)=E(e^{-s W_{n}})$ we have:
\begin{displaymath}
    \begin{array}{rl}
       Z_{n+1}(s):=E(e^{-s W_{n+1}}) = &E(e^{-s (V_{n}W_{n}+B_{n}-A_{n})})+1-E(e^{-s [V_{n}W_{n}+B_{n}-A_{n}]^{-}}) \vspace{2mm} \\
        = & E(e^{-s U_{n}W_{n}})E(e^{-s(B_{n}-A_{n})})+1-U^{-}_{V_{n}W_{n}}(s)\vspace{2mm}\\
        =&E(e^{-s U_{n}W_{n}})\frac{\mu\chi(-s)}{\mu+\psi(-s)+s}+1-U^{-}_{V_{n}W_{n}}(s),
    \end{array}
\end{displaymath}
where $U^{-}_{V_{n}W_{n}}(s):=E(e^{-s [V_{n}W_{n}+B_{n}-A_{n}]^{-}})$. Clearly, under the transformation $v=su$, we have:
\begin{displaymath}
    E(e^{-s U_{n}W_{n}})=\int_{0}^{1}E(e^{-s uW_{n}})du=\frac{1}{s}\int_{0}^{s}Z_{n}(v)dv.
\end{displaymath}
Thus, assuming that $\chi(s):=\frac{P_{1}(s)}{Q_{1}(s)}$, $\psi(s):=\frac{P_{2}(s)}{Q_{2}(s)}$, with $P_{2}(s),Q_{1}(s),Q_{2}(s)$ polynomials of degrees $L$, $M$, and $N$ respectively:
\begin{displaymath}
\begin{array}{rl}
     Z_{n+1}(s)=&\frac{P_{1}(-s)}{sQ_{1}(-s)}\frac{\mu Q_{2}(-s)}{(\mu+s)Q_{2}(-s)+P_{2}(-s)}\int_{0}^{s}Z_{n}(v)dv+1-U^{-}_{V_{n}W_{n}}(s)\vspace{2mm}\\
     =&\frac{\mu N_{Y}(s)}{s D_{Y}(s)}\int_{0}^{s}Z_{n}(v)dv+1-U^{-}_{V_{n}W_{n}}(s).
\end{array}
\end{displaymath}
Multiplying with $r^{n+1}$ and summing for zero to infinity we obtain
\begin{displaymath}
   s D_{Y}(s)[Z_{w}(r,s)-e^{-sw}]-r\mu N_{Y}(s)\int_{0}^{s}Z_{w}(r,v)dv=rs D_{Y}(s)(\frac{1}{1-r}-M_{w}(r,s)),
\end{displaymath}
where $M_{w}(r,s):=\sum_{n=0}^{\infty}r^{n}U^{-}_{V_{n}W_{n}}(s)$. Note that $D_{Y}(s):=Q_{1}(-s)((\mu+s)Q_{2}(-s)+P_{2}(-s))$ is a polynomial of degree $M+N+1$. Thus, following similar arguments and Liouville's theorem we have
\begin{displaymath}
    \begin{array}{rl}
        s D_{Y}(s)[Z_{w}(r,s)-e^{-sw}]-r\mu N_{Y}(s)\int_{0}^{s}Z_{w}(r,v)dv=& \sum_{l=0}^{M+N+L+2}C_{l}(r)s^{l},\,Re(s)\geq 0,  \\
        rs D_{Y}(s)(\frac{1}{1-r}-M_{w}(r,s))=& \sum_{l=0}^{M+N+L+2}C_{l}(r)s^{l},\,Re(s)\leq 0. 
    \end{array}
\end{displaymath}
For $s=0$, we can easily derive $C_{0}(r)=0$. Assuming that all the zeros of $D_{Y}(s)$, say $t_{j}$, $j=1,\ldots,M+N+1$ are all in positive half plane, we can derive a system of equations for the remaining coefficients $C_{j}(r)$:
\begin{displaymath}
    -r\mu N_{Y}(t_{j})\int_{0}^{t_{j}}Z_{w}(r,v)dv=\sum_{l=1}^{M+N+L+2}C_{l}(r)t_{j}^{l}.
\end{displaymath}
Now for $Re(s)\geq 0$, we have,
\begin{displaymath}
    Z_{w}(r,s)=r\mu \frac{N_{Y}(s)}{D_{Y}(s)}\int_{0}^{s}Z_{w}(r,v)dv+e^{-sw}-\frac{\sum_{l=1}^{M+N+L+2}C_{l}(r)s^{l}}{D_{Y}(s)}.
\end{displaymath}
The form of the above equation is the same as in \cite[eq (48), p. 239]{box2}, so it can be solved similarly. 

\section{On modified versions of a multiplicative Lindley recursion with dependences}
In the following, we focus on the recursion \eqref{rec22}, which generalizes the model in \cite{hoo}. More precisely, we assume that $V_{n}^{(1)}$ are such that $P(V_{n}^{(1)}\in [0,1))=1$, and $V_{n}^{(2)}$ such that $P(V_{n}^{(2)}<0)=1$. We further use $\mu$ to denote the probability measure on $[0,1)$, i.e., $\mu(A):=P(V_{n}^{(1)}\in A)$ for every Borel set $A$ on $[0,1)$.

%Assume also that $\{A_{n}^{(k)}\}$, $k=0,1,2$ are independent sequences of i.i.d. random variables with LST $\phi_{A_{k}}(s):=E(e^{-sA_{n}^{(k)}}):=\frac{N_{k}(s)}{D_{k}(s)}$, $D_{k}(s):=\prod_{j=1}^{K_{k}}(s+t_{j}^{(k)})$, with $Re(t_{j}^{(k)})>0$, and $N_{k}(s)$ polynomial of degree at most $K_{k}-1$, not sharing same zeros with $D_{k}(s)$, $k=0,1,2$. Moreover, $\{B_{n}\}$ are i.i.d. random variables with LST $\phi_{B}(s)=\frac{N}{}$

Assume also that $\{Y_{n}^{(0)}:=B_{n}-A_{n}^{(0)}\}_{n=0}^{\infty}$ are i.i.d. random variables and their LST, say $\phi_{Y_{0}}(s):=E(e^{-sY_{n}^{(0)}}):=\frac{N_{0}(s)}{D_{0}(s)}$, with $D_{0}(s):=\prod_{i=1}^{L}(s+m_{i})\prod_{j=1}^{K_{0}}(s+t_{j}^{(0)})$, with $Re(m_{i})<0$, $i=1,\ldots,L$, $Re(t_{j}^{(0)})>0$, $j=1,\ldots,K_{0}$. Assume also that $\{A_{n}^{(k)}\}$, $k=1,2$ are independent sequences of i.i.d. random variables with LST $\phi_{A_{k}}(s):=E(e^{-sA_{n}^{(k)}}):=\frac{N_{k}(s)}{D_{k}(s)}$, $D_{k}(s):=\prod_{j=1}^{K_{k}}(s+t_{j}^{(k)})$, with  $N_{k}(s)$ polynomial of degree at most $K_{k}-1$, not sharing same zeros with $D_{k}(s)$, and $Re(t_{j}^{(1)})>0$, $Re(t_{j}^{(2)})<0$, $j=1,\ldots,K_{k}$, $k=1,2$. Then,
\begin{displaymath}
    \begin{array}{l}
         E(e^{-s W_{n+1}}) = pE(e^{-s [W_{n}+B_{n}-A_{n}^{(0)}]^{+}})+q\left[E(e^{-s [V^{(1)}_{n}W_{n}+\widehat{B}_{n}-A_{n}^{(1)}]}1(\widehat{B}_{n}\leq T_{n}))+E(e^{-s [V^{(2)}_{n}W_{n}+T_{n}-A_{n}^{(2)}]}1(\widehat{B}_{n}> T_{n}))\right.\vspace{2mm}\\\left.+1-E(e^{-s [V^{(1)}_{n}W_{n}+\widehat{B}_{n}-A_{n}^{(1)}]^{-}}1(\widehat{B}_{n}\leq T_{n}))-E(e^{-s [V^{(2)}_{n}W_{n}+T_{n}-A_{n}^{(2)}]^{-}}1(\widehat{B}_{n}> T_{n}))\right] \vspace{2mm} \\
        = pE(e^{-s W_{n}})\frac{N_{0}(s)}{D_{0}(s)}+qE(e^{-sV_{n}^{(1)} W_{n}})\chi(s)\frac{N_{1}(-s)}{D_{1}(-s)}+qE(e^{-sV_{n}^{(2)} W_{n}})\psi(s)\frac{N_{2}(-s)}{D_{2}(-s)}+1-J_{n}^{-}(s),
    \end{array}
\end{displaymath}
where
\begin{displaymath}
\begin{array}{rl}
    J_{n}^{-}(s):=&pE(e^{-s [W_{n}+B_{n}-A_{n}^{(0)}]^{-}})+q[E(e^{-s [V^{(1)}_{n}W_{n}+\widehat{B}_{n}-A_{n}^{(1)}]^{-}}1(\widehat{B}_{n}\leq T_{n}))+E(e^{-s [V^{(2)}_{n}W_{n}+T_{n}-A_{n}^{(2)}]^{-}}1(\widehat{B}_{n}> T_{n}))],\vspace{2mm}\\
    \chi(s):= &E(e^{-s\widehat{B}_{n}}1(\widehat{B}_{n}\leq T_{n}))=\int_{0}^{\infty}e^{-sx}(1-T(x))d\widehat{B}(x),  \vspace{2mm}\\
        \psi(s):= &E(e^{-sT_{n}}1(\widehat{B}_{n}> T_{n}))=\int_{0}^{\infty}e^{-sx}(1-\widehat{B}(x))dT(x).
    \end{array}
\end{displaymath}
Letting $Z_{w}(r,s)=\sum_{n=0}^{\infty}r^{n}E(e^{-s W_{n}})$, $r\in[0,1)$ we have for $Re(s)=0$ that
\begin{equation}
    \begin{array}{l}
        D_{1}(-s)D_{2}(-s)\left[Z_{w}(r,s)(D_{0}(s)-rpN_{0}(s))-D_{0}(s)e^{-sw}\right]-rq\chi(s)D_{0}(s)D_{2}(-s)N_{1}(-s)\int_{[0,1)}Z_{w}(r,sy)P(V^{(1)}\in dy)\vspace{2mm}\\
        =D_{0}(s)[rq\psi(s)N_{2}(-s)\int_{(-\infty,0)}Z_{w}(r,sy)P(V^{(2)}\in dy)+rD_{1}(-s)D_{2}(-s)(\frac{1}{1-r}-J^{-}(r,s))],\label{mop1}
    \end{array}
\end{equation}
where $J^{-}(r,s)=\sum_{n=0}^{\infty}r^{n}J_{n}^{-}(s)$. It is readily seen that
\begin{itemize}
    \item The LHS in \eqref{mop1} is analytic for $Re(s)>0$, and continuous for $Re(s)\geq 0$,
    \item The RHS in \eqref{mop1} is analytic for $Re(s)<0$, and continuous for $Re(s)\leq 0$,
    \item For large $s$, both sides are $O(s^{L+K_{0}+K_{1}+K_{2}})$ in their respective half-planes.
\end{itemize}
Thus, Liouville implies that
\begin{displaymath}
    \begin{array}{l}
     D_{1}(-s)D_{2}(-s)\left[Z_{w}(r,s)(D_{0}(s)-rpN_{0}(s))-D_{0}(s)e^{-sw}\right]\vspace{2mm}\\-rq\chi(s)D_{0}(s)D_{2}(-s)N_{1}(-s)\int_{[0,1)}Z_{w}(r,sy)P(V^{(1)}\in dy)=    \sum_{l=0}^{L+K_{0}+K_{1}+K_{2}}C_{l}(r)s^{l},\,Re(s)\geq 0,\vspace{2mm}\\
     D_{0}(s)D_{1}(-s)[rq\psi(s)N_{2}(-s)\int_{(-\infty,0)}Z_{w}(r,sy)P(V^{(2)}\in dy)+rD_{2}(-s)(\frac{1}{1-r}-J^{-}(r,s))]\vspace{2mm}\\=         \sum_{l=0}^{L+K_{0}+K_{1}+K_{2}}C_{l}(r)s^{l},\,Re(s)\leq 0, 
    \end{array}
\end{displaymath}
where $C_{l}(r)$, $l=0,1,\ldots,L+K_{0}+K_{1}+K_{2}$, constants to be derived. For $s=0$, simple computations imply that
\begin{displaymath}
    C_{0}(r)=\frac{r}{1-r}(1-p-q\chi(0))\prod_{j=1}^{L}m_{j}\prod_{j=1}^{K_{0}}t_{j}^{(0)}\prod_{j=1}^{K_{1}}t_{j}^{(1)}\prod_{j=1}^{K_{2}}t_{j}^{(2)}.
\end{displaymath}
Thus, for $Re(s)\geq 0$, we have
\begin{equation}
    Z_{w}(r,s)=K(r,s)\int_{[0,1)}Z_{w}(r,sy_{1})\mu(dy_{1})+L(r,s),\label{kah}
\end{equation}
where
\begin{displaymath}
    K(r,s):=\frac{rq\chi(s)D_{0}(s)N_{1}(-s)}{D_{1}(-s)(D_{0}(s)-rpN_{0}(s))},\,\,\,L(r,s):=\frac{D_{0}(s)e^{-sw}+\sum_{l=0}^{L+K_{0}+K_{1}+K_{2}}C_{l}(r)s^{l}}{D_{0}(s)-rpN_{0}(s)}.
\end{displaymath}
\eqref{kah} has the same form as in \cite[eq. (13)]{hoo} and can be treated similarly. Note that in our case, for $r\in[0,1)$,
\begin{displaymath}
    |K(r,s)|\leq\frac{rq|\chi(s)||D_{0}(s)||N_{1}(-s)|}{|D_{1}(-s)|(|D_{0}(s)|-rp|N_{0}(s)|)}\to\frac{rq\chi(0)}{1-rp}<\frac{rq}{1-rp}\leq r<1,
\end{displaymath}
as $s\to 0$. Thus, there is a positive constant $\epsilon$ such that for $s$ satisfying $|s|\leq \epsilon$, we have $|K(r,s)|\leq \bar{r}:=\frac{1+r}{2}$. Note that $K(r,s)$, $L(r,s)$ satisfy the same properties as those in \cite{hoo}, thus, proceeding similarly and iterating $n$ times \eqref{kah} we obtain
\begin{equation}\begin{array}{rl}
     Z_{w}(r,s)=&L(r,s)+\sum_{j=1}^{n}\int\ldots\int_{[0,1)^{j}}K(r,s)\prod_{h=1}^{j-1}K(r,sy_{1}\ldots y_{h})L(r,sy_{1}\ldots y_{j})\mu(dy_{1})\ldots\mu(dy_{j})\\
     &+\int\ldots\int_{[0,1)^{n+1}}K(r,s)\prod_{h=1}^{n}K(r,sy_{1}\ldots y_{h})Z(r,sy_{1}\ldots y_{n+1})\mu(dy_{1})\ldots\mu(dy_{n+1})
     \end{array}
     \label{kah1}
\end{equation}
We will let $n\to\infty$ to obtain $Z_{w}(r,s)$, so we need to verify the convergence of the summation in the second term in \eqref{kah1}, as well as to obtain the limit of the second third term in \eqref{kah1}. Following the lines in \cite[pp. 9-10]{hoo} we can finally obtain,
\begin{equation}
    Z_{w}(r,s)=L(r,s)+\sum_{n=1}^{\infty}\int\ldots\int_{[0,1)^{n}}K(r,s)\prod_{j=1}^{n-1}K(r,sy_{1}\ldots y_{j})L(r,sy_{1}\ldots y_{n})\mu(dy_{1})\ldots\mu(dy_{n}).\label{nji}
    \end{equation}
We still need to derive the remaining coefficients $C_{l}(r)$, $l=1,\ldots,L+\sum_{k=0}^{3}K_{k}$: First, by using Rouch\'e's theorem, we can show that $D_{0}(s)-rpN_{0}(s)=0$ has $K_{0}$ roots, say $\delta_{1}(r),\ldots,\delta_{K_{0}}(r)$, with $Re(\delta_{j}(r))\geq 0$, $j=1,\ldots,K_{0}$. Thus, we can obtain $K_{0}$ equations:
\begin{displaymath}\begin{array}{l}
    -rq\chi(\delta_{j}(r))D_{0}(\delta_{j}(r))D_{2}(-\delta_{j}(r))N_{1}(-\delta_{j}(r))\int_{[0,1)}Z_{w}(r,\delta_{j}(r)y)P(V^{(1)}\in dy)\vspace{2mm}\\= D_{1}(-\delta_{j}(r))D_{2}(-\delta_{j}(r))D_{0}(\delta_{j}(r))e^{-\delta_{j}(r)w}   +\sum_{l=0}^{L+K_{0}+K_{1}+K_{2}}C_{l}(r)(\delta_{j}(r))^{l}.\end{array}
\end{displaymath}
Similarly, for $s=t_{j}^{(1)}$, $j=1,\ldots,K_{1}$,
\begin{displaymath}\begin{array}{l}
    -rq\chi(t_{j}^{(1)})D_{0}(t_{j}^{(1)})D_{2}(-t_{j}^{(1)})N_{1}(-t_{j}^{(1)})\int_{[0,1)}Z_{w}(r,t_{j}^{(1)}y)P(V^{(1)}\in dy)\vspace{2mm}\\= \sum_{l=0}^{L+K_{0}+K_{1}+K_{2}}C_{l}(r)(t_{j}^{(1)})^{l}.\end{array}
\end{displaymath}
For $s=t_{j}^{(2)}$, $j=1,\ldots,K_{2}$,
\begin{displaymath}\begin{array}{l}
    -rq\psi(t_{j}^{(2)})D_{0}(t_{j}^{(2)})D_{1}(-t_{j}^{(2)})N_{2}(-t_{j}^{(2)})\int_{(-\infty,0)}Z_{w}(r,t_{j}^{(2)}y)P(V^{(2)}\in dy)\vspace{2mm}\\= \sum_{l=0}^{L+K_{0}+K_{1}+K_{2}}C_{l}(r)(t_{j}^{(2)})^{l},\end{array}
\end{displaymath}
while for $s=m_{j}$, $j=1,\ldots,L$,
\begin{displaymath}\begin{array}{l}
\sum_{l=0}^{L+K_{0}+K_{1}+K_{2}}C_{l}(r)(m_{j})^{l}=0.\end{array}
\end{displaymath}
By inserting where is needed the expression in \eqref{nji}, we obtain a system of $L+K_{0}+K_{1}+K_{2}$ equations to obtain $C_{l}(r)$, $l=1,\ldots,L+K_{0}+K_{1}+K_{2}$.
\subsection{A mixed-autoregressive case}
Consider first a simple version of the recursion \eqref{rec23}, i.e., $W_{n+1}=[V_{n}W_{n}+B_{n}-A_{n}]^{+}$, where now $P(V_{n}=a)=p_{1}$, $P(V_{n}\in[0,1))=p_{2}$, and $P(V_{n}<0)=1-p_{1}-p_{2}$, with $a\in(0,1)$, $0\leq p_{1}\leq 1$, $0\leq p_{2}\leq 1$, $p_{1}+p_{2}\leq 1$ (the general version of \eqref{rec23} will be considered in Remark \ref{rem11} below). Note that the case where $a=1$ was analyzed in \cite{hoo}. In the following, we fill the gap in the literature, by analysing the case where $a\in(0,1)$, which we call it as \textit{mixed-autoregressive}, in the sense that in the obtained functional equation we will have the terms: $Z_{w}(r,as)$, and $\int_{[0,1)}Z(r,sy)P(V\in dy)$. Assume that $V^{+}\stackrel{\rm def}{=}(V|V\in[0,1))$, $V^{-}\stackrel{\rm def}{=}(V|V<0)$. Then, for $Re(s)=0$, $r\in[0,1)$ we have
\begin{equation}
\begin{array}{l}
Z_{w}(r,s)-e^{-sw}=rp_{1}\phi_{Y}(s)Z_{w}(r,as)+rp_{2}\phi_{Y}(s)\int_{[0,1)}Z_{w}(r,sy)P(V^{+}\in dy)\vspace{2mm}\\
+r(1-p_{1}-p_{2})\phi_{Y}(s)\int_{(-\infty,0)}Z_{w}(r,sy)P(V^{-}\in dy)+r(\frac{1}{1-r}-J^{-}(r,s)),\end{array}
    \label{eqoo}
\end{equation}
where $\{Y_{n}=B_{n}-A_{n}\}_{n=0}^{\infty}$ i.i.d. random variables with LST $\phi_{Y}(s):=\frac{N_{Y}(s)}{D_{Y}(s)}$, with $D_{Y}(s):=\prod_{j=1}^{L}(s-t_{j})\prod_{j=1}^{M}(s-s_{j})$. Without loss of generality we assume that $Re(t_{j})>0$, $j=1,\ldots,L$, $Re(s_{j})<0$, $j=1,\ldots,M$. Thus, \eqref{eqoo} becomes
\begin{equation}
\begin{array}{l}
D_{Y}(s)(Z_{w}(r,s)-e^{-sw})-rp_{1}N_{Y}(s)Z_{w}(r,as)-rp_{2}N_{Y}(s)\int_{[0,1)}Z_{w}(r,sy)P(V^{+}\in dy)\vspace{2mm}\\
=r(1-p_{1}-p_{2})N(s)\int_{(-\infty,0)}Z_{w}(r,sy)P(V^{-}\in dy)+rD_{Y}(s)(\frac{1}{1-r}-J^{-}(r,s)).\end{array}
    \label{eqoo1}
\end{equation}
It is readily seen that
\begin{itemize}
    \item The LHS in \eqref{eqoo1} is analytic for $Re(s)>0$, and continuous for $Re(s)\geq 0$,
    \item The RHS in \eqref{eqoo1} is analytic for $Re(s)<0$, and continuous for $Re(s)\leq 0$,
    \item For large $s$, both sides are $O(s^{L+M})$ in their respective half-planes.
\end{itemize}
Thus, Liouville implies that for $Re(s)\geq 0$,
\begin{equation}
    D_{Y}(s)(Z_{w}(r,s)-e^{-sw})-rp_{1}N_{Y}(s)Z_{w}(r,as)-rp_{2}N_{Y}(s)\int_{[0,1)}Z_{w}(r,sy)P(V^{+}\in dy)=\sum_{j=0}^{M+L}C_{l}(r)s^{j},\label{eq01}
\end{equation}
and for $Re(s)\leq 0$,
\begin{equation}
    r(1-p_{1}-p_{2})N(s)\int_{(-\infty,0)}Z_{w}(r,sy)P(V^{-}\in dy)+rD_{Y}(s)(\frac{1}{1-r}-J^{-}(r,s))=\sum_{j=0}^{M+L}C_{l}(r)s^{j}.\label{eq02}
\end{equation}
By using either \eqref{eq01} or \eqref{eq02} for $s=0$, we obtain
\begin{displaymath}
    C_{0}(r)=\frac{r(1-p_{1}-p_{2})}{1-r}\prod_{j=1}^{L}t_{j}\prod_{j=1}^{M}s_{j}.
\end{displaymath}
Denoting by $\mu$ the probability measure on $[0,1)$ induced by $V^{+}$, \eqref{eq01} is written as
\begin{equation}
    Z_{w}(r,s)=p_{1}K(r,s)Z_{w}(r,as)+p_{2}K(r,s) \int_{[0,1)}Z_{w}(r,sy_{1})\mu(dy_{1})+L_{w}(r,s),\label{eq03}
\end{equation}
where
\begin{displaymath}
    K(r,s):=r\phi_{Y}(s),\,\,L_{w}(r,s):=e^{-sw}+\frac{\sum_{j=0}^{M+L}C_{l}(r)s^{j}}{D_{Y}(s)}.
\end{displaymath}
Our aim is to solve \eqref{eq03}, which combines the model in \cite{box1}, with those in \cite{hoo,box3}, i.e., in the functional equation the unknown function $Z_{w}(r,s)$ arises also as $Z_{w}(r,as)$ as well as in $\int_{[0,1)}Z_{w}(r,sy)\mu(dy)$. Let for $i,j=0,1,\ldots,$
\begin{displaymath}
    \begin{array}{rl}
         f_{i,j}(s):=&a^{i}\prod_{k=1}^{j}y_{k}s,\,y_{k}\in[0,1),k=1,\ldots,j,  \vspace{2mm}\\
         F(r,f_{i,j}(s)):=&\left\{\begin{array}{ll}
              Z_{w}(r,a^{i}s),&j=0, \vspace{2mm} \\
              \int\ldots\int_{[0,1)^{j}}Z_{w}(r,f_{i,j}(s))\mu(dy_{1})\ldots\mu(dy_{j}),&j\geq 1,
         \end{array}\right.
    \end{array}
\end{displaymath}
where $f_{i,0}(s)=a^{i}s$ (i.e., $\prod_{k=1}^{0}y_{k}:=1$). Moreover, $f_{i,j}(f_{k,l}(s))=f_{i+k,j+l}(s)=f_{k,l}(f_{i,j}(s))$. Then, \eqref{eq03} becomes
\begin{equation}
    F(r,s)=p_{1}K(r,s)F(r,f_{1,0}(s))+p_{2}K(r,s)F(r,f_{0,1}(s))+L_{w}(r,s),\label{ere1}
\end{equation}
where $F(r,s)=F(r,f_{0,0}(s))=Z_{w}(r,s)$. Iterating \eqref{ere1} $n-1$ times yields,
\begin{equation}
    F(r,s)=\sum_{k=0}^{n}p_{1}^{k}p_{2}^{n-k}G_{k,n-k}(s)F(r,f_{k,n-k}(s))+\sum_{k=0}^{n-1}\sum_{m=0}^{k}p_{1}^{m}p_{2}^{m-k}G_{m,k-m}(s)\tilde{L}(r,f_{m,k-m}(s)),\label{ere2}
\end{equation}
where $G_{k,n-k}(r,s)$ are recursively defined as follows (with $G_{-1,.}(s) = G_{.,-1}(s) \equiv 0$, $G_{0,0}(s)=1$):
\begin{displaymath}
    \begin{array}{rl}
         G_{1,0}(s):=&K(r,s),\,\,\,\, G_{0,1}(s):=K(r,s),\vspace{2mm} \\
       G_{k+1,n-k}(s)=  &G_{k,n-k}(s) \tilde{K}(r,f_{k,n-k}(s))+G_{k+1,n-1-k}(s) \tilde{K}(r,f_{k+1,n-1-k}(s)),\vspace{2mm}\\
       G_{k,n+1-k}(s)=  &G_{k-1,n+1-k}(s) \tilde{K}(r,f_{k-1,n+1-k}(s))+G_{k,n-k}(s) \tilde{K}(r,f_{k,n-k}(s)),
    \end{array}
\end{displaymath}
where also,
\begin{displaymath}
    \begin{array}{rl}
         \tilde{K}(r,f_{i,j}(s)):=&\left\{\begin{array}{ll}
              K(r,a^{i}s),&j=0, \vspace{2mm} \\
              \int\ldots\int_{[0,1)^{j}}K(r,f_{i,j}(s))\mu(dy_{1})\ldots\mu(dy_{j}),&j\geq 1,
         \end{array}\right.
    \end{array}
\end{displaymath}
and
\begin{displaymath}
    \begin{array}{rl}
         \tilde{L}(r,f_{i,j}(s)):=&\left\{\begin{array}{ll}
              L_{w}(r,a^{i}s),&j=0, \vspace{2mm} \\
              \int\ldots\int_{[0,1)^{j}}L_{w}(r,f_{i,j}(s))\mu(dy_{1})\ldots\mu(dy_{j}),&j\geq 1.
         \end{array}\right.
    \end{array}
\end{displaymath}
It can easily verified that $G_{k,n-k}(r,s)$ is a sum of $\binom{n}{k}$ terms, and each of them is a product of $n$ terms of values of $\tilde{K}(r,f_{.,.})$, which are related to the LST $\phi_{Y}(.)$. We have to mention that our framework is related to the one developed in \cite{adan} with the difference that the functions $f_{i,j}(s)$ (for $j>0$) are more complicated compared to the corresponding $a_{i}(z)$ in \cite{adan}, and inherit difficulties in solving \eqref{ere1}. 

In what follows we will let $n\to\infty$ in \eqref{ere2} so that to obtain an expression for $F(r,s)$. In doing that, we have to verify the convergence of the summation in the second term in the right hand side of \eqref{ere2}, as well as to estimate the limit of the corresponding first term in the right hand side of \eqref{ere2}. The key ingredient shall be the boundness of $G_{k,n-k}(s)$. Similarly to \cite[p.8]{adan}, $G_{k,n-k}(s)$ can be interpreted as the total weight of of all $\binom{n}{k}$ paths from $(0,0)$ to $(k,n-k)$. Let $C_{k,n-k}$ the set of all paths leading from $(0,0)$ to $(k,n-k)$, where a path from $(0,0)$ to $(k,n-k)$ is defined as a sequence of grid points starting from $(0,0)$ and ending to $(k,n-k)$ by only taking unit steps $(1,0)$, $(0,1)$. Then, a typical term (one of the $\binom{n}{k}$ terms) of $G_{k,n-k}(s)$ should be the following:
\begin{displaymath}
    \int\ldots\int_{[0,1)^{m}}\prod_{(l,m)\in C_{k,n-k}}\tilde{K}(r,a^{l}y_{1}\ldots y_{m}s)\mu(dy_{1})\ldots\mu(dy_{n-k}),
\end{displaymath}
for $m=0,\ldots, n-k$, and $l=0,\ldots,k$ with $(l,m)\neq (k,n-k)$. For $Re(s)\geq 0$, $M_{1}(r,s):=sup_{y\in[0,1]}|K(r,sy)|<\infty$, $M_{2}(r,s):=sup_{y\in[0,1]}|L(r,sy)|<\infty$, and $|K(r,s)|\leq r<1$. Then, for $a\in(0,1)$, $M_{l}(r,a^{i}s)<M_{l}(r,s)$, $i\geq 1$, $l=1,2$. Following \cite{hoo}
\begin{displaymath}
    |\int\ldots\int_{[0,1)^{m}}\prod_{(l,m)\in C_{k,n-k}}\tilde{K}(r,a^{l}y_{1}\ldots y_{m}s)\mu(dy_{1})\ldots\mu(dy_{n-k})|\leq E\left[\prod_{(l,m)\in C_{k,n-k}}|\tilde{K}(r,a^{l}Z_{1}\ldots Z_{m}s)|\right],
\end{displaymath}
where $Z_{1},Z_{2},\ldots$ is a sequence of i.i.d. random variables with the same distribution as $V^{+}$. Following the same procedure as in \cite[pp. 8-9]{hoo} we can show that the each of the weights of the path is bounded, implying that $G_{k,n-k}(s)$ is also bounded. This result will imply as $n\to\infty$, that the first term in the right hand side of \eqref{ere2} vanishes, Thus,
\begin{equation}
F(r,s)=\sum_{k=0}^{\infty}\sum_{m=0}^{k}p_{1}^{m}p_{2}^{m-k}G_{m,k-m}(s)\tilde{L}(r,f_{m,k-m}(s)).\label{hpl}
\end{equation}

We are now ready to obtain the coefficients $C_{l}(r)$, $l=1,\ldots,M+L$. For $s=t_{j}$, $j=1,\ldots,L$, in \eqref{eq01}, we have
\begin{equation}
    -rp_{1}N_{Y}(t_{j})Z_{w}(r,at_{j})-rp_{2}N_{Y}(t_{j})\int_{[0,1)}Z_{w}(r,t_{j}y)\mu(dy)=\sum_{l=0}^{M+L}C_{l}(r)t_{j}^{l}.\label{eq11}
\end{equation}
For $s=s_{j}$, $j=1,\ldots,M$,
\begin{equation}
    r(1-p_{1}-p_{2})N(s_{j})\int_{(-\infty,0)}Z_{w}(r,s_{j}y)P(V^{-}\in dy)=\sum_{l=0}^{M+L}C_{l}(r)s_{j}^{l}.\label{eq12}
\end{equation}
\eqref{eq11}, \eqref{eq12} serve as a system of equations for the coefficients $C_{l}(r)$, $l=1,\ldots,M+L$.
\begin{remark}\label{rem11} We now turn back in the general case of recursion \eqref{rec23}. The analysis is still applicable when we assume that with probability $p_{1}$, $\{V_{n}^{(1)}\}_{n\in\mathbb{N}_{0}}$ follows a discrete random variable with finite number of values, i.e., with probability $p_{1}$, $V_{n}\in\{a_{1},\ldots,a_{M}\}$, with $a_{k}\in(0,1)$, $K=1,\ldots,M$, and $P(V_{n}=a_{k})=q_{k}$. Then, \eqref{eq03} takes the following form
    \begin{equation}
    Z_{w}(r,s)=p_{1}K(r,s)\sum_{k=1}^{M}q_{k}Z_{w}(r,a_{k}s)+p_{2}K(r,s) \int_{[0,1)}Z_{w}(r,sy_{1})\mu(dy_{1})+L_{w}(r,s),\label{eq031}
\end{equation}
Then, by setting $h_{j}:=p_{1}q_{j}$, $j=1,\ldots,M$, $h_{M+1}:=p_{2}$, $f_{i_{1},\ldots,i_{M},i_{M+1}}(s):=a_{1}^{i_{1}}\ldots a_{M}^{i_{M}}\prod_{j=1}^{i_{M+1}}y_{j}s$, and $e_{j}^{(M+1)}$ an $1\times (M+1)$ row vector with 1 at the $j$th position and all the other entries equal to zero, \eqref{ere1} becomes
\begin{equation}
    F(r,s)=K(r,s)\sum_{j=1}^{M+1}h_{j}F(r,f_{e_{j}^{(M+1)}}(s))+L_{w}(r,s).\label{ert}
\end{equation}
\eqref{ert} has the same form as the functional equations treated in \cite[eq. (2)]{adan}. After $n$ iterations \eqref{ert} becomes
\begin{displaymath}
\begin{array}{l}
      F(r,s)=\sum_{i_{1}+\ldots+i_{M}+i_{M+1}=n+1}h_{1}^{i_{1}}\ldots h_{M}^{i_{M}}h_{M+1}^{i_{M+1}}G_{i_{1},\ldots,i_{M},i_{M+1}}(s)F(r,f_{i_{1},\ldots,i_{M},i_{M+1}}(s))\vspace{2mm}\\
      +\sum_{k=0}^{n}\sum_{i_{1}+\ldots+i_{M}+i_{M+1}=k}h_{1}^{i_{1}}\ldots h_{M}^{i_{M}}h_{M+1}^{i_{M+1}}G_{i_{1},\ldots,i_{M},i_{M+1}}(s)\tilde{L}(r,f_{i_{1},\ldots,i_{M},i_{M+1}}(s)),
      
      \end{array}
\end{displaymath}
where now
\begin{displaymath}\begin{array}{rl}
    G_{i_{1},\ldots,i_{M},i_{M+1}}(s)=&\sum_{j=1}^{M+1}\tilde{K}(r,f_{i_{1},\ldots,i_{j-1},\ldots,i_{M+1}})G_{i_{1},\ldots,i_{j-1},\ldots,i_{M+1}}(s),\vspace{2mm}\\
    \tilde{K}(r,f_{i_{1},\ldots,i_{M+1}}(s)):=&\left\{\begin{array}{ll}
              K(r,a_{1}^{i_{1}}\ldots a_{M}^{i_{M}}s),&i_{M+1}=0, \vspace{2mm} \\
              \int\ldots\int_{[0,1)^{j}}K(r,f_{i_{1},\ldots,i_{M+1}}(s))\mu(dy_{1})\ldots\mu(dy_{i_{M+1}}),&i_{M+1}\geq 1,
         \end{array}\right.
         \vspace{2mm}\\
    \tilde{L}(r,f_{i_{1},\ldots,i_{M+1}}(s)):=&\left\{\begin{array}{ll}
              L(r,a_{1}^{i_{1}}\ldots a_{M}^{i_{M}}s),&i_{M+1}=0, \vspace{2mm} \\
              \int\ldots\int_{[0,1)^{j}}L(r,f_{i_{1},\ldots,i_{M+1}}(s))\mu(dy_{1})\ldots\mu(dy_{i_{M+1}}),&i_{M+1}\geq 1,
         \end{array}\right.
    \end{array}
\end{displaymath}
with $G_{0,\ldots,0,0}(s):=1$, $G_{i_{1},\ldots,i_{M},i_{M+1}}(s)=0$, in case one of the indices becomes $-1$. Following the approach above and having in mind that the function $f_{i_{1},\ldots,i_{M+1}}(s)$ are cummutative contraction mappings on $\{s\in\mathbb{C};Re(s)\geq 0\}$, $F(r,s):=Z_{w}(r,s)$ can be derived by using \cite[Theorem 3]{adan}.
\end{remark}
\begin{remark}
    Note that in this subsection we have not considered any dependence framework among $B_{n}$, $A_{n}$, since our major focus was on introducing this \textit{mixed-autoregressive} concept, and generalized the work in \cite{hoo}, by assuming $a\in(0,1)$, instead of $a=1$. However, the analysis is still applicable even when we lift this assumption. For example, assume the simple scenario where now with probability $p_{1}$, $A_{n}=cB_{n}+J_{n}$, i.e., the interarrival time is linearly dependent on the service time, with $c\in(0,1)$, $J_{n}\sim exp(\delta)$. Then, \eqref{eq03} becomes
    \begin{equation}
    Z_{w}(r,s)=p_{1}K_{1}(r,s)Z_{w}(r,as)+p_{2}K(r,s) \int_{[0,1)}Z_{w}(r,sy_{1})\mu(dy_{1})+L_{w}(r,s),\label{eq103}
\end{equation}
where now $K_{1}(r,s):=r\frac{\delta}{\delta-s}\phi_{B}(\bar{c}s)$. The rest of analysis can be treated similarly as above. Clearly, the analysis is still applicable if we consider $J_{n}$ to has distribution with rational LST, or more general dependence structure, e.g., $A_{n}=G_{n}(W_{n}+B_{n})+J_{n}$, $P(G_{n}=\beta_{i})=q_{i}$, $i=1,\ldots,M$, or the (random) threshold dependence structure analysed in Section \ref{thre}. Clearly, we can also apply the same steps when lifting independence assumptions for the general case analysed in Remark \ref{rem11}.
\end{remark}
\section{A more general dependence framework}
        In the following we consider a more general dependence structure among $\{B_{n}\}_{n\in\mathbb{N}_{0}}$, $\{A_{n}\}_{n\in\mathbb{N}_{0}}$. More precisely, assume that 
\begin{equation}
    E(e^{-s A_{n}}|B_{n}=t)=\chi(s)\sum_{i=1}^{N}p_{i}e^{-\psi_{i}(s)t},
    \label{vbh}
\end{equation}
thus, the interarrival times distribution depends on the service time of the previous customer, so that
{\small{\begin{displaymath}
    E(e^{-sA_{n}-zB_{n}})=\int_{0}^{\infty}e^{-zt}\chi(s)\sum_{i=1}^{N}p_{i}e^{-\psi_{i}(s)t}dF_{B}(t)=\chi(s)\sum_{i=1}^{N}p_{i}\phi_{B}(z+\psi_{i}(s)),
\end{displaymath}}}
with $Re(\psi_{i}(s)+z)>0$. Clearly $\chi(0)=1$, $\psi_{i}(0)=0$. The component $e^{-\psi_{i}(s)t}$ depends on the previous service time. The component $\chi(s)$ does not depend on the service time. 

Note that with the above framework we can recover some of the cases analysed above. In particular, the case $A_{n}=c B_{n}+J_{n}$: $N=1$, so that $p_{1}=1$,
    \begin{displaymath}\begin{array}{rl}
        E(e^{-s A_{n}}|B_{n}=t)=& E(e^{-s (c B_{n}+J_{n})}|B_{n}=t)=E(e^{-sJ_{n}})e^{-cst},
    \end{array}
    \end{displaymath}
    with $\chi(s):=E(e^{-sJ_{n}})$, $\psi(s)=cs$.

        The case $A_{n}=G_{n}B_{n}+J_{n}$, with $P(G_{n}=\beta_{i})=p_{i}$, $i=1,\ldots,N$. Then:
       \begin{displaymath}
            \chi(s)=E(e^{-sJ_{n}}),\,\,\psi_{i}(s)=\beta_{i}s,\,i=1,\ldots,N.
        \end{displaymath}
        
        Another interesting scenario: Given $B=t$, $A=\sum_{i=1}^{N_{i}(t)}H_{i,k}$, with probability $p_{i}$, $N_{i}(t)\sim Poisson(\gamma_{i}t)$, and $\{H_{i,k}\}$ sequences of i.i.d. r.v. with a rational LST, each of them distributed like $H_{i}$. Then,
        \begin{displaymath}\begin{array}{rl}
        E(e^{-s A_{n}}|B_{n}=t)=& \sum_{i=1}^{N}p_{i}E(e^{-s \sum_{k=1}^{N_{i}(t)}H_{i,k}}|B_{n}=t)=\vspace{2mm}\\
        =&\sum_{i=1}^{N}p_{i}\sum_{l_{i}=0}^{\infty}E(e^{-s \sum_{k=1}^{l_{i}}H_{i,k}}|B_{n}=t)\frac{e^{-\gamma_{i}t}(\gamma_{i}t)^{l_{i}}}{l_{i}!}\vspace{2mm}\\
        =&\sum_{i=1}^{N}p_{i}e^{-\gamma_{i}(1-E(e^{-sH_{i}}))},
    \end{array}
    \end{displaymath}
    so $\chi(s)=1$, $\psi_{i}(s)=\gamma_{i}(1-E(e^{-sH_{i}}))$. 
    
    So, turning back to the simpler general scenario for the stochastic recursion in \eqref{recu1}: $W_{n+1}=[aW_{n}+B_{n}-A_{n}]^{+}$, where the interarrival times distribution depends on the service time of the previous customer based on \eqref{vbh}, we have:
        \begin{displaymath}
            \begin{array}{rl}
                E(e^{-sW_{n+1}}) =&E(e^{-s(aW_{n}+B_{n}-A_{n})})+1-E(e^{-s[aW_{n}+B_{n}-A_{n}]^{-}}) \vspace{2mm}\\ 
                =&E(e^{-saW_{n}})E(e^{-s(B_{n}-A_{n})})+1-U_{n}(s)\vspace{2mm}\\
                = & E(e^{-saW_{n}})\chi(-s)\sum_{i=1}^{N}p_{i}\phi_{B}(s+\psi_{i}(-s))+1-U_{n}(s).
            \end{array}
        \end{displaymath}
        Assuming that the limit as $n\to\infty$ exists:
        {\small{\begin{displaymath}\begin{array}{c}
            Z(s)=Z(as)\chi(-s)\sum_{i=1}^{N}p_{i}\phi_{B}(s+\psi_{i}(-s))+1-U(s).\end{array}
        \end{displaymath}}}
    Assume that $\chi(s):=\frac{A_{1}(s)}{\prod_{i=1}^{K}(s+\lambda_{i})}$, $\psi_{i}(s):=\frac{B_{i}(s)}{\prod_{j=1}^{L_{i}}(s+\nu_{j})}$, with $A_{1}(s)$ a polynomial of order at most $K-1$, not sharing the same zeros with the denominator of $\chi(s)$, and similarly, $B_{i}(s)$ polynomial of order at most $L_{i}-1$, not sharing the same zeros with the denominator of $\psi_{i}(s)$, for $i=1,\ldots, N$. Then, for $Re(s)=0$,
\begin{displaymath}\begin{array}{c}
            \prod_{i=1}^{K}(\lambda_{i}-s)Z(s)-A_{1}(-s)Z(as)\sum_{i=1}^{N}p_{i}\phi_{B}(s+\psi_{i}(-s))=1-U(s).\end{array}
        \end{displaymath}
        By using similar arguments as above, Liouville implies that 
        \begin{equation}
            \prod_{j=1}^{K}(\lambda_{i}-s)Z(s)-A_{1}(-s)Z(as)\sum_{i=1}^{N}p_{i}\phi_{B}(s+\psi_{i}(-s))=\sum_{j=0}^{K}C_{j}s^{j},\,Re(s)\geq 0.\label{v1}
        \end{equation}
        Setting $s=0$, yields $C_{0}=0$. The rest $C_{j}$s are found by using the $K$ zeros $s=\lambda_{k}$, $k=1,\ldots,K$. Indeed, set $s=\lambda_{k}$, $k=1,\ldots,K$ in \eqref{v1} to obtain the following system:
        \begin{equation}
            -A_{1}(-\lambda_{k})Z(a\lambda_{k})\sum_{i=1}^{N}p_{i}\phi_{B}(\lambda_{k}+\psi_{i}(-\lambda_{k}))=\sum_{j=1}^{K}C_{j}\lambda_{k}^{j}.\label{v2}
        \end{equation}
        However, we still need to find $Z(a\lambda_{k})$, $k=1,\ldots,K$. This can be done by iterating
      \begin{displaymath}\begin{array}{c}
           Z(s)=Z(as)A(-s)\sum_{i=1}^{N}p_{i}\phi_{B}(s+\psi_{i}(-s))+\frac{\sum_{j=1}^{K}C_{j}s^{j}}{\prod_{j=1}^{K}(\lambda_{i}-s)},\end{array}
       \end{displaymath}
       This will result in expressions containing infinite products of the form $\prod_{m=0}^{\infty}A(-a^{m}s)\phi_{B}(a^{m}s+\psi_{i}(-a^{m}s))$. Indeed, after the iterations we get:
       \begin{equation}
           \begin{array}{l}
            Z(s)=\sum_{j=1}^{K}\sum_{n=0}^{\infty}\frac{C_{j}(a^{n}s)^{j}}{\prod_{j=1}^{K}(\lambda_{i}-a^{n}s)}\prod_{m=0}^{n-1}A(-a^{m}s)\sum_{i=1}^{N}p_{i}\phi_{B}(a^{m}s+\psi_{i}(-a^{m}s))\vspace{2mm}\\+\prod_{m=0}^{\infty}A(-a^{m}s)\sum_{i=1}^{N}p_{i}\phi_{B}(a^{m}s+\psi_{i}(-a^{m}s)).\end{array}\label{soly1}
       \end{equation}
       Note that for large $m$, $\phi_{B}(a^{m}s+\psi_{i}(-a^{m}s))$ approaches 1, since $a^{m}s+\psi_{i}(-a^{m}s)\to 0$.
       
       Substituting $s=a\lambda_{k}$, $k=1,\ldots,K$ in \eqref{soly1}, we obtain $Z(a\lambda_{k})$. Finally, by substituting the derived expression in \eqref{v2}, we get a system of equations for the constants $C_{j}$, $j=1,\ldots,K$.
\begin{remark}
    Note that in the independent case, i.e., $\psi(s)=0$, the situation is easy. In the linear dependent case, i.e., $A_{n}=\beta_{i}B_{n}+J_{n}$, $\psi_{i}(s)=\beta_{i}s$, the analysis is also easy to handle. If we additionally assume that $J_{n}\sim exp(\delta)$, then, we are interesting in the convergence of $\prod_{m=0}^{\infty}\frac{\delta\phi_{B}(a^{m}\bar{\beta}_{i}s)}{\delta-a^{m}s}$, which is also easy to handle.
\end{remark}
%       We need further conditions for the general case: Maybe the fact that $\psi_{i}(s)$ is rational functions will help.
      % -- Now we need to show that 
       %\begin{displaymath}
        %   g_{i}(s)=s+\psi_{i}(-s),
       %\end{displaymath}
       %is a contraction on $\{s\in\mathbb{C}:Re(s)\geq 0\}$.\\
       %-- Assuming that $\psi_{i}(s)$ has a completely monotone derivative will help.
        \subsection{Dependence on system time}
        Going one step further, we assume now that 
        \begin{displaymath}
            E(e^{-s A_{n}}|W_{n}+B_{n}=t)=\chi(s)e^{-\psi(s)t},
        \end{displaymath}
        so the interarrival times depends on the system time (sojourn time) of the previous customer. Thus,
        \begin{displaymath}
            E(e^{-sA_{n}-z(W_{n}+B_{n})})=\chi(s)\phi_{B}(z+\psi(s))Z(z+\psi(s)),
        \end{displaymath}
            with $Re(z+\psi(s))>0$. Then, for $Re(s)=0$, the functional equation becomes
       \begin{displaymath}
                Z(s)-\chi(-s)\phi_{B}(s+\psi(-s))Z(s+\psi(-s))=1-U(s).
            \end{displaymath}
            
            The case where $A_{n}=c(W_{n}+B_{n})+J_{n}$, $c\in(0,1)$, $J_{n}\sim exp(\lambda)$ was recently treated in \cite[Section 2]{boxman}. For that case $\chi(s)=\frac{\lambda}{\lambda+s}$, $\psi(s)=sc$.
        
          A more interesting case arise when we assume that the next interarrival time randomly depends on the sojourn time of the previous customer. More precisely,
          \begin{equation}
              E(e^{-s A_{n}}|W_{n}+B_{n}=t)=\chi(s)\sum_{i=1}^{N}p_{i}e^{-\psi_{i}(s)t}.\label{cs1}
          \end{equation}
    
         In such a case, 
           \begin{displaymath}
            E(e^{-sA_{n}-z(W_{n}+B_{n})})=\chi(s)\sum_{i=1}^{N}p_{i}\phi_{B}(z+\psi_{i}(s))Z(z+\psi_{i}(s)),
        \end{displaymath}
            with $Re(z+\psi_{i}(s))>0$, $i=1,\ldots,N$. Then, for $Re(s)\geq 0$, we have
         \begin{equation}
             \prod_{i=1}^{K}(\lambda_{i}-s) Z(s)-A_{1}(-s)\sum_{i=1}^{N}p_{i}\phi_{B}(s+\psi_{i}(-s))Z(s+\psi_{i}(-s))=\sum_{i=0}^{K}C_{i}s^{i}.\label{onb}
         \end{equation}
            
            A special case of of the dependence relation \eqref{cs1} arises when $A_{n}=G_{n}(W_{n}+B_{n})+J_{n}$, $P(G_{n}=\beta_{i})=p_{i}$, $i=1,\ldots,N$; see subsection \ref{cc1}. For such a case, $\chi(s)=\frac{\delta}{\delta+s}$, $\psi_{i}(s)=\beta_{i}s$, $i=1,\ldots,N$. In general, if $g_{i}(s)=s+\psi_{i}(-s)$ is a contraction, then following the lines in \cite{adan}, the functional equation \eqref{onb} can be handled.
            
\section{An integer-valued reflected autoregressive process and a novel retrial queueing system with dependencies}\label{integ}
In this section, we consider the following integer-valued stochastic process $\{X_{n};n=0,1,\ldots\}$ that is determined by the recursion \eqref{rec3}:
\begin{equation}
    X_{n+1}=\left\{\begin{array}{ll}
         \sum_{k=1}^{X_{n}}U_{k,n}+Z_{n}-Q_{n+1},&X_{n}>0 \vspace{2mm} \\
         Y_{n}-\tilde{Q}_{n+1},&X_{n}=0, 
    \end{array}\right.\label{ret1}
\end{equation}
with $Z_{1},Z_{2},\ldots$, $Y_{1},Y_{2},\ldots$, i.i.d. non-negative integer-valued random variables with probability generating function (pgf) $C(z)$, and $G(z)$ respectively, $U_{k,n}$ are i.i.d. Bernoulli distributed random variables with $P(U_{k,n}=1)=\xi_{n}$, $P(U_{k,n}=0)=1-\xi_{n}$ and $Q_{n}$, $\tilde{Q}_{n}$ are i.i.d. random variables such that
\begin{displaymath}\begin{array}{rl}
      P(Q_{n}=0|\sum_{k=1}^{X_{n}}U_{k,n}+Z_{n}=l,X_{n}>0):=&  \frac{\lambda_{1}}{\lambda_{1}+\alpha_{1}(1-\delta_{0,l})}, \vspace{2mm}\\
     P(Q_{n}=1|\sum_{k=1}^{X_{n}}U_{k,n}+Z_{n}=l,X_{n}>0):=&  \frac{\alpha_{1}(1-\delta_{0,l})}{\lambda_{1}+\alpha_{1}(1-\delta_{0,l})},
\end{array}
\end{displaymath}
\begin{displaymath}\begin{array}{rl}
      P(\tilde{Q}_{n}=0|Y_{n}=l,X_{n}=0):=&  \frac{\lambda_{0}}{\lambda_{0}+\alpha_{0}(1-\delta_{0,l})}, \vspace{2mm}\\
     P(\tilde{Q}_{n}=1|Y_{n}=l,X_{n}=0):=&  \frac{\alpha_{0}(1-\delta_{0,l})}{\lambda_{0}+\alpha_{0}(1-\delta_{0,l})}.
\end{array}
\end{displaymath}
Moreover, it is assumed that $\xi_{n}$ are also i.i.d. random variables with $P(\xi_{n}=a_{i})=p_{i}$, $i=1,\ldots,M$, with $\sum_{i=1}^{M}p_{i}=1$. As usual it is assumed that for all $n$, $Z_{n}$, $Y_{n}$, $U_{k,n}$, $Q_{n}$, $\tilde{Q}_{n}$ are independent of each other and of all preceding $X_{r}$.

Note that \eqref{ret1} can be interpreted as follows: Let $X_{n}$ be the number of waiting customers in an orbit queue just after the beginning of the $n$th service, $Q_{n+1}$ (resp. $\tilde{Q}_{n+1}$) be the number of orbiting customers that initiate the $(n+1)$th service when $X_{n}>0$ (resp. $X_{n}=0$). Note that when $X_{n}>0$ (resp. $X_{n}=0$), the first primary customer arrives according to a Poisson process with rate $\lambda_{1}$ (resp. $\lambda_{0}$). $Z_{n}$ (resp. $Y_{n}$) denotes the number of arriving customers during the $n$th service when $X_{n}>0$ with pgf $E(z^{Z_{n}}):=C(z)$ (resp. with pgf $E(z^{Y_{n}}):=G(z)$ $X_{n}=0$). Letting $U_{k,n}$ equal to 1, it is assumed that during the $n$th service each of the $X_{n}$ orbiting customers becomes impatient with probability $1-a$ and leaves without service then, $\{X_{n}\}_{n\in\mathbb{N}_{0}}$ satisfies \eqref{ret1} with $p_{1}=1$, $a_{1}=a$ with $\xi_{n}\equiv a$. Under such a setting the service time and/or the rate of the Poisson arriving process of the number of customers that join the orbit queue during a service time depend on the orbit size at the beginning of the service, as well as the retrieving times (exponentially distributed with rate $\alpha_{0}$ (resp. $\alpha_{1}$) when $X_{n}=0$ (resp. $X_{n}>0$)) depend on the whether the orbit queue is empty or not at the beginning of the last service. We have to note that to our best knowledge it is the first time that such a retrial model is considered in the related literature.

Then,
\begin{displaymath}
    \begin{array}{rl}
         E(z^{X_{n+1}})=&E(z^{\sum_{k=1}^{X_{n}}U_{k,n}+Z_{n}-Q_{n+1}}1(X_{n}>0))+E(z^{Y_{n}-\tilde{Q}_{n+1}}1(X_{n}=0))  \vspace{2mm}\\
         =&E(z^{Z_{n}})(\frac{\alpha_{1}}{z(\lambda_{1}+\alpha_{1})}+\frac{\lambda_{1}}{\lambda_{1}+\alpha_{1}})E(z^{\sum_{k=1}^{X_{n}}U_{k,n}}(1-1(X_{n}=0))) \vspace{2mm}\\&+E(z^{Y_{n}-\tilde{Q}_{n+1}}(1(X_{n}=0,Y_{n}>0)+1(X_{n}=0,Y_{n}>0))) \vspace{2mm}\\
         =&E(z^{Z_{n}})\frac{\alpha_{1}+z\lambda_{1}}{z(\lambda+\alpha_{1})}[E(z^{\sum_{k=1}^{X_{n}}U_{k,n}})-E(1_{(X_{n}=0)})]\vspace{2mm}\\
         &+E(z^{Y_{n}}1(Y_{n}>0))E(1(X_{n}=0))[\frac{\alpha_{0}}{z(\lambda_{0}+\alpha_{0})}+\frac{\lambda_{0}}{\lambda_{0}+\alpha_{0}}]+E(1(Y_{n}=0))E(1(X_{n}=0))\vspace{2mm}\\
         =&E(z^{Z_{n}})\frac{\alpha_{1}+z\lambda_{1}}{z(\lambda+\alpha_{1})}[E(z^{\sum_{k=1}^{X_{n}}U_{k,n}})-E(1(X_{n}=0))]\vspace{2mm}\\
         &+E(1(X_{n}=0))[\frac{\alpha_{0}+z\lambda_{0}}{z(\lambda_{0}+\alpha_{0})}E(z^{Y_{n}})+\frac{a_{0}(z-1)}{z(\lambda_{0}+\alpha_{0})}E(1(Y_{n}=0))],
    \end{array}
\end{displaymath}
where in the third equality we used the fact that when $Y_{n}=0$, then $\tilde{Q}_{n+1}=0$ with certainty. Let $f(z)$ be the pgf of the steady-state distribution of $\{X_{n}\}_{n\in\mathbb{N}_{0}}$ we have after some algebra,
\begin{equation}
    f(z)=\frac{\widehat{C}(z)}{z}\sum_{i=1}^{M}p_{i}f(\bar{a}_{i}+a_{i}z)+\frac{f(0)}{z}[G(0)\frac{\alpha_{0}(z-1)}{\alpha_{0}+\lambda_{0}}+\widehat{G}(z)-\widehat{C}(z)],\label{ret2}
\end{equation}
where $\widehat{C}(z)=C(z)\frac{\alpha_{1}+\lambda_{1}z}{\lambda_{1}+\alpha_{1}}$, $\widehat{G}(z)=G(z)\frac{\alpha_{0}+\lambda_{0}z}{\lambda_{0}+\alpha_{0}}$. After multiplying \eqref{ret2} with $z$ and letting $z=0$, we obtain
\begin{displaymath}
    f(0)=C(0)\sum_{i=1}^{M}p_{i}f(\bar{a}_{i}).
\end{displaymath}

Set $g(z)=\frac{\widehat{C}(z)}{z}$, $K(z)=\frac{f(0)}{z}[G(0)\frac{\alpha_{0}(z-1)}{\alpha_{0}+\lambda_{0}}+\widehat{G}(z)-\widehat{C}(z)]$, so that \eqref{ret2} is now written as
\begin{displaymath}
    f(z)=g(z)\sum_{i=1}^{M}p_{i}f(\bar{a}_{i}+a_{i}z)+K(z),
\end{displaymath}
which has the same form as the one in \cite[Section 5]{adan}. Note that $g(1)=1$, $K(1)=0$, thus, the functional equation in \eqref{ret2} can be solved following \cite[Theorem 6]{adan}, and further details are omitted.

\begin{remark}
Note that $\widehat{C}(z)$ (resp. $\widehat{G}(z)$) refers to the pgf of the number of primary customers that arrive
between successive service initiations when $X_{n}>0$ (resp. $X_{n}=0$). Moreover, we can further assume class dependent service times, i.e., when an orbiting (resp. primary) customer is the one that occupies the server, the pgf of the number of arriving customers during his/her service time equals $C_{o}(z)$ (resp. $C_{p}(z)$). In such a case $\widehat{C}(z)=\frac{\alpha_{1}C_{o}(z)+\lambda_{1}zC_{p}(z)}{\lambda_{1}+\alpha_{1}}$. Similarly, $\widehat{G}(z)=\frac{\alpha_{0}G_{o}(z)+\lambda_{0}zG_{p}(z)}{\lambda_{0}+\alpha_{0}}$ when $X_{n}=0$. 
\end{remark}
\begin{remark}
Moreover, some very interesting special cases may be deduced from \eqref{ret2}. In particular, when $\alpha_{k}\to\infty$, $k=0,1$, then $\widehat{C}(z)=C(z)$, and $\widehat{G}(z)=G(z)$, since $\frac{\alpha_{k}+\lambda_{k}z}{\lambda_{k}+\alpha_{k}}\to 1$ as $\alpha_{k}\to\infty$. Thus, \eqref{ret2} reduces to the functional equation that corresponds to the standard M/G/1 queue generalization in \cite[Section 5]{adan}. Moreover, one can further assume that one of $\alpha_{k}$s to tend to infinity, e.g., $\alpha_{0}\to\infty$ and $a_{1}>0$. In such a scenario, the server has the flexibility to treat the orbit queue as a typical queue, when at the beginning of the last service the orbit queue was empty.
\end{remark}

\subsection{An extension to a two-dimensional case: A priority retrial queue}\label{pri}
In the following, we go one step further towards a multidimensional case. In particular, we consider the two-dimensional discrete-time process $\{(X_{1,n},X_{2,n});n=0,1,\ldots\}$, and assume that only the component $\{X_{2,n};n=0,1,\ldots\}$ is subject to the autoregressive concept, i.e., we generalize the previous model to incorporate two classes of customers (primary and orbiting customers) and priorities, where orbiting customers are impatient. 

Primary customers arrive according to a Poisson process $\lambda_{1}$ and if they find the server busy form a queue waiting to be served. Retrial customers arrive according to a Poisson process $\lambda_{2}$, and upon finding a busy server join an infinite capacity orbit queue, from where they retry according to the constant retrial policy, i.e., only the first in orbit queue attempts to connect with the server after an exponentially distributed time with rate $\alpha$.

Let $X_{i,n}$ be the number of customers in queue $i$ (i.e., type $i$ customers) just after the beginning of the $n$th service, where with $i=1$ ($i=2$) we refer to the orbit (resp. primary) queue. As usual, the server becomes available to the orbiting customers only when there are no customers at the primary queue upon a service completion. We further assume that orbiting customers become impatient during the service of an orbiting customer, according to the machinery described above. 

Let also $A_{i,n}$ be the number of customers of type $i$ that join the system during the $n$th service, $i=1,2$ with pgf $A(z_{1},z_{2})$, and set $\lambda:=\lambda_{1}+\lambda_{2}$. Then $X_{n}:=\{(X_{1,n},X_{2,n});n=0,1,\ldots\}$ satisfies the following recursions:
\begin{displaymath}
    \left\{ \begin{array}{rl}
        X_{1,n+1}=&X_{1,n}+A_{1,n}-1,\,X_{1,n}>0,\,A_{1,n}\geq 0  \\
        X_{2,n+1}=&X_{1,n}+A_{1,n},\,X_{2,n}\geq 0,\,A_{2,n}\geq 0,
    \end{array}\right.
\end{displaymath}
\begin{displaymath}
  \left\{  \begin{array}{rl}
         X_{1,n+1}=&A_{1,n}-1,\,X_{1,n}>0,\,A_{1,n}>0,  \\
        X_{2,n+1}=&X_{1,n}+A_{1,n},\,X_{2,n}\geq 0,
    \end{array}\right.
\end{displaymath}
\begin{displaymath}
    \left\{\begin{array}{rlr}
         X_{1,n+1}=&0,\,X_{1,n}=A_{1,n}=0,&   \\
         &&\text{ with probability }\frac{\lambda}{\lambda+\alpha},\\
         X_{2,n+1}=& X_{2,n}+A_{2,n},\,X_{2,n}>0,\,A_{2,n}\geq 0,&\vspace{2mm}\\
         X_{1,n+1}=&0,\,X_{1,n}=A_{1,n}=0,&   \\
         &&\text{ with probability }\frac{\alpha}{\lambda+\alpha}.\\
         X_{2,n+1}=& \sum_{k=1}^{X_{2,n}}Y_{k,n}+A_{2,n}-1,\,X_{2,n}>0,\,A_{2,n}\geq 0,
    \end{array}\right.
\end{displaymath}
More precisely, the value of the impatience probability equals $\bar{a}_{i}:=1-a_{i}$ with probability $p_{i}$, $i=1,\ldots,M$, i.e., $P(\xi_{n}=a_{i})=p_{i}$, and $P(Y_{k,n}=1)=\xi_{n}$, $P(Y_{k,n}=0)=1-\xi_{n}$. 

\begin{displaymath}
    \left\{\begin{array}{rlr}
         X_{1,n+1}=&0,\,X_{1,n}=A_{1,n}=0,&   \\
         &&\text{ with probability }\frac{\lambda}{\lambda+\alpha},\\
         X_{2,n+1}=& X_{2,n}+A_{2,n},\,X_{2,n}=0,\,A_{2,n}> 0,&\vspace{2mm}\\
         X_{1,n+1}=&0,\,X_{1,n}=A_{1,n}=0,&   \\
         &&\text{ with probability }\frac{\alpha}{\lambda+\alpha},\\
         X_{2,n+1}=& A_{2,n}-1,\,X_{2,n}=0,\,A_{2,n}> 0,
    \end{array}\right.
\end{displaymath}
\begin{displaymath}
    \left\{\begin{array}{rl}
         X_{1,n+1}=&0,\,X_{1,n}=A_{1,n}=0,  \\
         X_{2,n+1}=& 0,\,X_{2,n}=A_{2,n}= 0.
    \end{array}\right.
\end{displaymath}
To our best knowledge it is the first time that such a priority retrial model is considered in the related literature.

Let $F(z_{1},z_{2}):=E(z_{1}^{X_{1,n}}z_{2}^{X_{2,n}})$. Then, using the recursions above, and after lengthy but straightforward calculations we come up with the following functional equation:
\begin{equation}
\begin{array}{rl}
    F(z_{1},z_{2})[z_{1}-A(z_{1},z_{2})]=&\frac{\alpha A(0,z_{2})z_{1}}{z_{2}(\lambda+\alpha)}\sum_{i=1}^{M}p_{i}F(0,\bar{a}_{i}+a_{i}z_{2})\vspace{2mm}\\&-\frac{F(0,z_{2})A(0,z_{2})(\alpha+\lambda(1-z_{1}))}{\lambda+\alpha}+\frac{F(0,0)A(0,0)\alpha(z_{2}-1)z_{1}}{z_{2}(\lambda+\alpha)}.\end{array}\label{plm}
\end{equation}
Then, it is readily seen by using Rouch\'e's theorem \cite[Theorem 3.42, p. 116]{tit} that $z_{1}-A(z_{1},z_{2})$ has for fixed $|z_{2}|\leq 1$, exactly one zero, say $z_{1}=q(z_{2})$ in $|z_{1}|<1$. Substitute $z_{1}=q(z_{2})$ in \eqref{plm} to obtain:
\begin{displaymath}
    F(0,z_{2})\frac{A(0,z_{2})(\alpha+\lambda(1-q(z_{2})))}{\lambda+\alpha}=\frac{\alpha A(0,z_{2})q(z_{2})}{z_{2}(\lambda+\alpha)}\sum_{i=1}^{M}p_{i}F(0,\bar{a}_{i}+a_{i}z_{2})+\frac{F(0,0)A(0,0)\alpha(z_{2}-1)q(z_{2})}{z_{2}(\lambda+\alpha)},
\end{displaymath}
or equivalently, by setting $\tilde{F}(z_{2}):=F(0,z_{2})$, $g(z_{2}):=\frac{\alpha q(z_{2})}{z_{2}(\alpha+\lambda(1-q(z_{2})))}$, $l(z_{2}):=\frac{A(0,0)\alpha(z_{2}-1)q(z_{2})}{A(0,z_{2})(\alpha+\lambda(1-q(z_{2})))z_{2}}$,
\begin{equation}
    \tilde{F}(z_{2})=g(z_{2})\sum_{i=1}^{M}p_{i}\tilde{F}(\bar{a}_{i}+a_{i}z_{2})+l(z_{2}).\label{plm1}
\end{equation}
Note that \eqref{plm1} has the same form as the one in \cite[Section 5, p. 19]{adan}, and $g(1)=1$, $l(1)=0$. Thus, from \cite[Theorem 2]{adan} we can solve \eqref{plm1} an get an expression for $F(0,z_{2})$. Using that expression in \eqref{plm}, we can finally get $F(z_{1},z_{2})$. Note also that from \eqref{plm1}, for $z_{2}=0$, 
\begin{displaymath}
    F(0,0)=\sum_{i=1}^{M}p_{i}F(0,\bar{a}_{i}).
\end{displaymath}
By substituting $z_{2}=\bar{a}_{i}$, $i=1,\ldots,M$, in the derived expression for $F(0,z_{2})$, we finally get $F(0,0)$. Then, by setting $\bar{a}_{i}+a_{i}z_{2}$ instead of $z_{2}$, in the expression for $F(0,z_{2})$, the function $F(z_{1},z_{2})$ is derived through \eqref{plm}.
\section{Conclusion}
In this work we investigated the transient and/or the stationary behaviour of various reflected autoregressive processes. These type of processes are described by stochastic recursions where various independence assumptions among the sequences of random variables that are involved there, are
lifted and for which, a detailed exact analysis can be also provided. This is accomplished by using Liouville's theorem \cite[Theorem 10.52]{tit}, as well as stating and solving a Wiener-Hopf boundary value problem or an integral equation. Various options for follow-up research arise. One of them is to concern multivariate extensions of the processes that we introduced. Such vector-valued counterparts are anticipated to be highly challenging. In subsection \ref{pri} we cope with a simple two-dimensional case, however, the autoregressive parameter was used only in one component. Other possible line of research concerns scaling limits and asymptotics. One also anticipates that, under certain appropriate scalings, a diffusion analysis similar to the one presented in \cite{box1} can be applied.
  
\section*{Acknowledgement} I.D. gratefully acknowledges the hospitality of SMACS group, Dept. of Telecommunication and Information Processing, Ghent University, where a part of this work was carried out.

\section*{Declarations}
\paragraph{Conflict of interest} The authors have no relevant financial or non-financial interests to disclose.
\paragraph{Ethical approval} Not applicable.
\paragraph{Consent for publication} Not applicable.
\paragraph{Consent to participate} Not applicable.
\bibliographystyle{abbrv}

\bibliography{mybibfile}

 \end{document}